\documentclass{compositio}
\usepackage{amsmath,amsfonts,amsthm,mathrsfs,amscd,amssymb,pb-diagram}

\title{Logarithmic nonabelian Hodge theory in characteristic~$p$}
\author{Daniel Schepler}
\email{schepler@math.berkeley.edu}
\address{3121 Maginn Dr.\\Beavercreek, OH 45434}
\newcommand\calA{\mathcal{A}}
\newcommand\calB{\mathcal{B}}
\newcommand\calC{\mathcal{C}}
\newcommand\calD{\mathcal{D}}
\newcommand\calE{\mathcal{E}}
\newcommand\calG{\mathcal{G}}
\newcommand\calH{\mathcal{H}}
\newcommand\calI{\mathcal{I}}
\newcommand\calJ{\mathcal{J}}
\newcommand\calK{\mathcal{K}}
\newcommand\calL{\mathcal{L}}
\newcommand\calR{\mathcal{R}}

\newcommand\scrE{\mathscr{E}}
\newcommand\scrF{\mathscr{F}}

\newcommand\scrM{\mathscr{M}}
\newcommand\scrN{\mathscr{N}}
\newcommand\scrO{\mathscr{O}}
\newcommand\scrT{\mathscr{T}}
\newcommand\bbA{\mathbb{A}}
\newcommand\bbC{\mathbb{C}}

\newcommand\bfT{\mathbf{T}}

\newcommand\frakZ{\mathfrak{Z}}
\newcommand\Jbar{\bar J}
\newcommand\nats{\mathbb{N}}

\newcommand\ints{\mathbb{Z}}
\newcommand\Fp{\mathbb{F}_p}
\newcommand\scrHom{\mathscr{H}om}
\newcommand\scrEnd{\mathscr{E}nd}
\newcommand\Hom{\operatorname{Hom}}
\newcommand\End{\operatorname{End}}
\newcommand\Ext{\operatorname{Ext}}
\newcommand\gp{\operatorname{gp}}
\newcommand\AXgp{\calA_X^{\gp}}
\newcommand\MXgp{\scrM_X^{\gp}}
\newcommand\Mbargp[1]{\overline{\scrM}_{#1}^{\gp}}
\newcommand\MXbargp{\Mbargp{X}}
\newcommand\tildeD{\tilde{\mathcal{D}}}
\newcommand\dlog{\mathop{d\log}}
\newcommand\diffop{differential operator}
\newcommand\PDop{PD differential operator}
\newcommand\id{\operatorname{id}}
\newcommand\Aut{\operatorname{Aut}}
\newcommand\Diff{\mathop{\operatorname{Diff}}}
\newcommand\PDDiff{\mathop{\operatorname{PD~Diff}}}
\newcommand\divpow{\mathcal{D}}
\newcommand\regpow{\mathcal{P}}
\newcommand\DXS{\calD_{X/S}}
\newcommand\tildeDXS{\tildeD_{X/S}}
\newcommand\calXS{\mathcal{X} / \mathcal{S}}
\newcommand\calYS{\mathcal{Y} / \mathcal{S}}
\newcommand\Spec{\operatorname{Spec}}
\newcommand\bfSpec{\operatorname{\mathbf{Spec}}}
\newcommand\Crys{\operatorname{Crys}}
\newcommand\Tor{\operatorname{Tor}}
\newcommand\Gr{\operatorname{Gr}}
\newcommand\FHIG{F\mbox{-}HIG}
\newtheorem{theorem}{Theorem}[section]
\newtheorem{proposition}[theorem]{Proposition}
\newtheorem{lemma}[theorem]{Lemma}
\newtheorem{corollary}[theorem]{Corollary}
\theoremstyle{definition}
\newtheorem{definition}[theorem]{Definition}
\theoremstyle{remark}
\newtheorem{remark}[theorem]{Remark}
\newtheorem{example}[theorem]{Example}
\newenvironment{myequation}{\setcounter{equation}{\arabic{theorem}}\addtocounter{theorem}{1}\begin{equation}}{\end{equation}}

\begin{document}

\begin{abstract}
  Given a morphism $X \to S$ of log schemes of characteristic $p > 0$
  and a lifting $\calXS$ of $X'$ over $S$ modulo $p^2$, we use
  Lorenzon's indexed algebras $\AXgp$ and $\calB_{X/S}$ to construct
  an equivalence between $\scrO_X$-modules with nilpotent integrable
  connection and indexed $\calB_{X/S}$-modules with nilpotent
  $\calB_{X/S}$-linear Higgs field.  If either satisfies a stricter
  nilpotence condition, we find an isomorphism between the de Rham
  cohomology of the connection and the Higgs cohomology of the Higgs
  field.

\end{abstract}

\maketitle
\section{Introduction}

The classical Riemann-Hilbert correspondence gives an equivalence
between the category of coherent modules $E$ with integrable
connection $\nabla$ on a complex manifold $X$ and the category of locally
constant sheaves of complex vector spaces $V$ on $X$.  Moreover, there
is a natural isomorphism
\[ H^n(X, V) \simeq H^n(X, E \otimes \Omega^\cdot_{X / \bbC}), \] where the right hand
side is the hypercohomology of the de Rham complex of $\nabla$.  However,
the maps in this complex are only $\bbC$-linear in general.  On the
other hand, if $V$ is constant and $X$ is the analytic space
associated to a projective scheme over $\bbC$, Hodge theory tells us
that in fact,
\[ H^n(X, V) \simeq \bigoplus_{i + j = n} H^i(X, E \otimes \Omega^j_{X / \bbC}). \] In general, if
$V$ is not constant, this last statement is no longer true.  Simpson
\cite{simpson} resolved this by defining a partial equivalence between
modules with integrable connection $(E, \nabla)$ and modules with Higgs
field $(E', \theta')$ (with a stability condition on both sides), such that
the de Rham complex of $(E, \nabla)$ and the Higgs complex of $(E', \theta')$
have the same hypercohomology.  Recall that a {\em Higgs field} on
$E'$ is a map $\theta' : E' \to E' \otimes_{\scrO_X} \Omega^1_{X/\bbC}$ which, instead of
satisfying the Leibniz rule, is simply {\em $\scrO_X$-linear} and
satisfies the integrability condition $\theta' \land \theta' = 0$, so that $\theta'$
induces a complex of $\scrO_X$-modules
\[ E' \to E' \otimes_{\scrO_X} \Omega^1_{X / \bbC} \to E' \otimes_{\scrO_X} \Omega^2_{X / \bbC}
\to \cdots. \]
(This is equivalent to an extension of the $\scrO_X$-module structure
on $E'$ to an $S^\cdot \scrT_{X/\bbC}$-module structure.)

On schemes in characteristic $p$, the straightforward generalization
of the Riemann-Hilbert correspondence requires that the connection $\nabla$
have vanishing $p$-curvature.  The construction of the $p$-curvature
is as follows: recall that given a morphism $X \to S$ of schemes of
characteristic $p$, for a derivation $D : \scrO_X \to \scrO_X$ over
$\scrO_S$, its $p$th iterate $D^{(p)}$ is again a derivation.  Now
given an $\scrO_X$-module $E$ with integrable connection $\nabla$, it turns
out that
\[ D \mapsto \psi_D := (\nabla_D)^p - \nabla_{D^{(p)}} \in \scrEnd_{\scrO_S} (E) \] induces
a Frobenius-linear map $\scrT_{X/S} \to F_{X/S *} \scrEnd_{\scrO_X}
(E)$, or equivalently an $\scrO_X$-linear map $\psi : E \to E \otimes_{\scrO_X}
F_{X/S}^* \Omega^1_{X'/S}$.  (The appendix to this paper gives a new proof
of this fact which is more conceptual than previous proofs, for
example that in \cite{katz}.)  The classical Cartier descent theorem
\cite{katz} then gives an equivalence between the category of modules
with integrable connection with vanishing $p$-curvature and the
category of $\scrO_{X'}$-modules.

However, many connections of interest do not have vanishing
$p$-curvature.  In their recent work \cite{ogus-vol}, Ogus and
Vologodsky define a more general equivalence between modules with
integrable connection $(E, \nabla)$ whose $p$-curvature is nilpotent of
order less than $p$ and modules with Higgs field $(E', \theta')$ with $\theta'$
nilpotent of order less than $p$.  This equivalence depends on a
lifting $\tilde X' \to \tilde S$ modulo $p^2$ of $X'$ over $S$, where
$X'$ is the target of the relative Frobenius map $F_{X/S} : X\to X'$.
Furthermore, they found that if the $p$-curvature of $\nabla$ is nilpotent
of sufficiently low order, then the de Rham complex of $(E, \nabla)$ and
the Higgs complex of $(E', \theta')$ are naturally isomorphic in the
appropriate derived category, so that they again have the same
hypercohomology.

Log geometry was created to deal with problems in compactification and
singularities; thus, connections on log schemes provide a language for
studying differential equations with log poles.  This is an important
case to study since many natural connections do have log poles.  A log
scheme is a scheme $X$ with a sheaf of commutative monoids $\scrM_X$
and a map $\alpha : \scrM_X \to \scrO_X^{\times}$, where $\scrO_X^{\times}$ is the
multiplicative monoid of $\scrO_X$, such that $\alpha$ induces an
isomorphism $\alpha^{-1}(\scrO_X^*) \to \scrO_X^*$.  For example, given a
divisor with normal crossings $D$, we may define $\scrM_X := \scrO_X \cap
i_* \scrO_Y^*$, where $Y := X \setminus D$ and $i : Y \hookrightarrow X$ is the open
immersion, and we define $\alpha : \scrM_X \to \scrO_X$ to be the natural
inclusion.  Then for a morphism $X \to S$ of log schemes, Kato defines a
sheaf of relative logarithmic differentials $\Omega^1_{X/S}$ \cite{kato},
which allows the natural extension of the notion of modules with
integrable connection; in the case above of a divisor with normal
crossings, this is just the classical sheaf $\Omega^1_{X/S} (\log D)$ of
differentials with log poles along $D$.

However, on log schemes of characteristic $p$, in addition to the
$p$-curvature there is another obstruction, called the residue, to the
classical Riemann-Hilbert correspondence; hence the straightforward
generalization of the Riemann-Hilbert correspondence requires that the
connection $\nabla$ have vanishing residue in addition to vanishing
$p$-curvature.  Lorenzon \cite{lorenzon} corrected this by introducing
an $\MXbargp$-indexed ring $\calA_X^{gp}$ with canonical connection
$d$, such that defining $\calB_{X/S} := (\calA_X^{gp})^{d = 0}$, one
gets an equivalence between indexed $\calA_X^{gp}$-modules with
integrable connection compatible with $d$ whose $p$-curvature vanishes
and indexed $\calB_{X/S}$-modules.  This work was inspired by Tsuji's
work generalizing the fundamental exact sequence
\[ 0 \to \scrO_{X'}^* \to \scrO_X^* \overset{d\log}{\longrightarrow} Z^1_{X/S}
\overset{\pi^* - C}{\longrightarrow} \Omega^1_{X'/S} \to 0 \] to the case of log schemes,
replacing $\scrO_X^*$ by $\scrM_X^{gp}$ \cite{tsuji}.

The aim of this paper is to extend the theory of Ogus and Vologodsky
in \cite{ogus-vol} to the case of log schemes.  The main result is
that given a lifting $\calXS$ of $X'$ over $S$ modulo $p^2$, we get an
equivalence $C_{\calXS}$ between the category of $\scrO_X$-modules
with integrable connection $(E, \nabla)$ whose $p$-curvature is nilpotent
of level less than $p$ and the category of $\MXbargp$-indexed
$\calB_{X/S}$-modules with $\calB_{X/S}$-linear Higgs field $(E', \theta')$
which is nilpotent of level less than $p$.  (Here $\MXbargp$ denotes
the quotient $\MXgp / \scrO_X^*$.)  Furthermore, this specializes to
give an equivalence between $\scrO_X$-modules with connection whose
$p$-curvature and residue are both nilpotent and $\scrO_{X'}$-modules
with nilpotent Higgs field.  In both cases, the de Rham complex of
$(E, \nabla)$ and the Higgs complex of $(E', \theta')$ are once again isomorphic
in the derived category.

We give here a brief sketch of the construction of $C_{\calXS}$ and
its pseudoinverse $C_{\calXS}^{-1}$, assuming for simplicity we have a
lifting $\tilde X \to \tilde S$ in addition to the lifting of $X'$.
Then the sheaf of liftings of $F_{X/S} : X \to X'$ to a map $\tilde X \to
\tilde X'$ is a torsor over $F_{X/S}^* \Omega^1_{X'/S}$.  This sheaf is
representable by an affine scheme $\calL_{\calXS} = \bfSpec
\calK_{\calXS}$ over $X$; in addition, the sheaf extends naturally to
the crystalline site $\Crys(X / S)$, which induces a natural
connection on $\calK_{\calXS}$.  Defining $\calK_{\calXS}^{\calA} :=
\calK_{\calXS} \otimes_{\scrO_X} \AXgp$, we then have
\[ C_{\calXS}(E) \simeq (E \otimes_{\scrO_X} \calK_{\calXS}^{\calA})^{\nabla_{tot}} \]
with Higgs field induced by the $p$-curvature of $\calK_{\calXS}$, and
\[ C_{\calXS}^{-1}(E') \simeq (E' \otimes_{\calB_{X/S}}
\calK_{\calXS}^{\calA})^{\theta_{tot}}_0
\]
with connection acting on $\calK_{\calXS}$.

A key technical result which is central to the proofs in
\cite{ogus-vol} is the fact that for a scheme $X$ of characteristic
$p$, the ring of \PDop{}s $\DXS$ is an Azumaya algebra over its
center, which is isomorphic to $S^\cdot \scrT_{X'/S}$ via the map $D \mapsto D^p
- D^{(p)}$.  (This means that locally, after a flat base extension, it
is isomorphic to a matrix algebra.)  This is no longer true in general
in the case of log schemes; however, what we find is that if we define
$\tildeDXS := \AXgp \otimes_{\scrO_X} \DXS$ with the appropriate
multiplication, then $\tildeDXS$ is an indexed Azumaya algebra over
its center, which is isomorphic to $\calB_{X/S} \otimes_{\scrO_{X'}} S^\cdot
\scrT_{X'/S}$.  The first section is mainly devoted to an elaboration
and proof of this result.  The first subsection discusses the theory
of indexed Azumaya algebras; perhaps one interesting point which
should be mentioned here is the fact that we need to define $\scrHom$
in the indexed case to allow morphisms which shift degree.  The second
subsection reviews the construction of $\AXgp$ and the connection $d$
and explores the structure of $\calB_{X/S}$.  The third subsection
then defines $\tildeDXS$ and proves that it is an indexed Azumaya
algebra over its center.

The second section begins with a construction of the scheme
representing a general torsor over a locally free sheaf.  We then
define the $F_{X/S}^* \Omega^1_{X'/S}$-torsor of liftings of Frobenius
$\calL_{\calXS}$ and its extension to the crystalline site, which
gives us a corresponding crystal of $\scrO_{X/S}$-algebras
$\calK_{\calXS}$.  We proceed to calculate the corresponding
connection and its $p$-curvature explicitly (the latter turns out to
be simply the map $d : \scrO_{\calL} \to \Omega^1_{\calL / X} \simeq \scrO_{\calL}
\otimes_{\scrO_X} F_{X/S}^* \Omega^1_{X'/S}$, where by abuse of notation we also
let $\calL$ be the scheme representing the torsor $\calL_{\calXS}$)
and to discuss a functorality property of this construction.

The third section gives a construction of the Cartier transform
$C_{\calXS}$ described above; the key observation is that defining
\[ \check \calK_{\calXS}^{\calA} := \scrHom_{\scrO_X} (\calK_{\calXS},
\scrO_X) \otimes_{\scrO_X} \AXgp, \] then $\check \calK_{\calXS}^{\calA}$ is
a splitting module for $\tildeDXS^\gamma := \tildeDXS \otimes_{S^\cdot \scrT_{X'/S}}
\hat \Gamma_\cdot \scrT_{X'/S}$.  We thus get an equivalence between the
slightly larger categories of $\tildeDXS^\gamma$-modules and $\hat \Gamma_\cdot
\scrT_{X'/S} \otimes_{\scrO_{X'}} \calB_{X/S}$-modules, which have the
advantage of being $\otimes$-categories.  We then prove that after a sign
change the equivalence given by Azumaya algebra theory gives the
formulas above in the case of quasinilpotent modules.  We also give a
local version of the transform which works more generally for
$\tildeDXS \otimes_{S^\cdot \scrT_{X'/S}} \hat S^\cdot \scrT_{X'/S}$-modules, but
depends on a lifting $\tilde F$ of $F_{X/S}$, which rarely exists
globally.  Finally, in the last section, we prove the analogue of
Simpson's formality theorem, which states that if the $p$-curvature of
$\nabla$ is nilpotent of level less than $p - d$, where $d$ is the relative
dimension of $X / S$, then the de Rham complex of $(E, \nabla)$ is
isomorphic in the derived category of complexes of
$\scrO_{X'}$-modules to the Higgs complex of its Cartier transform.

In
their paper Ogus and Vologodsky proved a theorem of Barannikov and
Kontsevich, which states that if $X$ is a quasi-projective smooth
scheme over $\bbC$ and $f \in \Gamma(X, \scrO_X)$ induces a proper map $X \to
\bbA^1$, then the hypercohomologies of the complexes
\[ \scrO_X \overset{d + \land df}\longrightarrow \Omega^1_{X / \bbC} \overset{d + \land df}\longrightarrow \Omega^2_{X
  / \bbC} \cdots \]
and
\[ \scrO_X \overset{\land df}\longrightarrow \Omega^1_{X / \bbC} \overset{\land df}\longrightarrow \Omega^2_{X /
  \bbC} \cdots \]
have the same finite dimension in every degree.  For further research, it would
be interesting to see whether by using log geometry one can relax the
condition on $f$.

\begin{acknowledgements}
I would like to thank Arthur Ogus and Vadim Vologodsky for their
paper, without which this thesis would not have existed.  I would
especially like to thank Arthur Ogus who, as my advisor, suggested
this problem to me and provided helpful input at several points.

Part of this research was supported by the European Union contract
``Arithmetic Algebraic Geometry.''
\end{acknowledgements}

\section{Indexed Azumaya Algebras}

Let $X$ be a topological space, and $\calI$ a sheaf of abelian groups
on $X$.  Then an {\em $\calI$-indexed ring} on $X$ is a sheaf $\calA$
along with a degree map $\deg : \calA \to \calI$, a zero section $\calI
\to \calA$ of the degree map, a global section $1 \in \calA(X)$, an
addition map $\calA \times_{\calI} \calA \to \calA$ compatible with the maps
to $\calI$, and a multiplication map $\calA \times \calA \to \calA$ fitting
into a commutative diagram
\[
\begin{CD}
  \calA \times \calA @> \cdot >> \calA \\
  @VVV @VVV \\
  \calI \times \calI @> + >> \calI,
\end{CD}
\]
satisfying the usual ring identities when either side is defined.  (In
other words, we only allow addition of sections of equal degree, and
require that multiplication adds degrees.  Also, for each section $i \in
\calI$, there is a corresponding zero section in $\calA_i :=
\deg^{-1}(i)$.)  Also, if $\calJ$ is a sheaf of $\calI$-sets, and
$\calA$ is an $\calI$-indexed ring, we have a similar definition of
$\calJ$-indexed $\calA$-modules.  A {\em $\calJ$-indexed homomorphism}
of $\calJ$-indexed $\calA$-modules is then a map $\scrE \to \scrF$ which
is $\calA$-linear in the usual sense and which respects degrees.  See
\cite{lorenzon} for more details.

\subsection{Basics}

Let $X$ be a topological space, $\calI$ a sheaf of abelian groups on $X$, and
$\calA$ an $\calI$-indexed ring on $X$.  In order to develop the
theory of indexed Azumaya algebras, we will need to allow
homomorphisms which shift degree.

\begin{definition}
  Let $\calJ$ be a sheaf of $\calI$-sets.
  \begin{enumerate}
  \item Let $\scrE$ be a $\calJ$-indexed $\calA$-module, and $j \in
    \calJ$ a section.  Then $\scrE(j) := \scrE \times_{\calJ} \calI$, where
    the map $\calI \to \calJ$ is addition of $j$.  This is an
    $\calI$-indexed object via the projection to $\calI$, and an
    $\calA$-module via the action on $\scrE$.  (For $i \in \calI$, we
    have $\scrE(j)_i \simeq \scrE_{i + j}$.  Thus, if $\calJ = \calI$ with
    the standard action, this agrees with the usual definition.)
  \item Let $\scrE$ and $\scrF$ be $\calI$ and $\calJ$-indexed
    $\calA$-modules, respectively.  Then $\scrHom_\calA(\scrE, \scrF)$
    is the $\calJ$-indexed $\calA$-module such that a section on an
    open set $U\subseteq X$ consists of a section $j \in \calJ(U)$ and
    an $\calI$-indexed $\calA$-linear homomorphism $\phi : \scrE|_U
    \to (\scrF|_U)(j)$, with the obvious restriction maps, and the
    degree map of projection to $j$.
  \end{enumerate}
\end{definition}

Now if $\calA$ is commutative and $M$ is a locally free
$\calI$-indexed $\calA$-module of finite rank (with generators not
necessarily of degree zero), then $\scrE := \scrEnd_\calA(M)$ is an
$\calA$-algebra with center $\calA$.  In fact, we have the following
result:

\begin{theorem}
  Let $\calJ$ be a sheaf of $\calI$-sets, and let $M$ and $\scrE$ be
  as above.
  \label{thm:azmbasic}
  \begin{enumerate}
  \item The functor $E \mapsto M \otimes E$ induces an equivalence of categories
    from the category of $\calJ$-indexed $\calA$-modules to the
    category of $\calJ$-indexed left $\scrE$-modules, with
    quasi-inverse $F \mapsto \scrHom_\scrE(M, F)$.
  \item If $M'$ is another locally free $\calA$-module of the same
    rank as $M$, with a structure of left $\scrE$-module, then the
    natural map $\scrE \to \scrEnd_\calA(M')$ is an isomorphism.
  \end{enumerate}

  \begin{proof}
    We have natural $\calJ$-indexed maps
    \begin{align*}
      \eta_E : E \to \scrHom_\scrE(M, M \otimes E), & ~ e \mapsto (m \mapsto m \otimes e) \\
      \epsilon_F : M \otimes \scrHom_\scrE(M, F) \to F, & ~ m \otimes \phi \mapsto \phi(m).
    \end{align*}
    It is sufficient to show these maps are isomorphisms locally;
    thus, let $\{ m_1, \ldots, m_r \}$ be a basis for $M$, and let $d_i =
    \deg(m_i) \in \calI$.  Also, let $\alpha_{ij} \in \scrE$ be the
    homomorphism which sends $m_k$ to $\delta_{jk} m_i$.  Note that
    $\deg(\alpha_{ij}) = d_i - d_j$.
    
    Now $\eta_E$ is clearly injective.  For surjectivity, let $\phi
    \in \scrHom_\scrE(M, M \otimes E)$.  Suppose $\phi(m_1) = \sum_{i=1}^r m_i
    \otimes e_i$; then since $\alpha_{i1} m_1 = m_i$, we must have $\phi(m_i) =
    \alpha_{i1} \phi(m_1) = m_i \otimes e_1$ for each $i$.  Therefore, $\phi(m) = m
    \otimes e_1$ for each $m \in M$, so $\phi = \eta_E(e_1)$.
    
    To show $\epsilon_F$ is an isomorphism, we define $\epsilon_F^{-1}$ as follows: for each
    $i$, let $\scrE_{\cdot i}$ denote the sub $\calA$-module of $\scrE$
    generated by $\alpha_{1i}, \ldots, \alpha_{ri}$.  Then we have an isomorphism
    of left $\scrE$-modules $M \to \scrE_{\cdot i}(-d_i)$ which sends $m_j$ to
    $\alpha_{ji}$.  For $f\in F$ of degree $d \in \calJ$, we now let
    $\phi_i(f) \in \scrHom_{\scrE} (M, F)$ be the composition
    \[ M \overset{\sim}{\to} \scrE_{\cdot i}(-d_i) \overset{\cdot f}{\to} F(-d_i
    + d), \] and define $\epsilon_F^{-1}(f) = \sum_{i=1}^r m_i \otimes
    \phi_i(f)$.  Again, this is a $\calJ$-indexed map.  To see it is
    indeed the inverse, note that any element of $M \otimes
    \scrHom_\scrE(M, F)$ can be written uniquely as a sum
    $\sum_{i=1}^r m_i \otimes \psi_i$.  This gets mapped by
    $\epsilon_F$ to $f := \sum_i \psi_i(m_i) \in F$; applying the
    above construction, we see that $\phi_i(f)$ sends $m_j$ to
    $\alpha_{ji} \sum_k \psi_k(m_k) = \sum_k \psi_k(\alpha_{ji} m_k) =
    \psi_i(m_j)$, so $\phi_i(f) = \psi_i$.  For the other direction,
    we see that $\epsilon_F$ maps $\sum_{i=1}^r m_i \otimes \phi_i(f)$
    to $\sum_i \phi_i(f)(m_i) = \sum_i \alpha_{ii} \cdot f = 1 \cdot
    f$.  This completes the proof of (i).
    
    For (ii), since $M$ and $M'$ are both locally free of rank $r$, it
    suffices to show that the localized map $\End_{\calA_x}(M_x) \to
    \End_{\calA_x}(M'_x)$ is an isomorphism for each $x \in X$; thus,
    we may assume $X$ is a one-point space, and $M$ and $M'$ are free.
    Also, by localizing at an arbitrary $\calI$-indexed prime ideal of
    $\calA$, we may assume that $\calA$ is local.  Now by (1) we have
    an isomorphism $M \otimes N \to M'$, where $N = \Hom_\scrE(M, M')$.
    Since $M$ is free, this implies $N$ is a projective
    $\calA$-module.  Also, applying $\alpha_{11}$ to both sides of the
    isomorphism $M \otimes N \to M'$ gives an isomorphism $N(-d_1) \simeq (\calA
    \cdot m_1) \otimes N \to \alpha_{11} M'$, so $N$ is finitely generated.  Now
    the standard proof that finitely generated projective modules over
    a local ring are free extends to the indexed case, so we see that
    $N$ is free.  Since $M$ and $M'$ are both of rank $r$, $N$ must
    have rank 1.  Therefore, if $N$ has a free generator of degree $d$,
    then $M' \simeq M(-d)$, and the result follows.
  \end{proof}
\end{theorem}

(In fact, this proof shows that if $\calA_x$ is local as an
$\calI_x$-indexed ring for each $x\in X$, then locally $M' \simeq M(i)$ as
left $\scrE$-modules for some $i \in \calI$.)

\begin{remark}
\label{rmk:azmadj}
It is easy to check that assuming only that $M$ is an $\scrE$-module,
where $\scrE$ is an $\calA$-algebra, the natural transformations $\eta$
and $\epsilon$ described in the above proof form the unit and counit,
respectively, of an adjunction which makes $M \otimes \cdot$ the left adjoint of
$\scrHom_{\scrE}(M, \cdot)$.
\end{remark}

\begin{definition}
  Let $\calA$ be a commutative $\calI$-indexed ring, and $\scrE$ an
  $\calI$-indexed $\calA$-algebra.  For a commutative $\calI$-indexed
  $\calA$-algebra $\calB$, we say $\scrE$ {\em splits over $\calB$} with
  splitting module $M$ if $\scrE \otimes_\calA \calB \simeq \scrEnd_\calB(M)$
  for some locally free $\calI$-indexed $\calB$-module $M$.  We say
  $\scrE$ is an {\em Azumaya algebra} over $\calA$ if there is some
  commutative $\calI$-indexed $\calA$-algebra $\calB$, faithfully flat
  over $\calA$, such that $\scrE$ splits over $\calB$.
\end{definition}

\begin{corollary}
  If $\scrE$ is an Azumaya algebra over $\calA$ of rank $r^2$, and
  there exists a locally free $\calI$-indexed $\calA$-module $M$ of
  rank $r$ with a structure of left $\scrE$-module, then $\scrE$ is
  split over $\calA$ with splitting module $M$.

  \begin{proof}
    Let $\calB$ be an $\calA$-algebra, faithfully flat over $\calA$,
    over which $\scrE$ splits, and let $M'$ be a splitting module.
    Then $M \otimes_\calA \calB$ has a structure of left $\scrE \otimes_\calA
    \calB$-module, and it is a locally free $\calB$-module of rank
    $r$.  On the other hand, if $M'$ has rank $s$, then $\scrE
    \otimes_\calA \calB \simeq \scrEnd_\calB(M')$ must have rank $s^2$ over
    $\calB$, so $s = r$.  Thus, by (ii) of theorem
    \ref{thm:azmbasic}, the natural map $\scrE \otimes_\calA \calB \to
    \scrEnd_\calB(M \otimes_\calA \calB) \simeq \scrEnd_\calA(M) \otimes_\calA
    \calB$ is an isomorphism.  Since $\calB$ is faithfully flat over
    $\calA$, this implies that the natural map $\scrE \to
    \scrEnd_\calA(M)$ is an isomorphism.
  \end{proof}
\end{corollary}

Note that since the transpose gives an isomorphism
$\scrE^{\operatorname{op}} \to \scrEnd_\calA(\check M)$, where $\check
M = \scrHom_\calA(M, \calA)$, the above theory works equally well for
right $\scrE$-modules.

\subsection{Indexed Algebras Associated to a Log Structure}

Recall from \cite{lorenzon} that to any log scheme $X$ we may
associate a canonical $\MXbargp$-indexed $\scrO_X$-algebra $\AXgp$,
which corresponds to the exact sequence of abelian groups $0 \to
\scrO_X^* \to \MXgp \to \MXbargp \to 0$.  In particular, for $\bar s \in
\MXbargp$, $(\AXgp)_{\bar s}$ is the invertible sheaf on $X$
corresponding to the $\scrO_X^*$-torsor given by the inverse image of
$\bar s$ in $\MXgp$; thus, for a section $s \in \MXgp$, we have a
corresponding basis element $e_s$ of $(\AXgp)_{\bar s}$, where $\bar
s$ is the image of $s$ in $\MXbargp$.  We then have a canonical
connection $d$ on $\AXgp$ characterized by the formula $d e_s = e_s \otimes
\dlog s$, and we define $\calB_{X/S}$ to be the kernel of $d$.  Note
that since $d$ is multiplicative, $\calB_{X/S}$ is a subring of
$\AXgp$.

The construction of $\AXgp$ is functorial in the following way: given
a morphism $f : X\to Y$ of fine log schemes, we have a commutative
diagram
\[
\begin{CD}
  0 @>>> \scrO_X^* @>>> (f^* \scrM_Y)^{\gp} @>>> f^{-1}
  \Mbargp{Y} @>>> 0 \\
  @. @VVV @VVV @VVV @. \\
  0 @>>> \scrO_X^* @>>> \MXgp @>>> \MXbargp @>>> 0.
\end{CD}
\]
This induces a natural cartesian diagram
\[
\begin{CD}
  f^* \calA_Y^{\gp} @>>> \AXgp \\
  @VVV @VVV \\
  f^{-1} \Mbargp{Y} @>>> \MXbargp.
\end{CD}
\]
For $s \in \scrM_Y^{\gp}$, the map $f^* \calA_Y^{\gp} \to \AXgp$ sends
$1 \otimes e_s$ to $e_{f^* s}$.

Now recall the relative Frobenius diagram for log schemes:

\[
\begin{diagram}
  \node{X} \arrow{e,t}{F_{X/S}} \arrow{se}
  \node{X'} \arrow{s} \arrow{se,t}{\pi_{X/S}} \\
  \node[2]{X''} \arrow{e} \arrow{s}
  \node{X} \arrow{s} \\
  \node[2]{S} \arrow{e,t}{F_S}
  \node{S.}
\end{diagram}
\]
Here the bottom square is cartesian in the category of fine log
schemes, and $F_{X/S}$ is uniquely determined by the requirement that
$F_{X/S}$ is purely inseparable, and the map $X' \to X''$ is \'etale.
Note that our notation for $X'$ and $X''$ is opposite that used by
Kato in \cite{kato} but agrees with \cite{ogus-higgs}.  We then have
the Cartier isomorphisms $C_{X/S} : \calH^q(F_{X/S *} \Omega^\cdot_{X/S})
\overset{\sim}\to \Omega^q_{X'/S}$ of $\scrO_{X'}$-modules.  In fact, we have
the following generalization:

\begin{proposition}
  \label{prop:cart-isom}
  Assume $X \to S$ is a smooth morphism of fine schemes.  Then there
  are canonical isomorphisms
  \[ \calH^i(\AXgp \otimes_{\scrO_X} \Omega^\cdot_{X/S}, d) \simeq \calB_{X/S}
  \otimes_{\scrO_{X'}} \Omega^i_{X'/S}, \] where $(\AXgp \otimes_{\scrO_X}
  \Omega^\cdot_{X/S}, d)$ is the de Rham complex of $\AXgp$ with the
  canonical connection $d$.

  \begin{proof}
    (This is a straightforward generalization of the proof in
    \cite[4.12]{kato}.)  We have a natural map $\calB_{X/S}
    \otimes_{\scrO_{X'}} \Omega^i_{X'/S} \to \calH^i(\AXgp \otimes_{\scrO_X} \Omega^\cdot_{X/S},
    d)$ which extends the inverse Cartier isomorphism $C_{X/S}^{-1} :
    \Omega^i_{X'/S} \to \calH^i(\Omega^\cdot_{X/S}, d)$.  Showing this map is an
    isomorphism is an \'etale local problem, so we may assume we have a
    chart $(Q \to \scrM_S, P \to \scrM_X, Q \to P)$ of $X \to S$ such that $P$
    and $Q$ are finitely generated integral monoids, the map $Q^{\gp}
    \to P^{\gp}$ is injective and the torsion part of its cokernel is a
    finite group of order prime to $p$, and the natural map $X \to S
    \times_{\Fp[Q]} \Fp[P]$ is an isomorphism.  Then $X' \simeq S \times_{\Fp[Q]}
    \Fp[H]$, where $H := P \cap (p P^{gp} + Q^{\gp})$.
    
    Now for $s \in P^{\gp}$, let $\bar s$ be the image of $s$ in
    $\MXbargp$; then using $e_s$ as a basis element for
    $(\AXgp)_{\bar s}$, the degree $\bar s$ part of $(\AXgp \otimes
    \Omega^\cdot_{X/S}, d)$ is isomorphic to the de Rham complex
    $(\Omega^\cdot_{X/S}, d + \dlog s \land)$.  We can now decompose this
    complex as a direct sum of complexes $C_\nu^\cdot$ for $\nu \in \Fp \otimes
    (P^{\gp} / Q^{\gp}) \simeq P^{\gp} / H^{\gp}$, where
    \[ C_\nu^i := \scrO_S \otimes_{\Fp[Q]} \Fp[P \cap \nu] \otimes_{\Fp} \land^i (\Fp \otimes
    (P^{\gp} / Q^{\gp})) \] with the differential $a \mapsto (s + \nu) \land
    a$.  (Here $\Fp[P \cap \nu]$ denotes the submodule of $\Fp[P]$
    generated by elements of $P \cap \nu$.)  This complex is exact if
    $s + \nu \neq 0$, while if $s + \nu = 0$, then $\calH^i(C_\nu^\cdot) \simeq
    C_\nu^i \simeq C_\nu^0 \otimes_{\scrO_{X'}} \Omega^i_{X'/S}$.  Hence
    $\calH^i(\Omega^\cdot_{X/S}, d + \dlog s \land) \simeq \calH^0(\Omega^\cdot_{X/S}, d +
    \dlog s \land) \otimes_{\scrO_{X'}} \Omega^i_{X'/S}$.  However,
    $\calH^0(\Omega^\cdot_{X/S}, d + \dlog s \land)$ is the kernel of $d +
    \dlog s$, which by definition corresponds to
    $(\calB_{X/S})_{\bar s}$.  Since the image of $P^{\gp}$
    generates $\MXbargp$, we are done.
  \end{proof}
\end{proposition}

\begin{remark}
  Since $(\AXgp \otimes \Omega^\cdot_{X/S}, d)$ has zero $p$-curvature, the Higgs
  complex is $(\AXgp \otimes_{\scrO_{X'}} \Omega^\cdot_{X'/S}, 0)$, and the kernel of
  the natural connection on $\AXgp \otimes_{\scrO_{X'}} \Omega^i_{X'/S}$ is
  $\calB_{X/S} \otimes_{\scrO_{X'}} \Omega^i_{X'/S}$.  Thus, this is a special
  case of \cite[3.1.1]{ogus-higgs}, which states that $\calH^i(E \otimes
  \Omega^\cdot_{X/S}, \nabla) \simeq \calH^i(E \otimes F_{X/S}^* \Omega^\cdot_{X'/S}, -\psi)^\nabla$.
  Similarly, for $i < p$, this is a special case of the main result of
  chapter 4, since in notation yet to be developed we will have
  $C_{\calXS}(\AXgp, d) \simeq (\calB_{X/S}, 0)$.
\end{remark}

The proof of this isomorphism also gives us insight into the structure
of $\calB_{X/S}$.  Namely, we see that given a local chart $(Q \to
\scrM_S, P \to \scrM_X, Q \to P)$ as in the proof, then for $m \in P^{\gp}$,
\[ (\calB_{X/S})_m = e_m \scrO_S \otimes_{\Fp[Q]} \Fp[(-m + H^{\gp}) \cap P] \subseteq
e_m \scrO_S \otimes_{\Fp[Q]} \Fp[P] = (\AXgp)_m. \]  This immediately gives
the following:

\begin{corollary}
  The natural map $\calA_{X'}^{\gp} \to F_{X/S *} \calA_X^{\gp}$ factors
  through $F_{X/S *} \calB_{X/S}$.  Furthermore, $F_{X/S *}
  \calB_{X/S}$ is generated as an $F_{X/S *} \MXbargp$-indexed module
  over $\calA_{X'}^{\gp}$ by elements of the form $\alpha(m) e_{-m}$ for $m
  \in P$.

  \begin{proof}
    Locally, $\Mbargp{X'}$ is generated by $H^{gp}$, and for $m \in
    H^{gp}$, we have $e_m \scrO_{X'} = e_m \scrO_S \otimes_{\Fp[Q]}
    \Fp[H^{gp}] \subseteq \calB_{X/S}$.  On the other hand, the above
    description of $\calB_{X/S}$ shows that it is generated as an indexed
    $\scrO_S$-module by elements of the form $\alpha(-m + h) e_m$, where $h
    \in H^{\gp}$ is such that $-m + h \in P$.  However, then $\alpha(-m +
    h) e_m = e_h [\alpha(-m + h) e_{-(-m + h)}]$, and $e_h \in
    \calA_{X'}^{\gp}$.
  \end{proof}
\end{corollary}

\begin{example}
  Consider the monoid $P = \langle a, b, c : a + b = 2c \rangle$, and let $X :=
  \Spec (P \to k[P])$ for some field $k$ of characteristic $p$.  We then
  claim that $F_{X/k}$ is not flat.  To see this, we can identify $P$
  with $\{ (m, n) \in \ints^2 : m \geq |n| \}$ by sending $a$ to $(1, 1)$,
  $b$ to $(1, -1)$, and $c$ to $(1, 0)$.  We then see that $\{ e^m : m
  \in P - P H^+ \}$ forms a basis for $k[P] / k[P H^+]$ as a
  vector space over $k[H] / k[H^+] \simeq k$.  However, for some cosets $S$
  of $H^{\gp} = p \ints^2$, $S \cap P$ has more than one minimal element;
  for example, $((-1, -1) + H^{\gp}) \cap P$ has minimal elements $(p-1,
  p-1)$ and $(p-1, -1)$.  Therefore, $k[P] \otimes_{k[H]} k$ has dimension
  strictly greater than $p^2$.  But localizing to $D(H^+)$, we get
  $k[P^{\gp}]$ over $k[H^{\gp}]$, which is free of rank $p^2$, since
  choosing a set of representatives $S$ of $P^{\gp} / H^{\gp}$, we get
  a basis $\{ e^s : s \in S \}$ of $k[P^{\gp}]$.

  On the other hand, the sets $S \cap P$ for $S$ a coset of $H^{\gp}$
  form a partition of $P$ into $p^2$ sets.  Therefore, if we choose a
  set of representatives $S$ of $P^{\gp} / H^{\gp}$, then $\{ e_s : s \in
  S \}$ forms a basis for $\AXgp$ as a $\calB_{X/S}$-module (since $-S$
  is also a set of representatives).
\end{example}

Generalizing this argument, we get:

\begin{proposition}
  (See also \cite[2.6]{lorenzon}.)
  Let $f : X\to S$ be a smooth morphism of log schemes of characteristic
  $p$, with $X$ of relative dimension $r$ over $S$.  Then $\AXgp$ is
  locally free of rank $p^r$ as a $\calB_{X/S}$-module.  In fact, if
  $m_1, \ldots, m_r \in \MXgp$ is a logarithmic system of coordinates, then
  letting $\theta_i := e_{m_i}$ and $\theta^I := \prod_{k=1}^r \theta_k^{I_k}$ for $I \in
  \ints^r$, then $\{ \theta^I : I \in \{ 0, \ldots, p-1 \}^r \}$ forms a basis.

  \begin{proof}
    Locally, if we have a chart $(Q \to \scrM_S, P \to \scrM_X, Q \to P)$,
    we get as in the example above that $\AXgp$ is free over
    $\calB_{X/S}$.  All we need to show is that $P^{\gp} / H^{\gp}$
    has $p^r$ elements.  However, $P^{\gp} / H^{\gp} \simeq \Fp \otimes (P^{\gp}
    / Q^{\gp})$ is a vector space over $\Fp$, and when we base extend
    to $\scrO_X$ we get $\Omega^1_{X/S}$, which is locally free of rank
    $r$.  Therefore, $\Fp \otimes (P^{\gp} / Q^{\gp})$ has dimension $r$,
    and the result follows.  Given a logarithmic system of
    coordinates, we may choose the chart so that $m_1, \ldots, m_r$ are in
    the image of $P^{\gp} \to \MXgp$, in which case $\{ \sum_{k=1}^r I_k m_k
    : I \in \{ 0, \ldots, p-1 \}^r \}$ forms a set of representatives of
    $P^{\gp} / H^{\gp}$.
  \end{proof}
\end{proposition}

\subsection{The Azumaya Algebra $\tildeDXS$}

Our main application of this theory will be to an extension of the
ring of PD differential operators on $X$, where $f : X\to S$ is a smooth
morphism of log schemes of characteristic $p$.  We briefly review the
standard notation: first, $D(1)$ is the logarithmic PD envelope of the
diagonal embedding $\Delta : X \to X \times_S X$, and $\calD(1) := \scrO_{D(1)}$,
which is considered to be an $\scrO_X$-algebra via the projection $p_1
: D(1) \to X$.  Then for $m \in \MXgp$, $\eta_m \in \calD(1)$ is the unique
element of the PD ideal in $\calD(1)$ such that $1 + \eta_m = \alpha_{D(1)}
(p_2^* m - p_1^* m)$.  If $m_1, \ldots, m_r \in \MXgp$ is a logarithmic
system of coordinates (i.e. $\dlog m_1, \ldots, \dlog m_r$ form a basis for
$\Omega^1_{X/S}$), then $\{ D_I : I \in \nats^r \}$ denotes the basis of the
ring of \PDop{}s dual to the basis $\{ \eta^{[I]} = \prod_{k=1}^r
\eta_{m_k}^{[I_k]} : I \in \nats^r \}$ of $\calD(1)$.  (In this paper we
will not be using the alternate basis $\{ \zeta^{[I]} \}$ of $\calD(1)$,
where $\zeta_m := \log(1 + \eta_m)$, except in the appendix; this basis has
the property that the dual basis $\{ D^I : I \in \nats^r \}$ satisfies
$D^I \circ D^J = D^{I + J}$.)

\begin{definition}
  Let $\DXS$ denote the ring of \PDop{}s on $X$ relative to $S$.  Then
  $\tildeDXS$ is the $\MXbargp$-indexed $\scrO_X$-algebra $\AXgp
  \otimes_{\scrO_X} \DXS$, given the unique multiplication such that $\AXgp
  \to \tildeDXS$, $a \mapsto a \otimes 1$ and $\DXS \to \tildeDXS$, $\phi \mapsto 1 \otimes \phi$ are
  ring homomorphisms, $(a \otimes 1) (1 \otimes \phi) = a \otimes \phi$, and for each log
  derivation $\partial$,
  \[ (1 \otimes \partial) (a \otimes 1) = (\partial . a) \otimes 1 + a \otimes \partial. \]
  (Here $\partial . a$ denotes
  the evaluation of $\partial$ at $a$ as opposed to composition of \PDop{}s,
  where we give $\AXgp$ the $\DXS$-module structure corresponding to
  the canonical connection $d$ on $\AXgp$.)
\end{definition}

A straightforward induction shows that in local coordinates, $(1 \otimes
D_N) (a \otimes 1) = \sum_{I \leq N} \binom{N}{I} (D_{N-I} . a) \otimes D_I$.  Thus,
this agrees with the definition of $\tildeD_{X,\calI}^{(0)}$ from
\cite{montagnon}.

For notational convenience, we will often treat the given maps $\AXgp
\to \tildeDXS$, $\DXS \to \tildeDXS$ as
embeddings, and thus write $a \phi$ in place of $a \otimes \phi$.  We now have
an inclusion
\[ \tildeDXS \hookrightarrow \PDDiff(\AXgp, \AXgp) =
\scrHom_{\scrO_X} (\divpow(1) \otimes_{\scrO_X} \AXgp, \AXgp), \] which
sends a section $a \in \AXgp$ to multiplication by $a$ and is
compatible with the map $\DXS \to \PDDiff(\AXgp, \AXgp)$
corresponding to the canonical connection $d$.  (The relation $\partial
\circ a = (\partial . a) + a \partial$
holds in $\PDDiff(\AXgp, \AXgp)$ since $d$ is multiplicative.)  In
fact, it is not hard to see that the image is exactly the set of \PDop{}s
$\phi$ such that $\phi(\tau \otimes e_s) = \phi(\tau (1 + \eta_s)) \cdot e_s$ for $\tau \in
\divpow(1)$, $s \in \MXgp$.

The following result indicates the significance of the ring
$\tildeDXS$.  Recall from \cite{lorenzon} that a connection $\nabla : \scrE
\to \scrE \otimes_{\scrO_X} \Omega^1_{X/S}$ on a $\calJ$-indexed $\AXgp$-module $E$
is {\em admissible} if and only if $\nabla(ae) = a \nabla(e) + e \otimes da$ whenever
$a \in \AXgp$, $e \in \scrE$.

\begin{proposition}
  Let $\calJ$ be a sheaf of $\MXbargp$-sets, and let $\scrE$ be a
  $\calJ$-indexed $\AXgp$-module with integrable connection $\nabla$ on
  $\scrE$.  Then $\nabla$ is admissible if and only if the $\AXgp$-module
  structure on $\scrE$ and the $\DXS$-module structure on
  $\scrE$ corresponding to $\nabla$ extend to a $\tildeDXS$-module
  structure on $\scrE$.

  \begin{proof}
    First, note that if we do have such an extension, it must be given
    by $(a \otimes \phi) e = a \nabla_\phi e$.  Now
    admissibility of $\nabla$ is equivalent to the
    condition that for every log derivation $\partial$ on $X$, $\nabla_\partial (a
    e) = a \nabla_\partial(e) + \partial(a) e$.  Rewriting, this tells us that
    \[ (1 \otimes \partial) [(a \otimes 1) e] = (\partial . a \otimes 1 + a \otimes
    \partial) e = [(1 \otimes \partial) (a \otimes 1)] e. \]
    From this the desired equivalence easily follows.
  \end{proof}
\end{proposition}

We can also provide a simple crystalline interpretation of the
category of $\tildeDXS$-modules.  In particular, if $g : T_1 \to T_2$ is
a morphism in $\Crys(X/S)$, then the natural map $g^{-1} \Mbargp{T_2}
\to \Mbargp{T_1}$ is an isomorphism, so $g^* \calA_{T_2}^{\gp} \to
\calA_{T_1}^{\gp}$ is also an isomorphism.  Hence the functor $T \mapsto
\calA_T^{\gp}$ gives a crystal of $\MXbargp$-indexed
$\scrO_{X/S}$-algebras.  Now letting $T$ be the logarithmic PD
envelope of the diagonal in $X \times_S X$ with the two projections $p_1,
p_2 : T \to X$, for $s \in \scrM_X^{\gp}$ we calculate
\[ p_2^*(e_s) = e_{(0, s)} = e_{(s, 0)} e_{(-s, s)} = p_1^*(e_s) (1 + \eta_s). \]
(Here $(-s, s) \in \scrM_T^{*}$, so $e_{(-s, s)} = \alpha_T(-s, s) = 1 +
\eta_s$.)  Thus the HPD stratification on $\AXgp$ sends $1 \otimes e_s$ to
$e_s \otimes (1 + \eta_s)$, and in particular, the connection this crystal
defines on $\AXgp$ agrees with $d$.  It is then easy to show:
\begin{proposition}
  Suppose we have a locally nilpotent $\calJ$-indexed $\AXgp$-module
  $E$ with connection $\nabla$.  Then $\nabla$ is admissible if and only if in
  the corresponding crystal of $\calJ$-indexed $\scrO_{X/S}$-modules,
  the transition maps $\theta_g : g^* E_{T_2} \to E_{T_1}$ are $g^*
  \calA_{T_2}^{\gp}$-linear.
\end{proposition}
Hence the category of locally nilpotent $\calJ$-indexed
$\tildeDXS$-modules is equivalent to the category of crystals of
$\calJ$-indexed $\AXgp$-modules on $\Crys(X / S)$.

\begin{theorem}
  Let $\frakZ = S^\cdot(\scrT_{X'/S})$, which we embed into $\DXS$ via
  the map $\partial \mapsto \partial^p - \partial^{(p)}$.  Let $\tilde \frakZ$ be the center
  of $\tildeDXS$.  Then $\tilde \frakZ = \calB_{X/S}
  \otimes_{\scrO_{X'}} \frakZ$ as a subring of
  $\tildeDXS$.

  \begin{proof}
    If $f \in \calB_{X/S}$, then since $\partial f = 0$ for any log
    derivation $\partial$, we see that $f$ commutes with $\DXS$, and it
    obviously commutes with $\AXgp$; thus, since $\DXS$ and $\AXgp$
    generate $\tildeDXS$ as a ring, we get $f \in \tilde
    \frakZ$.
    
    Now choose a system of logarithmic coordinates $m_1, \ldots, m_r \in
    \MXgp$, and let $\theta_i = e_{m_i}$.  Finally, let $M = \{ 0, 1, \ldots, p-1
    \}^r$.  Then $\AXgp$ is a locally free $\calB_{X/S}$-module with
    local basis $\{ \theta^I : I \in M \}$.  Thus, $\tildeDXS$ is generated as
    a $\calB_{X/S}$-algebra by $\{ \theta_i, D_{\epsilon_i} \}$.  Expressing a
    section $\phi \in \tildeDXS$ as a sum $\sum_N f_N D_N$, and recalling that
    $\phi^\flat(1, e_s) = \phi^\flat(s, 1) e_s$ for $\phi \in \DXS$ and $s \in \MXgp$
    \cite[4.1.1]{montagnon}, we now calculate
    \begin{align*}
      [\phi, \theta_i] & = \sum_N N_i f_N \theta_i D_{N - \epsilon_i}; \\
      [\phi, D_{\epsilon_i}] & = - \sum_N (D_{\epsilon_i} f_N) \cdot D_N.
    \end{align*}
    
    Therefore, $\phi \in \tilde\frakZ$ if and only if $f_N \neq 0$ only for
    $p \mid N$, and $d f_N = 0$ for each $N$, i.e.~$f_N \in \calB_{X/S}$.
    Since the embedding $\frakZ \hookrightarrow \DXS$ sends
    $D_{\epsilon_i}$ to $D_{\epsilon_i}^p - D_{\epsilon_i} = D_{p \epsilon_i}$, so that its image is
    generated as an $\scrO_{X'}$-module by $\{ D_{p N} : N \in \nats^r
    \}$, this completes the proof.
  \end{proof}
\end{theorem}

Note that we also get that the centralizer of $\AXgp$ is $\AXgp
\otimes_{\scrO_X} \frakZ$.  Denote this by $\calC_X$, and
consider $\tildeDXS$ as a right $\calC_X$-module by
multiplication on the right.

\begin{theorem}
  The map $\tildeDXS \otimes_{\tilde \frakZ} \calC_X \to
  \scrEnd_{\calC_X} (\tildeDXS)$ by multiplication on the left
  by $\tildeDXS$ and on the right by $\calC_X$ is an
  isomorphism.

  \begin{proof}
    It suffices to work locally, so choose a local system of
    logarithmic coordinates $m_1, \discretionary{}{}{}\ldots,
    \discretionary{}{}{} m_r$.  Then locally,
    $\tildeDXS \otimes_{\tilde \frakZ} \calC_X$ has a basis $\{ 1
    \otimes \theta^I : I \in M = \{ 0, 1, \ldots, p-1 \}^r \}$ as a left
    $\tildeDXS$-module.  Also, $\scrEnd_{\calC_X}
    (\tildeDXS)$ has a basis $\{ \alpha_I : I \in M \}$ as a left
    $\tildeDXS$-module, where $\alpha_I$ is the unique
    homomorphism which sends $D_J$ to $\delta_{IJ}$.  (This is because
    $\tildeDXS$ has a basis $\{ D_I : I\in M \}$ as a right
    $\calC_X$-module.)
    
    We now calculate that $1 \otimes \theta_i$ acts on $\tildeDXS$ by
    sending $D_N$ to $\theta_i D_N + N_i \theta_i D_{N - \epsilon_i}$.  Thus, setting
    $\beta_i = \theta_i^{-1} \otimes \theta_i - 1 \otimes 1$,
    $\beta_i$ acts by sending $D_N$ to $N_i D_{N - \epsilon_i}$.  Now letting $\beta^K := \prod_{i=1}^r \beta_i^{K_i}$
    for $K \in M$, this implies that $\beta^K$ acts on $\tildeDXS$
    by sending $D_I$ to $\frac{I!}{(I-K)!} D_{I-K}$ if $I \geq K$, and 0
    otherwise.

    We thus have
    \[ \beta^J = \sum_{I \leq J} (-1)^{|J - I|} \binom{J}{I} (\theta^{-I} \otimes \theta^I). \]
    Thus, enumerating $M$ in some order compatible with the product
    partial order, the transition matrix from the set $\{ \beta^J \}$ to
    the basis $\{ 1 \otimes \theta^I : I \in M \}$ of
    $\tildeDXS \otimes_{\tilde \frakZ} \calC_X$ is upper
    triangular, with units on the diagonal, so $\{ \beta^J : J\in M \}$ is
    also a basis for $\tildeDXS \otimes_{\tilde \frakZ} \calC_X$.
    Similarly, letting $\beta^J$ also denote the corresponding
    endomorphism on $\tildeDXS$, we have
    \[ \beta^J = \sum_{I \geq J, I\in M} \left( \frac{I!}{(I - J)!} D_{I-J}
      \right) \alpha_I. \]
    Thus, the transition matrix from the set $\{ \beta^J : J \in M \}$ to
    the basis $\{ \alpha_I : I \in M \}$ of $\scrEnd_{\calC_X}
    (\tildeDXS)$ is lower triangular with units
    on the diagonal, so $\{ \beta^J : J\in M \}$ is also a basis for
    $\scrEnd_{\calC_X} (\tildeDXS)$.
  \end{proof}
\end{theorem}

\begin{corollary}
  The indexed ring $\tildeDXS$ is an Azumaya algebra over its
  center $\tilde \frakZ$ of rank $p^{2r}$, where $r$ is the relative
  dimension of $X$ over $S$.

  \begin{proof}
    Since $\AXgp$ is locally free and thus faithfully flat as a
    $\calB_{X/S}$-module, $\calC_X \simeq \AXgp \otimes_{\calB_{X/S}} \tilde
    \frakZ$ is a faithfully flat extension of $\tilde \frakZ$, and by
    the previous theorem $\tildeDXS$ splits over $\calC_X$
    with splitting module $\tildeDXS$.  Since
    $\tildeDXS$ has a basis $\{ D_I : I \in \{ 0, \ldots, p-1
    \}^r \}$ as a right $\calC_X$-module, the rank must be $p^r$, so
    $\tildeDXS$ has rank $p^{2r}$ over $\tilde\frakZ$.
  \end{proof}
\end{corollary}

As an application, we can recover the main result of \cite{lorenzon}
as follows.  Consider $\calB_{X/S}$ as a $\tilde\frakZ$-algebra via
base extension of the quotient map $S^\cdot(\scrT_{X'/S}) \to
S^\cdot(\scrT_{X'/S}) / S^+(\scrT_{X'/S}) \simeq \scrO_{X'}$, and let
$\tildeD_0 := \tildeDXS \otimes_{\tilde\frakZ} \calB_{X/S}$.  Then
$\AXgp$ is a locally free $\calB_{X/S}$-module of rank $p^r$
(compatible with the $\tilde\frakZ$-algebra structure of $\calB_{X/S}$
since the connection $d$ is $p$-integrable), which has a structure of
left module over $\tildeD_0$.  Therefore, $\tildeDXS$ splits
over $\calB_{X/S}$, and for $\calJ$ a sheaf of $\MXbargp$-sets we get an
equivalence between the $\calJ$-indexed $\calB_{X/S}$-modules and the
$\calJ$-indexed left $\tildeD_0$-modules, or equivalently the
$\calJ$-indexed $\AXgp$-modules with integrable, $p$-integrable,
admissible connections.  This equivalence sends a $\calB_{X/S}$-module
$\scrF$ to $\AXgp \otimes_{\calB_{X/S}} \scrF$ with the connection acting
on $\AXgp$, and it sends an $\AXgp$-module $\scrE$ with connection
$\nabla$ to $\scrHom_{\tildeD_0} (\AXgp, \scrE) \simeq \scrE^\nabla$.

\section{The Fundamental Extension}

\begin{definition}
  Suppose we are given a morphism $f : X\to S$ of fine log schemes of
  characteristic $p$.  Then a {\em lifting of $f$ modulo $p^n$} is a
  map $\tilde f : \tilde X\to \tilde S$ of fine log schemes flat over
  $\ints / p^n$ which fits into a cartesian square
  \[
  \begin{CD}
    X @>>> \tilde X \\
    @V f VV @VV \tilde f V \\
    S @>>> \tilde S,
  \end{CD}
  \]
  where $S \to \tilde S$ is the closed immersion defined by $p$.
\end{definition}

Note that since $\Spec(\ints / p) \to \Spec(\ints / p^n)$ is an exact
closed immersion, so are the base extensions $X \to \tilde X$ and $S \to
\tilde S$.  If $f$ is smooth, respectively \'etale, resp.~integral,
resp.~exact, so is $\tilde f$.  We will mostly be interested in
liftings modulo $p^2$, which is what we will mean if we do not specify
$n$.

For the rest of the paper, we suppose we are given an integral smooth
morphism $f : X\to S$ of fine log schemes of characteristic $p$, with a
given lifting of $X' \to S$ modulo $p^2$, $\tilde X' \to \tilde S$.  We
denote these data by $\calXS := (X\to S, \tilde X' \to \tilde S)$.  For
example, if $f : X\to S$ has a given lifting $\tilde f : \tilde X\to
\tilde S$, and $F_S$ has a lifting $F_{\tilde S} : \tilde S\to \tilde
S$, we may define $\tilde X''$ to be the pullback of $F_{\tilde S}$
and $\tilde f$.  Then since $X' \to X''$ is \'etale, there is a unique
lifting $\tilde X' \to \tilde X''$.  Alternately, a smooth lifting
exists if $H^2(X', \scrT_{X'/S}) = 0$, in particular if $X$ is affine
or if $X$ is a curve over a field $k$, and since we assume $X \to S$
integral, a smooth lifting is automatically a lifting in the sense
defined above.

As in \cite{ogus-vol}, we will use a lifting $\tilde X' / \tilde S$ to
construct a canonical sheaf of $\scrO_X$-algebras $\calK_{\calXS}$
with multiplicative connection, along with a natural filtration $N_\cdot$
on $\calK_{\calXS}$ such that $N_i \calK_{\calXS} \simeq S^i (N_1
\calK_{\calXS})$ and a short exact sequence
\[ 0 \to \scrO_X \to N_1 \calK_{\calXS} \to F_{X/S}^* \Omega^1_{X'/S} \to 0 \]
in which the maps are horizontal.

\subsection{Torsors over Locally Free Sheaves}

In this section we give the general construction which gives rise to
$\calK_{\calXS}$.

\begin{lemma}
  \label{lemma:crys-tors}
  Let $X$ be a ringed topos, with a locally free $\scrO_X$-module
  $\scrT$ of finite rank; let $\Omega := \check\scrT = \scrHom_{\scrO_X}
  (\scrT, \scrO_X)$.  Suppose we have a $\scrT$-torsor $\calL$.

  \begin{enumerate}
  \item There exists an $\scrO_X$-algebra $\calK$ such that $\calL$ is
    isomorphic to the sheaf of $\scrO_X$-algebra homomorphisms $\calK
    \to \scrO_X$.  Furthermore, there is a canonical cotorsor
    structure on $\calK$ consisting of a coaction map $\calK \to \calK
    \otimes S^\cdot \Omega$ and a cosubtraction map $S^\cdot \Omega
    \to \calK \otimes \calK$, such that the action of $D \in \scrT$ on
    $\calL$ corresponds to composition with the automorphism
    \[ \calK \to \calK \otimes S^\cdot \Omega \overset{\id \otimes
      S^\cdot D} \longrightarrow \calK. \]
  \item There is a canonical
    filtration $N_\cdot$ on $\calK$ such that letting $\calE := N_1
    \calK$, we have $S^i \calE \simeq N_i \calK$ via the natural map
    $S^i \calE \to \calK$.
  \item There is a canonical locally split exact sequence
    \[ 0 \to \scrO_X \to \calE \to \Omega \to 0. \]
  \item If $X$ is the Zariski topos on a scheme, then $\calK$ is
    quasicoherent, and the affine
    scheme $\bfT(\calL) := \bfSpec \calK$ over $X$ represents the
    functor $\calL$.
  \item If $X$ is the crystalline topos on a scheme over $S$, and
    $\scrT$ is a crystal of $\scrO_{X/S}$-modules, then $\calE$ and
    $\calK$ are also crystals of $\scrO_{X/S}$-modules and
    $\scrO_{X/S}$-algebras, respectively.
  \end{enumerate}

  \begin{proof}
    We first construct $\calE$ as the subsheaf of $\scrHom(\calL,
    \scrO_X)$ consisting of morphisms $\phi : \calL \to \scrO_X$ such that
    for any local section $a$ of $\calL$, the function $\scrT \to
    \scrO_X$, $D \mapsto \phi(a + D) - \phi(a)$, is $\scrO_X$-linear.  Note that
    for any object $U$ of $X$, if this is true for one $a \in \calL(U)$,
    it is true for all $a \in \calL(U)$, and in this case, the function
    is independent of the choice of $a$.

    To get the exact sequence, note that any constant function on
    $\calL$ is in $\calE$; this gives an injection $\scrO_X \to \calE$.
    Also, given $\phi \in \calE$, as we remarked above, the map $\scrT \to
    \scrO_X$, $D \mapsto \phi(a + D) - \phi(a)$, is independent of the choice of
    $a$; thus, gluing gives a global map $\scrT \to \scrO_X$, which
    induces a map $\calE \to \Omega$.  To see this is a locally split
    surjection, given a section $a \in \calL$, we define a splitting
    $\sigma_a : \Omega\to \calE$ by $\sigma_a(\omega)(b) = \langle \omega, b - a \rangle$.  Now it is clear
    that the sequence in (iii) is exact, since the kernel of $\calE \to \Omega$
    comprises exactly the constant functions $\calL \to \scrO_X$.
    
    Now the injection $\scrO_X \to \calE$ induces an injection
    $S^n(\calE) \to S^{n+1}(\calE)$ for each $n \geq 0$; we now define
    $\calK := \underset{\longrightarrow}{\lim} ~ S^n(\calE)$.  Then
    $\calK$ has a natural filtration where $N_i \calK$ is the image of
    $S^i \calE$ in the direct limit.  Also, an $\scrO_X$-algebra
    homomorphism $\calK \to \scrO_X$ is equivalent to an
    $\scrO_X$-linear map $\calE \to \scrO_X$ which maps the constant
    function with value 1 to 1; that is, a splitting of the above
    exact sequence on the left.  The sheaf of such splittings forms a
    $\scrT$-torsor, where $D \in \scrT$ acts on the set of splittings
    by addition of the composition of $D : \Omega \to \scrO_X$ with
    the projection $\calE \to \Omega$.  Moreover, there is a morphism
    of torsors from $\calL$ to the torsor of left splittings given by
    $a \mapsto \epsilon_a$, where $\epsilon_a : \calE \to \scrO_X$ is
    evaluation at $a$.  This establishes the desired property of
    $\calK$.  This, in turn, easily implies that $\bfSpec \calK$
    represents the functor $\calL$ on a Zariski topos.

    We now have a map $\calE \to \calE \otimes \Omega$ which sends
    $\phi \in \calE$ to $\phi \otimes 1 + 1 \otimes \omega_\phi$,
    where $\omega_\phi$ is the image of $\phi$ in $\Omega$.  Since a
    constant function $a$ maps to $a \otimes 1$, this extends to a
    coaction map $\calK \to \calK \otimes S^\cdot \Omega$.  Similarly,
    a local splitting $\sigma_a$ of the sequence in (iii) induces a map
    $1 \otimes \sigma_a - \sigma_a \otimes 1 : \Omega \to \calE
    \otimes \calE$, which is independent of $a$ since the image of
    $\sigma_{a_2} - \sigma_{a_1}$ consists of constant functions.
    Thus, these local maps glue to a global map $\Omega \to \calE
    \otimes \calE$, which extends to a cosubtraction map $S^\cdot
    \Omega \to \calK \otimes \calK$.
    
    Finally, if $X$ is a crystalline topos, then any extension of
    crystals is automatically another crystal.  Thus, $\calE$ and
    therefore $\calK$ are also crystals.
  \end{proof}
\end{lemma}

\subsection{Liftings of Frobenius as a Crystal}

We consider the crystalline sites $\Crys(X / S)$ and $\Crys(X / \tilde
S)$.  If $(U, T)$ is an object of $\Crys(X / S)$, then since $U$ is
defined by a divided power ideal $I_T$ in $T$, $a^p = 0$ for any $a \in
I_T$.  Therefore, the map $T\to T''$ factors through $U''$.  We thus get
a commutative diagram
\[
\begin{diagram}
\node{U} \arrow{e} \arrow{s}
\node{U'} \arrow{s} \\
\node{T} \arrow{e} \arrow{ne,..}
\node{U''}
\end{diagram}
\]
However, since the map $U' \to U''$ is \'etale, and $U\to T$ is an exact
thickening, there exists a unique map $T\to U'$ making the above diagram
commute.  We let $f_{T/S} : T \to X'$ denote the composition of this map
with the open embedding $U' \to X'$.  Note that the pullback map
$f_{T/S}^* : \Omega^1_{X'/S} \to f_{T/S *} \Omega^1_{T/S}$ is zero.  To see this,
we may assume $X = U$; in this case, let $i' : U'\to T'$ be the closed
immersion.  Then after composing with the natural surjection $i'^* :
i'^* \Omega^1_{T'/S} \to \Omega^1_{U'/S}$, we get the map $f_{T/S}^* i'^* =
F_{T/S}^* = 0$.

\begin{definition}
  Let $\tilde T$ be an object of $\Crys(X / \tilde S)$ which is flat
  over $\tilde S$, and let $T$ be the closed subscheme defined by
  $p$.
  \begin{enumerate}
  \item A {\em lifting of $f_{T/S}$ to $\tilde T$} is a lifting
    $\tilde F : \tilde T\to \tilde X'$ over $\tilde S$ modulo $p^2$.
    $\calL_{\calXS}(\tilde T)$ is the set of all such liftings, and
    $\calL_{\calXS, \tilde T}$ is the sheaf on $\tilde T$ of local
    liftings of $f_{T/S}$.
  \item A {\em lifting of $f_{T/S}$} is a pair $(\tilde T', \tilde F)$
    where $\tilde T'$ is a flat object of $\Crys(X / \tilde S)$ and
    $\tilde F : \tilde T' \to \tilde X'$ is a lifting of $f_{T/S}$ to
    $\tilde T'$.  An isomorphism $(\tilde T_1, \tilde F_1) \to (\tilde
    T_2, \tilde F_2)$ is a map $\tilde g : \tilde T_1 \to \tilde T_2$
    which reduces to the identity on $T$ and such that $\tilde F_1 =
    \tilde F_2 \circ \tilde g$.  $\calL_{\calXS, T}$ is the sheaf on $T$
    associated to the presheaf of isomorphism classes of local
    liftings of $f_{T/S}$.
  \end{enumerate}
\end{definition}
Since $\tilde X'$ is smooth over $\tilde S$, liftings of $f_{T/S}$ to
$\tilde T$ exist locally on $T$, so $\calL_{\calXS, \tilde T}$ has
nonempty stalks.  A map $\tilde g : \tilde T_1 \to \tilde T_2$ induces a
map $\calL_{\calXS}(\tilde g) : \calL_{\calXS}(\tilde T_2) \to
\calL_{\calXS}(\tilde T_1)$, defined by $\tilde F \mapsto \tilde F \circ \tilde
g$.  Now $\calL_{\calXS, \tilde T}$ is a torsor over
$\scrHom_{\scrO_T} (f_{T/S}^* \Omega^1_{X'/S}, \scrO_T) \simeq f_{T/S}^*
\scrT_{X'/S}$.  We will show that this torsor has a natural structure
of a sheaf on $\Crys(X / S)$, compatible with the crystal structure of
the crystal $F_{X/S}^* \scrT_{X'/S}$, which maps $T \in \Crys(X / S)$ to
$f_{T/S}^* \scrT_{X'/S}$.  (The corresponding connection is the
Frobenius descent connection on $F_{X/S}^* \scrT_{X'/S}$.)

\begin{lemma}
  Let $\tilde T_1$ and $\tilde T_2$ be two flat objects of $\Crys(X /
  \tilde S)$, and let $\tilde g_1, \tilde g_2 : \tilde T_1 \to \tilde
  T_2$ be two morphisms with the same reduction $g : T_1 \to T_2$ modulo
  $p$.
  \begin{enumerate}
  \item $\calL_{\calXS}(\tilde g_1) = \calL_{\calXS}(\tilde g_2)$ as
    maps $\calL_{\calXS}(\tilde T_2) \to \calL_{\calXS}(\tilde T_1)$.
  \item For any flat object $\tilde T \in \Crys(X / \tilde S)$, the
    natural map $\calL_{\calXS, \tilde T} \to \calL_{\calXS, T}$ is an
    isomorphism, where $T$ is the reduction of $\tilde T$ modulo $p$.
  \end{enumerate}

  \begin{proof}
    For $\tilde F \in \calL_{\calXS}(\tilde T_2)$, consider the
    commutative diagram
    \[
    \begin{CD}
      T_1 @> inc \, \circ \, g >> \tilde T_2 @> \tilde F >> \tilde X' \\
      @VVV @VV \Delta V @VV \Delta' V \\
      \tilde T_1 @> (\tilde g_1, \tilde g_2) >> \tilde T_2 \times_{\tilde
        S} \tilde T_2 @> \tilde F \times \tilde F >> \tilde X' \times_{\tilde S}
      \tilde X'.
    \end{CD}
    \]
    Then the left square induces the map $h : g^* \Omega^1_{T_2 / S} \to
    \scrO_{T_1}$ expressing the difference between $\tilde g_1$ and
    $\tilde g_2$, and the right square induces the map $\tilde F^* :
    \tilde F^* \Omega^1_{\tilde X' / \tilde S} \to \Omega^1_{\tilde T_2 /
      \tilde S}$.  Thus, overall the diagram induces the map
    $h \circ g^* (f_{T_2/S}^*) : g^* f_{T_2/S}^* \Omega^1_{X'/S} \to \scrO_{T_1}$,
    which is zero since $f_{T_2/S}^* = 0$.  This implies that $(\tilde
    F \circ \tilde g_1, \tilde F \circ \tilde g_2) : \tilde T_1 \to \tilde X'
    \times_{\tilde S} \tilde X'$ factors through $\Delta'$, so $\tilde F
    \circ \tilde g_1 = \tilde F \circ \tilde g_2$.
    
    To prove (ii) it suffices to prove the map is an isomorphism on
    stalks.  Thus, for $t \in \tilde T$, if $(\tilde T', \tilde F)$ is
    a lifting of $f_{T/S}$ on a neighborhood of $t$, then locally at
    $t$ we have an isomorphism $\tilde T \simeq \tilde T'$, showing
    surjectivity.  For injectivity, if $\tilde F, \tilde F'\in
    \calL_{\calXS}(\tilde T)$ become equal in $\calL_{\calXS, T}$,
    then there is an automorphism of $\tilde T$ which reduces to the
    identity modulo $p$ and which carries $\tilde F$ to $\tilde F'$.
    But then by (i), $\tilde F = \tilde F'$.
  \end{proof}
\end{lemma}

Hence given a map $g : T_1 \to T_2$, we can define a map $\theta_g : g^{-1}
\calL_{\calXS, T_2} \to \calL_{\calXS, T_1}$ by gluing the maps
$\calL_{\calXS}(\tilde g)$ for local liftings $\tilde g$ of $g$.  It
is easy to see this map satisfies the cocycle condition $\theta_{hg} = \theta_h
\circ h^{-1} \theta_g$.  Also, since $\tilde F_2 \circ \tilde g - \tilde F_1 \circ
\tilde g = g^* (\tilde F_2 - \tilde F_1) : g^* f_{T_2 / S}^*
\Omega^1_{X'/S} \to \scrO_{T_1}$ for $\tilde F_i \in \calL_{\calXS, T_2}$ if
$\tilde g : \tilde T_1 \to \tilde T_2$ lifts $g$, $\theta_g$ respects the
torsor structure.

Thus from (\ref{lemma:crys-tors}) we get a crystal $\calK_{\calXS}$
with filtration $N_\cdot$, along with an extension
\begin{myequation}
  \label{eqn:fundext}
  0 \to \scrO_{X/S} \to \calE_{\calXS} \to F_{X/S}^* \Omega^1_{X'/S} \to 0
\end{myequation}
of crystals, where $\calE_{\calXS} := N_1 \calK_{\calXS}$.

We have an alternate construction of $\calE_{\calXS}$, and hence
$\calK_{\calXS} \simeq \underset{\longrightarrow}{\lim} ~ S^n(\calE_{\calXS})$, as follows:
given an object $\tilde T \in \Crys(X / \tilde S)$ which is flat over
$\tilde S$, we define $\calE_{\calXS, \tilde T}$ to be the logarithmic
conormal sheaf of the closed immersion
\[ T \overset{\Gamma}\longrightarrow T \times_S X' \overset{inc}\longrightarrow \tilde T \times_{\tilde
  S} \tilde X', \] where $T$ is the reduction of $\tilde T$ modulo
$p$, and $\Gamma$ is the graph of $f_{T/S} : T \to X'$.  The functoriality
properties of the conormal sheaf allow us to define maps $\theta_{\tilde g}
: g^* \calE_{\calXS, \tilde T_2} \to \calE_{\calXS, \tilde T_1}$ for
$\tilde g : \tilde T_1 \to \tilde T_2$, which it is possible to show
depend only on the reduction $g : T_1 \to T_2$ of $\tilde g$ modulo $p$.
We then glue together $\calE_{\calXS, \tilde T}$ to get the sheaves
$\calE_{\calXS, T}$.  The exact sequence above comes from the exact
sequence $\Gamma^* \check \scrN_{\tilde T \times \tilde X' / T \times X'} \to \check
\scrN_{\tilde T \times \tilde X' / T} \to \check \scrN_{T \times X' / T} \to 0$.
To connect this to the previous construction, observe that given
$\tilde F \in \calL_{\calXS, \tilde T}$, we can refactor the closed
immersion above as
\[ T \overset{inc}\longrightarrow \tilde T \overset{\Gamma_{\tilde F}}\longrightarrow \tilde T
\times_{\tilde S} \tilde X'. \] The corresponding map $\check \scrN_{\tilde
  T \times \tilde X' / T} \to \check \scrN_{\tilde T / T}$ gives a map
$\calE_{\calXS, \tilde T} \to \scrO_T$.  Thus, for each element of
$\calE_{\calXS, \tilde T}$, we get a map $\calL_{\calXS, \tilde T} \to
\scrO_T$; it is then possible to show this map is in $N_1
\calK_{\calXS, T}$, and that the induced map $\calE_{\calXS, \tilde T}
\to N_1 \calK_{\calXS, T}$ is an isomorphism of extensions of $f_{T/S}^*
\Omega^1_{X'/S}$ by $\scrO_T$.

\subsection{Explicit Formulas: Connection and $p$-curvature}

We now calculate the corresponding connection on $\calE_{\calXS}$ and
its $p$-curvature.  The answer involves the following construction:
suppose we have a lifting $\tilde F : \tilde X\to \tilde X'$ of
$F_{X/S}$.  Then since $F_{X/S}^* : \Omega^1_{X'/S} \to F_{X/S *}
\Omega^1_{X/S}$ is the zero map, $\tilde F^* : \Omega^1_{\tilde X' / \tilde S}
\to \tilde F_* \Omega^1_{\tilde X / \tilde S}$ is a multiple of $p$.  Hence
there is a unique map $\zeta_{\tilde F} : \Omega^1_{X' / S} \to F_{X/S *}
\Omega^1_{X / S}$ making the following diagram commute:
\[
\begin{CD}
  \Omega^1_{\tilde X' / \tilde S} @> \tilde F^* >> \tilde F_* \Omega^1_{\tilde
    X / \tilde S} \\
  @VVV @A [p] AA \\
  \Omega^1_{X' / S} @> \zeta_{\tilde F} >> F_{X/S *} \Omega^1_{X / S}.
\end{CD}
\]
Now since (\ref{eqn:fundext}) is an exact sequence of crystals, the
only nontrivial part of the connection $\nabla$ on $\calE_{\calXS} \simeq
\scrO_X \oplus F_{X/S}^* \Omega^1_{X'/S}$ is the part which maps $F_{X/S}^*
\Omega^1_{X'/S}$ to $\scrO_X \otimes \Omega^1_{X/S}$, and similarly for the
$p$-curvature $\psi(\nabla)$.  (Note that the splitting
$\calE_{\calXS} \simeq \scrO_X \oplus F_{X/S}^* \Omega^1_{X'/S}$
depends on a lifting $\tilde F$ of $F_{X/S}$.)

\begin{proposition}
  \label{prop:econn}
  Let $\tilde F$ and $\zeta_{\tilde F}$ be as above.
  \begin{enumerate}
  \item (Mazur's Formula) $\zeta_{\tilde F}$ factors through $F_{X/S *}
    Z^1_{X/S}$ and induces a splitting of the exact sequence
    \[ 0 \to F_{X/S *} B^1_{X/S} \to F_{X/S *} Z^1_{X/S}
    \overset{C_{X/S}}{\longrightarrow} \Omega^1_{X'/S} \to 0, \]
    where $C_{X/S}$ is the Cartier operator.
  \item The connection $\nabla$ on $\calE_{\calXS}$ satisfies
    \[
    (\epsilon_{\tilde F} \otimes \id) \circ \nabla \circ \sigma_{\tilde F} = -\zeta_{\tilde F}
    : F_{X/S}^* \Omega^1_{X'/S} \to \Omega^1_{X/S}. \]
  \item The $p$-curvature $\psi$ of $\nabla$ is equal to the composition
    \[ \calE_{\calXS} \to F_{X/S}^* \Omega^1_{X'/S} \overset{i \otimes \id}{\longrightarrow}
    \calE_{\calXS} \otimes F_{X/S}^* \Omega^1_{X'/S}, \]
    where $i : \scrO_X \to \calE_{\calXS}$ is the inclusion map.
  \end{enumerate}

  \begin{proof}
    For $m \in \scrM_X^{\gp}$, let $\tilde m \in \scrM_{\tilde X}^{\gp}$
    and $\tilde m' \in \scrM_{\tilde X'}^{\gp}$ be liftings of $m$ and
    $\pi^* m \in \scrM_{X'}^{\gp}$, respectively.  Then $\tilde u :=
    \tilde F^*(\tilde m') - p \tilde m$ reduces to 0 in
    $\scrM_X^{\gp}$, so $\tilde u \in \scrM_{\tilde X}^*$ and
    $\alpha_{\tilde X}(\tilde u) = 1 + [p] b$ for some $b \in \scrO_X$.
    Therefore,
    \[ [p] \zeta_{\tilde F}(\dlog(\pi^* m)) = \tilde F^*(\dlog(\tilde m')) =
    \dlog(p \tilde m) + \frac{d\alpha(\tilde u)}{\alpha(\tilde u)} = [p]
    \dlog(m) + [p] db. \] Thus $\zeta_{\tilde F}(\dlog(\pi^* m)) =
    \dlog(m) + db \in Z^1_{X/S}$ and $C_{X/S}(\zeta_{\tilde F}(\dlog(\pi^*
    m))) = C_{X/S}(\dlog(m)) = \dlog(\pi^* m)$ as required.  Since $\{
    \dlog(\pi^* m) : m \in \scrM_X^{\gp} \}$ generates $\Omega^1_{X'/S}$ as
    an $\scrO_{X'}$-module, this proves (i).
    
    To prove (ii), let $T$ be the first infinitesimal neighborhood of
    $X$ in $X \times_S X$ and $\tilde T$ the first infinitesimal
    neighborhood of $\tilde X$ in $\tilde X \times_{\tilde S} \tilde X$,
    with its two natural projections $\tilde p_1, \tilde p_2 : \tilde
    T \to \tilde X$.  Then for $\omega \in F_{X/S}^* \Omega^1_{X'/S}$, it is easy to
    see that $\theta_{p_i} \sigma_{\tilde F}(\omega) = \sigma_{\tilde F \circ \tilde p_i}
    (p_i^* \omega) = \sigma_{\tilde F \circ \tilde p_i} (1 \otimes \omega)$, so
    \begin{align*}
      (\epsilon_{\tilde F} \otimes \id) \nabla \sigma_{\tilde F}(\omega) & = [\theta_{p_2} \sigma_{\tilde F}(\omega) -
      \theta_{p_1} \sigma_{\tilde F}(\omega)](\tilde F \circ \tilde p_1) = [\sigma_{\tilde F \circ
        \tilde p_2}(\omega) -
      \sigma_{\tilde F \circ \tilde p_1}(\omega)](\tilde F \circ \tilde p_1) \\
      & = \langle \omega, \tilde F \circ \tilde p_1
      - \tilde F \circ \tilde p_2 \rangle.
    \end{align*}
    However, by definition, the map $\tilde F \circ \tilde p_1 - \tilde F
    \circ \tilde p_2 : f_{T/S}^* \Omega^1_{X'/S} \to \scrO_T$ is induced by
    pullback along the map $(\tilde F \circ \tilde p_2, \tilde F \circ \tilde
    p_1) : \tilde T \to \tilde X_1'$, where $\tilde X_1'$ is the first
    infinitesimal neighborhood (in the logarithmic sense) of the
    diagonal in $\tilde X' \times_{\tilde S} \tilde X'$.  Now this map
    can be factored as the transposition map $(\tilde p_2, \tilde p_1)
    : \tilde T \to \tilde T$, which induces the map $-1 : \Omega^1_{\tilde X
      / \tilde S} \to \Omega^1_{\tilde X / \tilde S}$, followed by the map $\tilde F
    \times \tilde F : \tilde T \to \tilde X_1'$, which induces the map
    $\tilde F^* : \tilde F^* \Omega^1_{\tilde X' / \tilde S} \to \Omega^1_{\tilde
      X / \tilde S}$.  Therefore, $\tilde F \circ \tilde p_1 - \tilde F \circ
    \tilde p_2$ is the map induced by $-\tilde F^*$, which is exactly
    $-\zeta_{\tilde F}$.

    Now for the $p$-curvature, since $\psi = 0$ on $\scrO_X$ and on the
    quotient $F_{X/S}^* \Omega^1_{X'/S}$, $\psi$ induces an endomorphism on
    $F_{X/S}^* \Omega^1_{X'/S}$; we need to show this endomorphism is the
    identity.  It suffices to show that for $D' \in \scrT_{X'/S}$ and
    $\omega \in \Omega^1_{X'/S}$, we have $\psi_{D'}(\sigma_{\tilde F}(1 \otimes \omega)) =
    F_{X/S}^* \langle \omega, D' \rangle$.  Furthermore, it suffices to consider the
    case in which $\omega = \dlog(\pi^* m)$ for $m \in \scrM_X^{\gp}$ and
    $D' = \pi^* D$ for $D = (\partial, \delta) \in \scrT_{X/S}$.  Thus, suppose
    $\zeta_{\tilde F}(\omega) = \dlog(m) + db$ for $b \in \scrO_X$.  Then
    $\nabla_D(\sigma_{\tilde F}(1 \otimes \omega)) = -\delta m - \partial b \in \scrO_X$, so
    \[ \nabla_D^p(\sigma_{\tilde F}(1 \otimes \omega)) = -\partial^{p-1} (\delta m) -
    \partial^{(p)} b. \]
    On the other hand, since $D^{(p)} = (\partial^{(p)}, \partial^{p-1} \circ \delta +
    F_X^* \circ \delta)$, we have
    \[ \nabla_{D^{(p)}}(\sigma_{\tilde F}(1 \otimes \omega)) = -\partial^{p-1}(\delta m) -
    F_X^*(\delta m) - \partial^{(p)} b. \] Thus $\psi_{D'}(\sigma_{\tilde F}(1 \otimes
    \dlog(m))) = F_X^*(\delta m) = F_{X/S}^* \langle \omega, D' \rangle$ as required.
  \end{proof}
\end{proposition}

\begin{proposition}
  Let $\psi : \calK_{\calXS} \to \calK_{\calXS} \otimes F_{X/S}^* \Omega^1_{X'/S}$
  be the $p$-curvature of the canonical connection on
  $\calK_{\calXS}$.  Then there is a canonical commutative diagram
  \[
  \begin{CD}
    \calK_{\calXS} @> \psi >> \calK_{\calXS} \otimes F_{X/S}^* \Omega^1_{X'/S} \\
    @VVV @VVV \\
    \scrO_{\bfT(\calL)} @> d >> \Omega^1_{\bfT(\calL) / X},
  \end{CD}
  \]
  in which the vertical arrows are isomorphisms.  (Here we give
  $\bfT(\calL_{\calXS})$ the log structure induced by that of $X$.)

  \begin{proof}
    Since $\calK_{\calXS}$ is a crystal of $\scrO_{X/S}$-algebras, the
    corresponding connection on $\calK_{\calXS}$ satisfies the Leibniz
    rule, so its $p$-curvature does also.  Hence $\psi$ is a derivation
    which annihilates $\scrO_X$; it will suffice to show it is the
    universal derivation on $\calK_{\calXS}$ over $\scrO_X$.  Given
    any derivation $D : \calK_{\calXS} \to \scrF$ over $X$, the
    restriction of $D$ to $\calE_{\calXS}$ factors through the
    quotient map to $F_{X/S}^* \Omega^1_{X'/S}$.  The $\scrO_X$-linear map
    $F_{X/S}^* \Omega^1_{X'/S} \to \scrF$ induces a $\calK_{\calXS}$-linear
    map $\calK_{\calXS} \otimes F_{X/S}^* \Omega^1_{X'/S} \to \scrF$.  Using the
    formula for $\psi$ from (\ref{prop:econn}), we easily see that the
    composition of this map with $\psi$ agrees with $D$ on
    $\calE_{\calXS}$.  Since $\calE_{\calXS}$ generates
    $\calK_{\calXS}$ as a ring and both $D$ and $\psi$ satisfy the
    Leibniz rule, this shows that $D$ factors through $\psi$.
  \end{proof}
\end{proposition}

Global liftings of $F_{X/S}$ rarely exist; this limits the
effectiveness of generating elements
of $\calE_{\calXS}$ using $\sigma_{\tilde F}$ and the inclusion of
$\scrO_X$.  On the other hand, global
liftings of $\pi_{X/S} : X' \to X$ exist more often.  For example, in
the situation in which the lifting $\tilde X' \to \tilde S$ comes from
a lifting $F_{\tilde S}$ of $F_S$, we can define $\tilde \pi$ to be the
map $\tilde X' \to \tilde X''$ composed with the projection $\tilde X''
= \tilde X \times_{F_{\tilde S}} \tilde S \to \tilde X$.  The following
result gives an alternate way to generate elements of $\calE_{\calXS}$
using $\tilde \pi$ instead of $\tilde F$.

\begin{proposition}
  Suppose that $\tilde \pi : \tilde X' \to \tilde X$ lifts $\pi_{X/S} : X'\to
  X$.  Then there is a unique map $\beta_{\tilde \pi}(\tilde m) :
  \scrM_{\tilde X}^{\gp} \to \calE_{\calXS}$ such that for each section
  $\tilde m \in \scrM_{\tilde X}^{\gp}$ and each local lifting $\tilde
  F$ of $F_{X/S}$, we have
  \[ 1 + [p] \beta_{\tilde \pi}(\tilde m) (\tilde F) = \alpha_{\tilde X}(\tilde
  F^* \tilde \pi^* \tilde m - p \tilde m). \]
  Furthermore:

  \begin{enumerate}
  \item $\beta_{\tilde \pi}$ induces a surjective map $\scrM_{\tilde
      X}^{\gp} \otimes_{\ints} \scrO_X \to \calE_{\calXS}$.
  \item We have
    \[ [p] \alpha(m)^p \beta_{\tilde \pi}(\tilde m)(\tilde F) = \tilde F^* \tilde
    \pi^* \alpha(\tilde m) - \alpha(\tilde m)^p. \]
    (Thus if we define $\delta_{\tilde \pi} : \scrO_{\tilde X} \to
    \calE_{\calXS}$ so that $[p] \delta_{\tilde \pi}(\tilde a)(\tilde F) =
    \tilde F^* \tilde \pi^* \tilde a - \tilde a^p$, then
    \[ \alpha(m)^p
    \beta_{\tilde \pi}(\tilde m) = \delta_{\tilde \pi}(\alpha(\tilde m)). \]
    This map
    $\delta_{\tilde \pi}$ matches the construction of $\delta_{\tilde \pi}$ in the
    nonlogarithmic case \cite{ogus-vol}.)
  \item The following diagram commutes:
    \[
    \begin{CD}
      \scrO_X @> \alpha_{\tilde X}^{-1}(1 + [p] \cdot) >> \scrM_{\tilde
        X}^{\gp} @>>> \scrM_X^{\gp} \\
      @VV F_X^* V @VV \beta_{\tilde \pi} V @VV 1 \otimes (\dlog \circ \pi^*) V \\
      \scrO_X @>>> \calE_{\calXS} @>>> F_{X/S}^* \Omega^1_{X'/S}.
    \end{CD}
    \]
  \item If $\tilde F : \tilde X \to \tilde X'$ lifts $F_{X/S}$, then for
    each $\tilde m \in \scrM_{\tilde X}^{\gp}$,
    \[ \beta_{\tilde \pi}(\tilde m) = \beta_{\tilde \pi}(\tilde m)(\tilde F) +
    \sigma_{\tilde F}(1 \otimes \dlog(\pi^* m)), \]
    where $m \in \MXgp$ is the reduction of $\tilde m$.
  \item For every local section $\tilde m$ of $\scrM_{\tilde X}^{\gp}$
    lifting some $m \in \scrM_X^{\gp}$,
    \begin{align*}
    \nabla \beta_{\tilde \pi}(\tilde m) & = -1 \otimes \dlog m, \\
    \psi \beta_{\tilde \pi}(\tilde m) & = 1 \otimes (1 \otimes \dlog(\pi^* m)).
    \end{align*}
  \end{enumerate}
  \begin{proof}
    The pullback of $\tilde F^* \tilde \pi^* \tilde m - p \tilde m$ to
    $X$ is $F_{X/S}^* \pi^* m - p m = 0$, which implies $\tilde F^*
    \tilde \pi^* \tilde m - p \tilde m \in \scrM_{\tilde X}^*$ since the
    inclusion $X \to \tilde X$ is exact, and $\alpha_{\tilde X}(\tilde F^*
    \tilde \pi^* \tilde m - p \tilde m)$ pulls back to 1 in $X$ and is
    therefore in $1 + [p] \scrO_X$.  Now for $\tilde F_1, \tilde F_2 \in
    \calL_{\calXS}$, we have
    \[ \alpha_{\tilde X}(\tilde F_2^* \tilde \pi^* \tilde m - \tilde F_1^*
    \tilde \pi^* \tilde m) = 1 + [p] \langle \dlog(\pi^* m), \tilde F_2 - \tilde
    F_1 \rangle. \]
    Multiplying this by the equation $\alpha_{\tilde X}(\tilde F_1^* \tilde
    \pi^* \tilde m - p \tilde m) = 1 + [p] \beta_{\tilde \pi}(\tilde m)(\tilde
    F_1)$, we get
    \[ \alpha_{\tilde X}(\tilde F_2^* \tilde \pi^* \tilde m - \tilde F_1^*
    \tilde \pi^* \tilde m) = 1 + [p] (\beta_{\tilde \pi}(\tilde m)(\tilde F_1)
    + \langle \dlog(\pi^* m), \tilde F_2 - \tilde F_1 \rangle) = 1 + [p] \beta_{\tilde
      \pi}(\tilde m)(\tilde F_2). \] Hence $\beta_{\tilde \pi}(\tilde
    m)(\tilde F_2) - \beta_{\tilde \pi}(\tilde m)(\tilde F_1) = \langle \pi^*(\dlog
    m), \tilde F_2 - \tilde F_1 \rangle$; this proves that $\beta{\tilde
      \pi}(\tilde m) \in \calE_{\calXS}$ and also shows the commutativity
    of the right square in (iii).

    Now multiplying both sides of the equation defining $\beta_{\tilde
      \pi}(\tilde m)$ by $\alpha(p \tilde m)$, we get
    \[ \alpha_{\tilde X}(\tilde m)^p + [p] \alpha(m)^p \beta_{\tilde \pi}(\tilde F) =
    \alpha_{\tilde X}(\tilde F^* \tilde \pi^* \tilde m). \] Subtracting
    $\alpha_{\tilde X}(\tilde m)^p$ from both sides gives (ii).  In
    particular, if $\alpha(\tilde m) = 1 + [p] a$, then $\alpha(\tilde m)^p =
    1$, while $\tilde F^* \tilde \pi^* \alpha(\tilde m) = 1 + [p] F_X^* a$,
    so (ii) reduces to $[p] \beta_{\tilde \pi}(\tilde m) = [p] F_X^* a$,
    which shows the commutativity of the left square in (iii).
    
    To prove (i), we first show $\beta_{\tilde \pi}$ is additive.  To see
    this, we multiply the equations defining $\beta_{\tilde \pi}(\tilde
    m_1)$ and $\beta_{\tilde \pi}(\tilde m_2)$ to get
    \[ 1 + [p] (\beta_{\tilde \pi}(\tilde m_1) + \beta_{\tilde \pi}(\tilde m_2))
    (\tilde F) =
    \alpha_{\tilde X}(\tilde F^* \tilde \pi^* (\tilde m_1 + \tilde m_2) -
    p (\tilde m_1 + \tilde m_2)) = 1 + [p] \beta_{\tilde \pi}(\tilde m_1 +
    \tilde m_2) (\tilde F). \]
    Combining this fact with (iii) easily proves (i).
    
    The formula for the $p$-curvature in (v) follows directly from
    (iii).  Also, it suffices to verify the formula for $\nabla
    \beta_{\tilde \pi}(\tilde m)$ locally, so we may assume we have a
    lifting $\tilde F$ of $F_{X/S}$.  Let $g_{\tilde m} := \beta_{\tilde
      \pi}(\tilde m) (\tilde F)$; then from (iii) we conclude $\nabla
    \beta_{\tilde \pi}(\tilde m) = 1 \otimes d g_{\tilde m} - 1 \otimes \zeta_{\tilde
      F}(\dlog(\pi^* m))$.  On the other hand,
    \begin{align*}
      \tilde F^* (\dlog (\tilde \pi^* \tilde m)) & = [p] \dlog(m) +
      \dlog(\tilde F^* \tilde \pi^* \tilde m - p \tilde m) \\
      & = [p] \dlog(m) + \frac{d(1 + [p] g_{\tilde m})}{1 + [p]
        g_{\tilde m}} = [p] \dlog(m) + [p] dg_{\tilde m}.
    \end{align*}
    Hence $\zeta_{\tilde F}(\dlog(\pi^* m)) = \dlog(m) + dg_{\tilde m}$, and
    $\nabla \beta_{\tilde \pi}(\tilde m) = -1 \otimes \dlog(m)$ as required.
  \end{proof}
\end{proposition}

In terms of the construction of $\calE_{\calXS}$ as a conormal sheaf,
we can characterize $\beta_{\tilde \pi}$ as follows: for $\tilde m \in
\scrM_{\tilde X}^{\gp}$, the pullback of $(-p \tilde m, \tilde \pi^*
\tilde m) \in \scrM_{\tilde X \times \tilde X'}^{\gp}$ to $\MXgp$ is 0, which
implies that it becomes a unit in the first infinitesimal neighborhood
$X_1$ of $X$ in $\tilde X \times_{\tilde S} \tilde X'$.  Then $\beta_{\tilde
  \pi}(\tilde m) = \alpha_{X_1}(-p \tilde m, \tilde \pi^* \tilde m) - 1$.

\subsection{Functoriality}

The geometric construction of $\calE_{\calXS}$ given above makes it
straightforward to check its functoriality.  First we treat the
general functoriality properties of $\bfT(\calL)$.

\begin{lemma}
  Suppose we have a morphism $h : X\to Y$ of ringed topoi.
  \begin{enumerate}
  \item Let $\scrT_1 \to \scrT_2$ be a map of locally free sheaves on
    $X$, and let $\Omega_i := \check \scrT_i$.  Then for $\calL_1$ a
    $\scrT_1$-torsor, let $\calL_2 := \calL_1 \land_{\scrT_1} \scrT_2$ be
    the $\scrT_2$-torsor induced from $\calL$.  (Note that this is the
    unique $\scrT_2$-torsor $\calL_2$ with a map $\calL_1 \to \calL_2$
    respecting the actions.)  Then letting $\calK_i$ be the
    $\scrO_X$-algebra representing $\calL_i$, we have a natural map
    $\calK_2 \to \calK_1$ respecting the filtration $N_\cdot$, and the
    restriction to $\calE_2 \to \calE_1$ fits into a commutative diagram
    \[
    \begin{CD}
      0 @>>> \scrO_X @>>> \calE_2 @>>> \Omega_2 @>>> 0 \\
      @. @VV \id V @VVV @VVV \\
      0 @>>> \scrO_X @>>> \calE_1 @>>> \Omega_1 @>>> 0.
    \end{CD}
    \]
  \item Let $\scrT$ be a locally free sheaf on $Y$, and let $\calL$ be
    an $\scrT$-torsor on $Y$.  Then we have a natural isomorphism $h^*
    \calK_{\calL} \simeq \calK_{h^* \calL}$ respecting the filtration
    $N_\cdot$, where $h^* \calL := h^{-1} \calL \land_{h^{-1} \scrT} h^*
    \scrT$ is the induced $h^* \scrT$-torsor.  Furthermore, the
    restriction $h^* \calE_{\calL} \simeq \calE_{h^* \calL}$ is an
    isomorphism of extensions of $h^* \Omega$ by $\scrO_X$.
  \end{enumerate}

  \begin{proof}
    For (i), the map $\calE_2 \to \calE_1$ is induced by composition
    with the map $\calL_1 \to \calL_2$.  Since the map $\calL_1 \to
    \calL_2$ respects the actions of $\scrT_1$ and $\scrT_2$, we see
    that the image of $\calE_2$ is in fact contained in $\calE_1$, and
    we also get the commutativity of the right square in the diagram
    above.  The commutativity of the left square is obvious.  This
    induces the map $\calK_2 \to \calK_1$.

    For (ii), we have a map from $h^{-1} \calL$ to the $h^*
    \scrT$-torsor of splittings of
    \[ 0 \to h^* \scrO_Y \to h^* \calE_{\calL} \to h^* \Omega \to 0 \]
    on the left, which for $a \in \calL$ sends $h^{-1} a \in h^{-1} \calL$
    to $h^* \epsilon_a$, where $\epsilon_a$ is evaluation at $a$.  This induces a
    morphism of $h^* \scrT$-torsors from $h^* \calL$ to the torsor of
    splittings.  However, the torsor of splittings determines an
    extension uniquely; to see this, given an extension $0 \to \scrO_X \to
    \calE \to \Omega \to 0$, the dual is an extension
    \[ 0 \to \scrT \to \check \calE \to \scrO_X \to 0. \]
    Then the torsor of splittings is just the inverse image of $1 \in
    \scrO_X$ in $\check \calE$; however, taking the inverse image of 1
    is exactly the canonical isomorphism $\Ext^1(\scrO_X, \scrT) \simeq
    H^1(X, \scrT)$.  Therefore, we get an isomorphism of extensions of
    $h^* \Omega$ by $\scrO_X$, $h^* \calE_{\calL} \simeq \calE_{h^* \calL}$.
    Again, this easily extends to an isomorphism $h^* \calK_{\calL} \simeq
    \calK_{h^* \calL}$.
  \end{proof}
\end{lemma}

By a morphism $\calXS \to \calYS$ we mean a map $h : X \to Y$ along with a
lifting $\tilde h' : \tilde X' \to \tilde Y'$ of $h'$.

\begin{proposition}
  A morphism $(h, \tilde h') : \calXS \to \calYS$, where $X$ and $Y$
  are smooth $S$-schemes, induces a horizontal morphism $\theta_{h, \tilde
    h'} : h^* \calK_{\calYS} \to \calK_{\calXS}$ respecting the
  filtration $N_\cdot$, such that the restriction to $h^* \calE_{\calYS} \to
  \calE_{\calXS}$ fits into a commutative diagram
  \[
  \begin{CD}
    0 @>>> h^* \scrO_Y @>>> h^* \calE_{\calYS} @>>> h^* F_{Y/S}^*
    \Omega^1_{Y'/S} @>>> 0 \\
    @. @VVV @VV \theta_{h, \tilde h'} V @VVV @. \\
    0 @>>> \scrO_X @>>> \calE_{\calXS} @>>> F_{X/S}^* \Omega^1_{X'/S} @>>>
    0.
  \end{CD}
  \]
  Here the map $h^* \scrO_Y \to \scrO_X$ is the standard isomorphism,
  and the map $h^* F_{Y/S}^* \Omega^1_{Y'/S} \simeq F_{X/S}^* h'^* \Omega^1_{Y'/S}
  \to F_{X/S}^* \Omega^1_{X'/S}$ is the pullback by $F_{X/S}^*$ of the
  natural map $h'^* \Omega^1_{Y'/S} \to \Omega^1_{X'/S}$.
  
  For fixed $h$, let $\tilde h_1', \tilde h_2'$ be two liftings of
  $h'$, and let $D : h'^* \Omega^1_{Y'/S} \to \scrO_{X'}$ express their
  difference.  Then $\theta_{h, \tilde h_2'} - \theta_{h, \tilde h_1'} : h^*
  \calE_{\calYS} \to \calE_{\calXS}$ is the composition of the
  projection $h^* \calE_{\calYS} \to h^* F_{Y/S}^* \Omega^1_{Y'/S}$, the map
  \[ F_{X/S}^*(D) : h^* F_{Y/S}^* \Omega^1_{Y'/S} \simeq F_{X/S}^* h'^* \Omega^1_{Y'/S} \to
  \scrO_X, \]
  and the inclusion $\scrO_X \to \calE_{\calXS}$.

  \begin{proof}
    Let $T_1, T_2$ be objects of $\Crys(X / S)$, $\Crys(Y / S)$,
    respectively, and let $g : T_1 \to T_2$ be a PD morphism extending
    $h : X \to Y$.  (Thus $g$ is a section of $h^{-1} T_2 (T_1)$.)
    Then as in the construction of the crystalline structure of
    $\calL_{\calXS}$, composing with liftings $\tilde g$ of $g$
    induces a morphism of $h^* F_{Y/S}^* \scrT_{Y'/S}$-torsors $g^*
    \calL_{\calYS, T_2} \to \calL_{\calYS, X, T_1}$, where
    $\calL_{\calYS, X, T_1}$ denotes the sheaf of liftings of $f_{T_2
      / S} \circ g : X\to Y'$.  This is clearly compatible with the
    transition maps for $\calL_{\calYS, X}$ and $\calL_{\calYS}$, so
    we get a morphism of sheaves $h^* \calL_{\calYS} \to
    \calL_{\calYS, X}$ on the crystalline site $\Crys(X / S)$.  On the
    other hand, composing a lifting $\tilde T_1 \to \tilde X'$ with
    $\tilde h'$ gives a morphism of torsors $\calL_{\calXS} \to
    \calL_{\calYS, X}$.  Thus, overall, $h$ and $\tilde h'$ induce a
    morphism of torsors $\calL_{\calXS} \land_{\scrT_1} \scrT_2 \to
    h^* \calL_{\calYS}$, which induces the desired map $\theta_{h,
      \tilde h'}$.

    Now given $\tilde h_1', \tilde h_2' : \tilde X' \to \tilde Y'$ with
    difference $D : h'^* \Omega^1_{Y'/S} \to \scrO_{X'}$, for $\tilde F \in
    \calL_{\calXS}$, we have $\tilde h_2' \circ \tilde F - \tilde h_1' \circ
    \tilde F = F_{X/S}^* (D) : F_{X/S}^* h'^* \Omega^1_{Y'/S} \to \scrO_X$.
    Thus, for $\phi \in \calE_{\calYS}$ with image $\omega \in F_{Y/S}^*
    \Omega^1_{Y'/S}$,
    \[ (\theta_{h, \tilde h_2'} \phi - \theta_{h, \tilde h_1'} \phi)(\tilde F) =
    (\theta_h \phi)(\tilde h_2' \circ \tilde F) - (\theta_h \phi)(\tilde h_1' \circ \tilde F)
    = \langle h^* \omega, \tilde h_2' \circ \tilde F - \tilde h_1' \circ \tilde F \rangle, \]
    which implies that $\theta_{h, \tilde h_2'} \phi - \theta_{h, \tilde h_1'} \phi$
    is the constant $(F_{X/S}^* D)(\omega)$.
  \end{proof}
\end{proposition}

\begin{corollary}
  A morphism $(h, \tilde h') : \calXS \to \calYS$, where $X$ and $Y$
  are smooth $S$-schemes, induces an exact sequence
  \[ h^* \calE_{\calYS} \to \calE_{\calXS} \to F_{X/S}^* \Omega^1_{X'/Y'} \to
  0. \] If $h$ is smooth, this sequence is short exact and locally
  split.  Similarly, if $h$ is a closed immersion and $F_{X/S}$ is
  flat, there is an exact and locally split sequence
  \[ 0 \to F_{X/S}^* \check\scrN_{Y' / X'} \to h^* \calE_{\calYS} \to
  \calE_{\calXS} \to 0. \]
\end{corollary}


\section{The Cartier Transform}

\subsection{PD Higgs Fields and $p$-crystals}

Suppose we are given a lifting $\tilde X' / \tilde S$ of $X' / S$; we
shall see that $\calXS := (X/S, \tilde X' / \tilde S)$ determines a
splitting of the Azumaya algebra $\tildeDXS$ over $\calB_{X/S}
\otimes_{\scrO_{X'}} \scrO_{\calG}$, where $\calG$ is the nilpotent divided
power envelope of the zero section of the cotangent bundle of $X'/S$,
so that $\scrO_{\calG} = \hat \Gamma_\cdot \scrT_{X'/S}$.  We thus let
$\scrO_{\calG}^{\calB} := \calB_{X/S} \otimes_{\scrO_{X'}} \scrO_{\calG} \simeq
\tilde\frakZ \otimes_{S^\cdot \scrT_{X'/S}} \scrO_{\calG}$, where
$\tilde\frakZ$ is the center of $\tildeDXS$.  We also let
$\scrO_{\calG}^{\calA} := \AXgp \otimes_{\scrO_{X'}} \scrO_{\calG} \simeq \AXgp
\otimes_{\calB_{X/S}} \scrO_{\calG}^{\calB}$, and $\tildeDXS^\gamma :=
\tildeDXS \otimes_{\tilde\frakZ} \scrO_{\calG}^{\calB} \simeq \tildeDXS
\otimes_{S^\cdot \scrT_{X'/S}} \scrO_{\calG}$.

Fix an $\MXbargp$-set $\calJ$.  We denote by $HIG_{PD}^{\calB}(X'/S)$
the category of $\calJ$-indexed $\scrO_{\calG}^{\calB}$-modules, or
equivalently the category of $\calJ$-indexed $\calB_{X/S}$-modules
$E'$ equipped with a $\calB_{X/S}$-linear $\calG$-Higgs field
\[ \theta : \scrO_{\calG}
\to \scrEnd_{\calB_{X/S}}(E'). \] Similarly, we denote by
$MIC_{PD}^{\calA}(X/S)$ the category of $\calJ$-indexed
$\tildeDXS^\gamma$-modules.  An object of $MIC_{PD}^{\calA}(X/S)$ is
equivalent to a $\calJ$-indexed $\AXgp$-module $E$ with an integrable
and admissible connection $\nabla$, along with an $\AXgp$-linear
$\calG$-Higgs field
\[ \theta : \scrO_{\calG} \to F_{X/S*} \scrEnd_{\AXgp}(E) \]
extending the $F$-Higgs field
\[ \psi : S^\cdot \scrT_{X'/S} \to F_{X/S*} \scrEnd_{\AXgp}(E) \]
given by the $p$-curvature of $\nabla$.  We also let $MIC_{PD
  \cdot}^{\calA}(X/S)$ and $HIG_{PD \cdot}^{\calB}(X'/S)$ be the full
subcategories of locally nilpotent objects of $MIC_{PD}^{\calA}(X/S)$
and $HIG_{PD}^{\calB}(X'/S)$, respectively, that is, the objects such
that every element of $E$ or $E'$ is annihilated by a sufficiently
high power of $\Gamma_+ \scrT_{X'/S}$.

For example, let $MIC_\ell^\calA(X/S)$ be the category of
$\AXgp$-modules with integrable, admissible connection $\nabla$ such that
the $p$-curvature $\psi$ is nilpotent of level $\leq \ell$; then for $\ell
< p$ we may consider $MIC_\ell^\calA(X/S)$ as a full subcategory of
$MIC_{PD}^\calA(X/S)$ by letting $\Gamma_{\geq p} \scrT_{X'/S}$ act as zero.
Now the convolution product coming from the group scheme structure of
$\calG$ allows us to make $MIC_{PD}^{\calA}(X/S)$ and
$HIG_{PD}^{\calB}(X'/S)$ into tensor categories; in particular, the
total $\calG$-Higgs field on a tensor product is given by
\[ \theta_{(D')^{[n]}} = \sum_{i+j=n} \theta_{(D')^{[i]}} \otimes \theta_{(D')^{[j]}}. \]
Note that, in particular, $\theta_{(D')^{[p]}}$ can be nonzero on $E_1 \otimes
E_2$ even if $E_1, E_2 \in MIC_\ell^\calA(X/S)$ for $\ell < p$.

We shall see that the Azumaya algebra $\tildeDXS$ splits over
$\scrO_{\calG}^{\calB}$, which gives a Riemann-Hilbert correspondence
between $MIC_{PD \cdot}^{\calA}(X / S)$ and $HIG_{PD \cdot}^{\calB}(X' / S)$.
We shall also use this splitting to get a Riemann-Hilbert
correspondence for certain $\scrO_X$-modules with connection.  For this,
let $\DXS^\gamma := \DXS \otimes_{S^\cdot \scrT_{X/S}} \scrO_{\calG}$, and let
$MIC_{PD}(X/S)$ and $HIG_{PD}(X'/S)$ be the categories of
$\DXS^\gamma$- and $\scrO_{\calG}$-modules, respectively.  We also let
$MIC_{PD \cdot}(X/S)$ and $HIG_{PD \cdot}(X'/S)$ be the corresponding
full subcategories of locally nilpotent objects.

We can provide interpretations of the categories $MIC_{PD \cdot}(X/S)$
and $MIC_{PD \cdot}^\calA(X/S)$ which are crystalline in nature.  Our
motivation for the following theory is the fact that the ring of
\diffop{}s on $X$ over $S$ is very similar to $\Gamma_\cdot \scrT_{X/S}$,
except for the fact that it is noncommutative.  Thus, in our theory we
will allow $p$th divided powers of elements of the ideal of the
diagonal in $X ×_S X$, in order to preserve the action of the image of
$c' : F_{X/S}^* \scrT_{X'/S} \to \DXS$, but no higher divided powers.

In particular, define a $p$-PD ring $A$ to be an
algebra over a given ring $B$ of characteristic $p$, with an ideal $I$
of $B$ such that $I^{(p)} := \{ x^p : x \in I \} = (0)$, and a map
$\gamma_p : I \to A$ satisfying
\begin{align*}
  \gamma_p(x + y) & = \gamma_p(x) + \gamma_p(y) + \sum_{i=1}^{p-1} \frac{x^{p-i} y^i}
  {(p-i)! i!}, \\
  \gamma_p(\lambda x) & = \lambda^p \gamma_p(x).
\end{align*}
As usual, we write $x^{[p]} := \gamma_p(x)$.  Similarly to the case of
full PD ideals, we can show that there is a left adjoint $\Gamma_{R, (p)}$
to the functor $(A, B, I, \gamma) \mapsto I$ from the category of $p$-PD rings
over $R$ to the category of $R$-modules.  Similarly, given a closed
immersion $Y \to X$, there is a right universal $p$-PD scheme $D_{Y,
  (p)}$ with an exact closed immersion $Y \to D_{Y, 0, (p)}$ defined by
the $p$-PD ideal of $D_{Y, (p)}$ and a map $D_{Y, 0, (p)} \to X$; we
call $D_{Y, (p)}$ the $p$-PD envelope of $Y$.  (Here $D_{Y, (p)}$
corresponds to $A$, and $D_{Y, 0, (p)}$ corresponds to $B$.)  We now
define $\Crys_{(p)}(X/S)$ as the category of tuples $(U, T, V, i, f,
\gamma_p)$, where $U$ is an open subscheme of $X$, $f : V\to T$ is a strict
affine map, and $i : U\to T$ is an exact closed immersion defined by a
$p$-PD ideal $(I_T, \gamma_p : I_T \to f_* \scrO_V)$.  For example, for each
$n$ we may define $D_{(p)}(n) \in \Crys_{(p)}(X/S)$ to be the $p$-PD
envelope of the diagonal embedding $\Delta : X \to X^{n+1}$.  Also, for
$(U, T, i, \gamma) \in \Crys(X/S)$, we have an object $(U, T, T, i, \id,
\gamma_p) \in \Crys_{(p)}(X/S)$, so that we may regard $\Crys(X/S)$ as a
(full) subcategory of $\Crys_{(p)}(X/S)$.  We let $\scrO_{X/S}$ be
the sheaf of rings on $\Crys_{(p)}(X/S)$ which to $(U, T, V, i, f,
\gamma_p)$ associates $\scrO_V$.

Locally on $X$, if we have a logarithmic system of coordinates $m_1,
\ldots, m_r$, then as usual define $\eta_i := \alpha_{D_{0, (p)}(1)} (-m_i, m_i)
- 1$.  Then a basis for $\divpow_{(p)}(1) := \scrO_{D_{(p)}(1)}$ as an
$\scrO_X$-module is $\{ \eta^I (\eta^{[p]})^J : I \in M := \{ 0, \ldots, p-1
\}^r, J \in \nats^r \}$, and $\scrO_{D_{0, (p)}(1)}$ has basis $\{ \eta^I
: I \in M \}$.  Also, we have a natural $p$-PD morphism $\delta :
\divpow_{(p)}(1) \to \divpow_{(p)}(1) \otimes_{\scrO_X} \divpow_{(p)}(1)$
such that $\delta(1 + \eta_i) = (1 + \eta_i) \otimes (1 + \eta_i)$, which allows us
to define composition on the ring $\hat \calD_{(p), X/S} :=
\scrHom_{\scrO_X} (\divpow_{(p)}(1), \scrO_X)$ of $p$-HPD differential
operators $\scrO_X \to \scrO_X$.

\begin{proposition}
  The following data are equivalent:
  \begin{enumerate}
  \item An object $E$ of $MIC_{PD \cdot}(X/S)$.
  \item A ring homomorphism $\hat \calD_{(p), X/S} \to
    \scrHom_{\scrO_X}(\divpow_{(p)}(1) \otimes E, E)$.
  \item A $p$-HPD stratification $\varepsilon : \divpow_{(p)}(1) \otimes E \to E \otimes
    \divpow_{(p)}(1)$ satisfying the natural cocycle condition.
  \item A crystal $E$ of $\scrO_{X/S}$-modules on the site
    $\Crys_{(p)}(X/S)$.
  \end{enumerate}

  \begin{proof}
    The proof follows as in the case of regular crystals once we
    establish an isomorphism $D^\gamma_{X/S} \simeq \hat \calD_{(p), X/S}$.
    To define this isomorphism, first note that we have a natural map
    $D(1) \to D_{(p)}(1)$ since $D(1)$ is a member of
    $\Crys_{(p)}(X/S)$; this induces a ring homomorphism $\DXS \to \hat
    \calD_{(p), X/S}$ by pushforward.  Also, if $I$ is the $p$-PD
    ideal of $\divpow_{0, (p)}(1)$, then we have a natural isomorphism
    $S^\cdot F_{X/S}^* \Omega^1_{X'/S} \overset{\sim}{\to} \divpow_{(p)}(1) / I
    \divpow_{(p)}(1)$ defined by $1 \otimes \pi^* \omega \mapsto \omega^{[p]}$.  This
    gives a natural isomorphism $\hat \Gamma_\cdot F_{X/S}^* \scrT_{X'/S} \simeq
    (I \divpow_{(p)}(1))^\perp$.  Furthermore, since $\delta(\omega^{[p]}) \equiv
    \omega^{[p]} \otimes 1 + 1 \otimes \omega^{[p]} \pmod{I \divpow_{(p)}(1)}$ for $\omega
    \in \Omega^1_{X/S}$, $\delta$ becomes the comultiplication on $S^\cdot
    F_{X/S}^* \Omega^1_{X'/S}$, and since $1 \otimes f - f \otimes 1 \in I$ for $f \in
    \scrO_X$, composition is dual to this comultiplication map.  This
    implies that the map $\hat \Gamma_\cdot F_{X/S}^* \scrT_{X'/S} \simeq (I
    \divpow_{(p)}(1))^\perp$ is in fact a ring homomorphism.  It is easy
    to see the two maps agree on $S^\cdot F_{X/S}^* \scrT_{X'/S}$, and a
    local calculation below will show that the image of $\DXS$
    commutes with $(I \divpow_{(p)}(1))^\perp$, so we get a map
    $D_{X/S}^\gamma \to \hat \calD_{(p), X/S}$.
    
    To prove this gives an isomorphism, it suffices to work locally,
    so assume we have a system of logarithmic coordinates $m_1, \ldots,
    m_r$.  Define $\zeta^{\langle I + pJ \rangle} := \zeta^I (\zeta^{[p]})^J / I!$ for $I
    \in M$ and $J \in \nats^r$, and let $\{ D_{\langle N \rangle} \}$ be the
    dual basis of $\hat \calD_{(p), X/S}$ to $\{ \zeta^{\langle N \rangle} \}$.
    Now $\{ D_I \otimes (D')^{[J]} : I\in M, J \in \nats^r \}$ forms a
    completion basis for $D^\gamma_{X/S}$; we claim that our map sends
    $D_I \otimes (D')^{[J]}$ to $D_{\langle I + pJ \rangle}$, which will complete the
    proof.  It is easy to check that the map $\DXS \to \hat \calD_{(p),
      X/S}$ sends $D_I$ to $D_{\langle I \rangle}$ for $I \in M$, and the map
    $\hat \Gamma_\cdot F_{X/S}^* \scrT_{X'/S} \to \hat \calD_{(p), X/S}$ sends
    $(D')^{[J]}$ to $D_{\langle pJ \rangle}$; thus, it will suffice to show
    $D_{\langle I \rangle} \circ D_{\langle pJ \rangle} = D_{\langle pJ \rangle} \circ D_{\langle I \rangle} = D_{\langle
      I + pJ \rangle}$.  However, for $A \in M$, we have $\delta(\eta^{\langle A \rangle}) =
    (\eta \otimes 1 + 1 \otimes \eta + \eta \otimes \eta)^A / A!$.  If we reduce modulo $1 \otimes
    I \divpow_{(p)}(1)$, the only nonzero term in the expansion of
    this expression is $(\eta \otimes 1)^A / A! = \eta^{\langle A \rangle} \otimes 1$.  Thus,
    for $A \in M$ and $B \in \nats^r$,
    \begin{align*}
      \delta(\eta^{\langle A + pB \rangle}) & \equiv (\eta^{\langle A \rangle} \otimes 1) (\eta^{[p]} \otimes 1 + 1
      \otimes \eta^{[p]})^B \\
      & \equiv \sum_{D \leq B} \binom{B}{D} \eta^{\langle A + p
        (B - D) \rangle} \otimes \eta^{\langle p D \rangle} \pmod{1 \otimes I \divpow_{(p)}(1)}.
    \end{align*}
    Applying $\id \otimes D_{\langle pJ \rangle}$ gives $\eta^{\langle A + p (B - J) \rangle}$,
    and applying $D_{\langle I \rangle}$ to this gives $\delta_{AI} \delta_{BJ} = \delta_{A
      + pB, I + pJ}$; this shows that $D_{\langle I \rangle} \circ D_{\langle pJ \rangle} =
    D_{\langle I + pJ \rangle}$.  Similarly, reducing $\delta(\eta^{\langle A \rangle})$ modulo
    $I \divpow_{(p)}(1) \otimes 1$ gives $1 \otimes \eta^{\langle A \rangle}$, so reducing
    $\delta(\eta^{\langle A + pB \rangle})$ gives $\sum_{D\leq B} \binom{B}{D} \eta^{\langle p (B
      - D) \rangle} \otimes \eta^{\langle A + pD \rangle}$.  Hence $D_{\langle pJ \rangle} \circ D_{\langle I
      \rangle}(\eta^{\langle A + pB \rangle}) = \delta_{AI} \delta_{BJ}$, and $D_{\langle pJ \rangle} \circ
    D_{\langle I \rangle} = D_{\langle I + pJ \rangle}$.
  \end{proof}
\end{proposition}

Note that in particular the constructions from the previous chapter
easily extend to give crystals $\calE_{\calXS}$ and $\calK_{\calXS}$
on $\Crys_{(p)}(X/S)$, so they have canonical structures of
$D_{X/S}^\gamma$-module.  In particular, for $(U, \tilde T, \tilde V) \in
\Crys_{(p)}(X / \tilde S)$, $\calE_{\calXS, \tilde V}$ is the conormal
sheaf of $V \overset{\Gamma}{\to} V ×_S X' \to \tilde V ×_{\tilde S} \tilde
X'$, where $\Gamma$ is the graph of $f_{V/S}$, the composition of the map
$V \to T$ with $f_{T/S}$.  Similarly, $MIC_{PD \cdot}^{\calA}(X/S)$ is
equivalent to the category of $p$-crystals of $\AXgp$-modules.

\subsection{The Global Cartier Transform}

Given a lifting $\tilde X' \to \tilde S$ of $X'$ over $S$, we locally
have isomorphisms $\calK_{\calXS} \simeq S^\cdot F_{X/S}^* \Omega^1_{X'/S}$, and
thus $\check \calK_{\calXS} = \scrHom_{\scrO_X} (\calK_{\calXS},
\scrO_X) \simeq \hat \Gamma_\cdot F_{X/S}^* \scrT_{X'/S} \simeq F_{X/S}^*
\scrO_{\calG}$.  Now the total connection on $\AXgp \otimes_{\scrO_X}
\check \calK_{\calXS}$ extends naturally to a $\tildeDXS^\gamma$-module
structure, and $\AXgp \otimes_{\scrO_X} \check \calK_{\calXS}$ is a locally
free $\scrO_{\calG}^{\calA}$-module of rank 1 (with basis $1 \otimes
\xi_{\tilde F}$ for $\tilde F$ a lifting of $F_{X/S}$).  However,
$\scrO_{\calG}^{\calA} \simeq \scrO_{\calG}^{\calB} \otimes_{\calB_{X/S}}
\AXgp$ is locally free of rank $p^d$ over $\scrO_{\calG}^{\calB}$.
Hence $\AXgp \otimes_{\scrO_X} \check \calK_{\calXS}$ is a splitting module
for $\tildeDXS$ over $\scrO_{\calG}^{\calB}$, immediately giving the
following equivalence.

\begin{theorem}
  Let $\check \calK_{\calXS}^{\calA} := \AXgp \otimes_{\scrO_X} \check
  \calK_{\calXS}$, and define functors
  \begin{align*}
    C_{\calXS} : MIC_{PD}^{\calA}(X/S) \to HIG_{PD}^{\calB}(X'/S), & ~
    E \mapsto \iota_* \scrHom_{\tildeDXS^\gamma} (\check \calK_{\calXS}^{\calA}, E), \\
    C_{\calXS}^{-1} : HIG_{PD}^{\calB}(X'/S) \to MIC_{PD}^{\calA}(X/S),
    & ~ E' \mapsto \check \calK_{\calXS}^{\calA} \otimes_{\scrO_\calG^\calB}
    \iota_* E'.
  \end{align*}
  (Here $\iota : \scrO_{\calG} \to \scrO_{\calG}$ is the map corresponding
  to the inverse map $\calG \to \calG$, so that $\iota_* E'$ is $E'$ with a
  sign change of the Higgs field.)
  Then $C_{\calXS}$ and $C_{\calXS}^{-1}$ are quasi-inverse
  equivalences of categories.
\end{theorem}
\begin{remark}
We will see later that the sign changes in the above definitions are
necessary in order for the isomorphism of de Rham and Higgs complexes
to be compatible with the standard Cartier isomorphism.
\end{remark}

In the case that $\calJ = \MXbargp$ with the standard action, we have
that the category of $\calJ$-indexed $\AXgp$-modules is equivalent to
the category of $\scrO_X$-modules; indeed, it is easy to see that for
any $\AXgp$-module $E$, the natural map $E_0 \otimes_{\scrO_X} \AXgp \to E$ is
an isomorphism.  Similarly, a connection on $E$ is admissible if and
only if this map is horizontal.  We thus get:

\begin{theorem}
  Given a lifting $\calXS$ of $X'$ over $S$, we have an equivalence of
  categories
  \[ C_{\calXS} : MIC_{PD}(X / S) \to HIG_{PD}^{\calB} (X' / S). \]
\end{theorem}

However, $\check \calK_{\calXS}^\calA$ has the unfortunate property
that it is not locally nilpotent.  We thus reformulate the functors
described above in terms of $\calK_{\calXS}^{\calA} := \AXgp
\otimes_{\scrO_X} \calK_{\calXS}$ instead.  For $E$ an object of $MIC_{PD
  \cdot}^{\calA}(X / S)$, let $E^\theta := \scrHom_{\scrO_{\calG}^{\calA}}
(\AXgp, E)$ be the subsheaf of $E$ annihilated by $\Gamma_+
\scrT_{X'/S}$, and let $E^{\nabla, \gamma} := \scrHom_{\tildeDXS^\gamma} (\AXgp,
E)$ be the subsheaf of $E^\theta$ annihilated by $\nabla$.

\begin{theorem}
  \label{thm:cart-rh}
  \begin{enumerate}
  \item For $E$ an object of $MIC_{PD \cdot}^{\calA}(X / S)$, there is a
    natural isomorphism
    \[ C_{\calXS}(E) \simeq (\calK_{\calXS}^{\calA} \otimes_{\AXgp} E)^{\nabla, \gamma}, \]
    where $\scrO_{\calG}^{\calB}$ acts on the right hand side via the
    action on $\calK_{\calXS}^{\calA}$.
  \item For $E'$ an object of $HIG_{PD \cdot}^{\calB}(X' / S)$, there is
    a natural isomorphism
    \[ C_{\calXS}^{-1}(E') \simeq (\calK_{\calXS}^{\calA} \otimes_{\calB_{X/S}}
    E')^\theta, \]
    where $\tildeDXS^\gamma$ acts on the right hand side via the action on
    $\calK_{\calXS}^{\calA}$.
  \end{enumerate}
\end{theorem}

\begin{proof}
  First, since $E$ is locally nilpotent, and $\check
  \calK_{\calXS}^{\calA}$ is locally free of rank 1 as an
  $\scrO_{\calG}^{\calA}$-module, we have an isomorphism
  \[ \scrHom_{\scrO_{\calG}^{\calA}} (\check \calK_{\calXS}^{\calA},
  E) \simeq (\calK_{\calXS}^{\calA} \otimes_{\AXgp} E)^\theta. \] Under this
  isomorphism, the submodule $C_{\calXS}(E) = \scrHom_{\tildeDXS^\gamma}
  (\check \calK_{\calXS}^{\calA}, E)$ corresponds to
  $(\calK_{\calXS}^{\calA} \otimes_{\AXgp} E)^{\nabla, \gamma}$ with
  $\scrO_{\calG}^{\calB}$ acting on $E$.  Changing the sign is the
  same as letting $\scrO_{\calG}^{\calB}$ act on
  $\calK_{\calXS}^{\calA}$ instead, proving (i).

  To prove (ii), we begin with a lemma.

\begin{lemma}
  There is a canonical isomorphism
  \[
  \begin{matrix}
    \check\calK_{\calXS} \otimes_{F^* \scrO_{\calG}} \iota_*
    \check\calK_{\calXS} \simeq F_{X/S}^* \scrO_{\calG}.
  \end{matrix}
  \]

  \begin{proof}
    Suppose we have a lifting $\tilde F$ of $F_{X/S}$.  Then
    $\xi_{\tilde F} \otimes \xi_{\tilde F}$ forms a basis for $\check
    \calK_{\calXS} \otimes_{F^* \scrO_\calG} \iota_* \check\calK_{\calXS}$ as
    an $F_{X/S}^* \scrO_{\calG}$-module.  We claim that, in fact, this
    basis element is independent of the choice of $\tilde F$, so the
    local bases glue to a global basis of $\check \calK_{\calXS}
    \otimes_{F^* \scrO_{\calG}} \iota_* \check \calK_{\calXS}$.  Thus, suppose
    we have two liftings $\tilde F_1$ and $\tilde F_2$, and let $h :
    F_{X/S}^* \Omega^1_{X'/S} \to \scrO_X$ express the difference between
    them.  Also, suppose we have a system of logarithmic coordinates
    $m_1', \ldots, m_r' \in \scrM_{X'}^{\gp}$, and let $D^{[N]}$ ($N \in
    \nats^r$) be the corresponding basis of $F_{X/S}^* \scrO_{\calG} =
    F_{X/S}^* \hat \Gamma_\cdot \scrT_{X'/S}$.  Then for $N \in \nats^r$,
    \[ \prod_{i=1}^r \sigma_{\tilde F_2}(\dlog(m_i'))^{N_i} = \prod_{i=1}^r
    (\sigma_{\tilde F_1}(\dlog(m_i')) + h(\dlog(m_i')))^{N_i}. \] Now
    $\xi_{\tilde F_1}$ maps this element of $\calK_{\calXS}$ to $h^N :=
    \prod_{i=1}^r h(\dlog(m_i'))^{N_i}$, while $\theta_{D^{[I]}} \xi_{\tilde
      F_2}$ maps it to $(-1)^{|I|} \delta_{IN}$.  In other words,
    \[ \xi_{\tilde F_1} = \sum_{N \in \nats^r} (-1)^{|N|} h^N \theta_{D^{[N]}}
    \xi_{\tilde F_2}. \] Thus, defining $g := \sum_N (-1)^{|N|} h^N
    D^{[N]} \in F_{X/S}^* \scrO_{\calG}$, we see that $\xi_{\tilde F_1}
    = \theta_g \xi_{\tilde F_2}$; also, plugging in $-h$ in place of $h$,
    we get $\xi_{\tilde F_2} = \theta_{\iota(g)} \xi_{\tilde F_1}$.  Therefore,
    \[ \xi_{\tilde F_1} \otimes \xi_{\tilde F_1} = \theta_g \xi_{\tilde F_2} \otimes
    \xi_{\tilde F_1} = \xi_{\tilde F_2} \otimes \theta_{\iota(g)} \xi_{\tilde F_1} =
    \xi_{\tilde F_2} \otimes \xi_{\tilde F_2}. \]
  \end{proof}
\end{lemma}

  Now for $E'$ an object of $HIG_{PD \cdot}^{\calB} (X'/S)$, let $E_0 := E'
  \otimes_{\calB_{X/S}} \AXgp$.  Using the above result, we have $\check
  \calK_{\calXS}^{\calA} \otimes_{\scrO_{\calG}^{\calA}} \iota_* \check
  \calK_{\calXS}^{\calA} \simeq \scrO_{\calG}^{\calA}$.  Hence
  \begin{align*}
    \check \calK_{\calXS}^{\calA} \otimes_{\scrO_{\calG}^{\calB}} \iota_* E' & \simeq
    \check \calK_{\calXS}^{\calA} \otimes_{\scrO_{\calG}^{\calA}} \iota_* E_0 \simeq
    \scrHom_{\scrO_{\calG}^{\calA}} (\iota_* \check
      \calK_{\calXS}^{\calA}, \iota_* E_0)
      \simeq \scrHom_{\scrO_{\calG}^{\calA}}
      (\check \calK_{\calXS}^{\calA}, E_0) \\
      & \simeq
      (\calK_{\calXS}^{\calA} \otimes_{\AXgp} E_0)^\theta
      \simeq (\calK_{\calXS}^{\calA} \otimes_{\calB_{X/S}} E')^\theta.
    \end{align*}
\end{proof}

As an application, we can use this to calculate
$C_{\calXS}(\calK_{\calXS}^{\calA})$.  First, we claim there is a
canonical isomorphism
\[ (\calK_{\calXS} \otimes_{\scrO_X} \calK_{\calXS})^\theta \simeq S^\cdot F_{X/S}^*
\Omega^1_{X'/S} \simeq F_{X/S}^*(\iota_* \check \scrO_{\calG}). \]
To construct this isomorphism, we begin with the cosubtraction map
$S^\cdot F_{X/S}^* \Omega^1_{X'/S} \to \calK_{\calXS} \otimes_{\scrO_X}
\calK_{\calXS}$.  Then locally on $X$, if we have a lifting $\tilde F$
of $F_{X/S}$, the restriction to $F_{X/S}^* \Omega^1_{X'/S} \to
\calE_{\calXS} \otimes_{\scrO_X} \calE_{\calXS}$ is equal to $1 \otimes \sigma_{\tilde
  F} - \sigma_{\tilde F} \otimes 1$, which factors through $(\calE_{\calXS}
\otimes_{\scrO_X} \calE_{\calXS})^\theta$.
Since the total $\calG$-Higgs field on $\calK_{\calXS}
\otimes_{\scrO_X} \calK_{\calXS}$ satisfies the Leibniz rule, this implies
that the cosubtraction map factors through $(\calK_{\calXS}
\otimes_{\scrO_X} \calK_{\calXS})^\theta$.

We may check the map $S^\cdot F_{X/S}^* \Omega^1_{X'/S} \to (\calK_{\calXS}
\otimes_{\scrO_X} \calK_{\calXS})^\theta$ is an isomorphism
locally, so we may assume there is a lifting $\tilde F$ of $F_{X/S}$.
Then $(\calK_{\calXS} \otimes_{\scrO_X} \calK_{\calXS})^\theta \simeq \scrHom_{F^*
  \scrO_{\calG}}(\check \calK_{\calXS}, \calK_{\calXS})$, and
$\xi_{\tilde F}$ forms a basis for $\check \calK_{\calXS}$ as an $F^*
\scrO_{\calG}$-module.  Thus, $\xi_{\tilde F} \otimes \id$ gives an
isomorphism $(\calK_{\calXS} \otimes_{\scrO_X} \calK_{\calXS})^\theta \simeq
\calK_{\calXS}$, and the composition with the above map is just the
isomorphism $S^\cdot F_{X/S}^* \Omega^1_{X'/S} \overset{\sim}{\to}
\calK_{\calXS}$ induced by $\tilde F$.  In fact, this shows that the
canonical structure of $\calG$-Higgs field on $F_{X/S}^*(\iota_*
\check\scrO_{\calG})$ corresponds to the action of $\scrO_\calG$ on
the second $\calK_{\calXS}$.  Note also that since $\nabla_{tot}(1 \otimes
\sigma_{\tilde F}(1 \otimes \omega) - \sigma_{\tilde F}(1 \otimes \omega) \otimes 1) = -(1 \otimes 1) \otimes
\zeta_{\tilde F}(\omega) + (1 \otimes 1) \otimes \zeta_{\tilde F}(\omega) = 0$, the total
connection on $(\calK_{\calXS} \otimes_{\scrO_X} \calK_{\calXS})^\theta$
corresponds to the Frobenius descent connection on $F_{X/S}^* (\iota_*
\check \scrO_{\calG})$.
From this it follows that there is a canonical isomorphism
$C_{\calXS}(\calK_{\calXS}^{\calA}) \simeq (\calK_{\calXS}^{\calA}
\otimes_{\AXgp} \calK_{\calXS}^{\calA})^\theta \simeq \iota_* \check
\scrO_{\calG} \otimes_{\scrO_{X'}} \calB_{X/S}$.

\subsection{Nilpotent Residue}

The goal of this section is to apply the preceding theory to get
information about objects of $MIC_{PD}(X / S)$ and of $HIG_{PD}(X' /
S)$.  To do this, note that for $E$ an object of $MIC_{PD}(X / S)$, $E
\otimes_{\scrO_X} \AXgp$ with the total connection is an object of
$MIC_{PD}^{\calA}(X / S)$; similarly, for $E'$ an object of
$HIG_{PD}(X' / S)$, $E' \otimes_{\scrO_{X'}} \calB_{X/S}$ is an object of
$HIG_{PD}^{\calB}(X' / S)$.  Our first task is to calculate what the
Cartier transform does to these objects.

\begin{proposition}
  \begin{enumerate}
  \item Let $E'$ be an object of $HIG_{PD}(X' / S)$.  Then there is a
    natural isomorphism
    \[ C_{\calXS}^{-1}(E' \otimes_{\scrO_{X'}} \calB_{X/S}) \simeq (\check
    \calK_{\calXS} \otimes_{\scrO_{\calG}} \iota_* E') \otimes_{\scrO_X} \AXgp. \]
  \item Let $E$ be an object of $MIC_{PD}(X / S)$.  Then there is a
    natural map
    \[ \iota_* \scrHom_{\DXS^\gamma}(\check \calK_{\calXS}, E) \otimes_{\scrO_{X'}}
    \calB_{X/S} \to C_{\calXS}(E \otimes_{\scrO_X} \AXgp). \] Furthermore,
    denoting $E_0 := \scrHom_{F^* \scrO_{\calG}} (\check
    \calK_{\calXS}, E)$ with the internal Hom connection, this map is
    injective, resp.~surjective, if and only if the natural map
    $E_0^\nabla \otimes_{\scrO_{X'}} \scrO_X \to E_0$ is.
  \end{enumerate}
\end{proposition}

\begin{example}
  To illustrate the second part, let $X$ be $\bbA^1_k$ minus one point
  with $k$ a field of characteristic $p$; that is, $X := \Spec(\nats \to
  k[t])$, where the map $\nats \to k[t]$ sends 1 to $t$, and $S := \Spec
  k$.  Consider the $\scrO_X$-module $\scrO_X \cdot e$ with the connection
  $\nabla$ such that $\nabla(e) = e \otimes \dlog(t)$.  Then since $\nabla$ has zero
  $p$-curvature, any $F_{X/S}^* \scrO_{\calG}$-linear map $\check
  \calK_{\calXS} \to E$ factors uniquely through the map $\check
  \calK_{\calXS} \to \scrO_X$ induced by the inclusion $\scrO_X \to
  \calK_{\calXS}$.  We thus get an isomorphism $E_0 \simeq E$, and we will
  see below that this isomorphism is horizontal.  Thus
  $\scrHom_{\DXS^\gamma}(\check \calK_{\calXS}, E) \simeq E^\nabla = \scrO_{X'} \cdot
  t^{p-1} e$.  Similarly, $C_{\calXS}(E \otimes_{\scrO_X} \AXgp) \simeq (E
  \otimes_{\scrO_X} \AXgp)^\nabla = \calB_{X/S} \cdot (e \otimes e_{-1})$.  Now the map
  described in the proof of (ii) below sends $t^{p-1} e$ to $t^{p-1} e
  \otimes 1 = (t^{p-1} e_1) \cdot (e \otimes e_{-1})$, where $t^{p-1} e_1 \in
  \calB_{X/S}$.  This map is clearly not surjective, and neither is
  the map $E_0^\nabla \otimes_{\scrO_{X'}} \scrO_X \to E_0$.
\end{example}

\begin{proof}
  For (i), we have
  \begin{align*}
    C_{\calXS}^{-1}(E' \otimes_{\scrO_{X'}} \calB_{X/S})
    & = \check \calK_{\calXS}^{\calA} \otimes_{\scrO_{\calG}^{\calB}} (\iota_* E'
    \otimes_{\scrO_{X'}} \calB_{X/S}) \simeq \check \calK_{\calXS}^{\calA}
    \otimes_{\scrO_{\calG}} \iota_* E' \\
    & \simeq (\check \calK_{\calXS} \otimes_{\scrO_{\calG}}
    \iota_* E') \otimes_{\scrO_X} \AXgp,
  \end{align*}
  as required.

  For (ii), we have $C_{\calXS}(E \otimes_{\scrO_X} \AXgp) \simeq
  \iota_* \scrHom_{\tildeDXS^\gamma} (\check \calK_{\calXS}^{\calA}, E
  \otimes_{\scrO_X} \AXgp)$.  However, since
  \[ \scrHom_{\scrO_{\calG}^{\calA}} (\check \calK_{\calXS}
  \otimes_{\scrO_X} \AXgp, E \otimes_{\scrO_X} \AXgp) \simeq \scrHom_{F^*
    \scrO_{\calG}} (\check \calK_{\calXS}, E) \otimes_{\scrO_X} \AXgp =
  E_0 \otimes_{\scrO_X} \AXgp, \]
  applying $\iota_*$ to the natural map $E_0^\nabla \otimes_{\scrO_{X'}}
  \calB_{X/S} \to (E_0 \otimes_{\scrO_X} \AXgp)^\nabla$ gives the desired map.

  We now have a commutative diagram
  \[
  \begin{CD}
    (E_0^\nabla \otimes_{\scrO_{X'}} \calB_{X/S}) \otimes_{\calB_{X/S}} \AXgp @>>>
    (E_0 \otimes_{\scrO_X} \AXgp)^\nabla \otimes_{\calB_{X/S}} \AXgp \\
    @VV \sim V @VV \sim V \\
    (E_0^\nabla \otimes_{\scrO_{X'}} \scrO_X) \otimes_{\scrO_X} \AXgp @>>> E_0
    \otimes_{\scrO_X} \AXgp.
  \end{CD}
  \]
  Here the map on the right is an isomorphism since $E_0$ has zero
  $p$-curvature by definition.  Now the top row is injective,
  resp.~surjective, if and only if $E_0^\nabla \otimes_{\scrO_{X'}} \calB_{X/S}
  \to (E_0 \otimes_{\scrO_X} \AXgp)^\nabla$ is, since $\AXgp$ is locally free
  over $\calB_{X/S}$.  Similarly, the bottom row is injective,
  resp.~surjective, if and only if $E_0^\nabla \otimes_{\scrO_{X'}} \scrO_X \to
  E_0$ is, since $(\AXgp)_s$ is invertible for each section $s \in
  \MXbargp$.
\end{proof}

We now need to study the map $E_0^\nabla \otimes_{\scrO_{X'}} \scrO_X \to E_0$;
however, the $p$-curvature of $\nabla$ on $E_0$ is zero since each
homomorphism in $E_0$ commutes with $\scrT_{X'/S} \subseteq \scrO_{\calG}$.
Recall that for a log scheme $X$ over $S$, we define the {\em residue
  sheaf} $\calR_{X/S} := \Omega^1_{X / X^*}$, where $X^*$ denotes the log
scheme with the same underlying scheme as $X$ but with log structure
$f^* \scrM_S$; we then have an exact sequence
\[ \Omega^1_{X^* / S} \to \Omega^1_{X/S} \to \calR_{X/S} \to 0. \] Also, for a
connection $\nabla : E \to E \otimes_{\scrO_X} \Omega^1_{X/S}$ we define the residue to
be the induced map $\rho : E \to E \otimes_{\scrO_X} \calR_{X/S}$; for each
derivation $D \in \scrT_{X/S}$ this induces an endomorphism $\rho_D$ on $E
/ D(\scrO_X) E$, where $D(\scrO_X)$ is the ideal of $\scrO_X$
generated by the image of $D$.  Now according to \cite{lorenzon}, the
map $E_0^\nabla \otimes_{\scrO_{X'}} \scrO_X \to E_0$ is surjective if the residue
of the connection on $E_0$ vanishes, and it is an isomorphism if in
addition $\Tor_1(E_0, \calR_{X/S}) = 0$.  However, since $\check
\calK_{\calXS}$ has basis $\xi_{\tilde F}$ as an $F_{X/S}^*
\scrO_{\calG}$-module, $E_0$ is locally isomorphic to $E$; we now
calculate the connection on $E_0$ in terms of this isomorphism.

\begin{proposition}
  \label{lemma:nablap}
  Suppose we are given a lifting $\tilde F$ of $F_{X/S}$, and let
  $\nabla'$ be the connection on $E$ which makes the following diagram
  commute:
  \[
  \begin{CD}
    \scrHom_{F^* \scrO_{\calG}} (\check \calK_{\calXS}, E) @>
    \cdot(\xi_{\tilde F}) > \sim > E \\
    @VV \nabla V @VV \nabla' V \\
    \scrHom_{F^* \scrO_{\calG}} (\check \calK_{\calXS}, E)
    \otimes \Omega^1_{X/S} @> \cdot(\xi_{\tilde F}) \otimes \id > \sim > E \otimes
    \Omega^1_{X/S}.
  \end{CD}
  \]
  \begin{enumerate}
    \item We have \[ \nabla' = \nabla + (\id \otimes \zeta_{\tilde F}) \circ \psi, \]
      where the second term is the composition $E \overset{\psi}{\to}
      E \otimes_{\scrO_X} F_{X/S}^* \Omega^1_{X'/S} \overset{\id \otimes \zeta_{\tilde
          F}}{\longrightarrow} E \otimes_{\scrO_X} \Omega^1_{X/S}$.
    \item Also suppose we have a logarithmic system of coordinates
      $m_1, \ldots, m_r \in \MXgp$ for $X$, and let $D_1, \ldots, D_r$ be the
      corresponding basis of $\scrT_{X/S}$.  Then $\rho'_{D_i} =
      \rho_{D_i}^p$ for $i = 1, \ldots, r$.
  \end{enumerate}

  \begin{proof}
    To prove (i), we first calculate the connection on $\check
    \calK_{\calXS}$: we see that $(\nabla \xi_{\tilde F})(1) = 0$; $(\nabla
    \xi_{\tilde F})(\sigma_{\tilde F}(\pi^* \dlog(m_i))) = -1 \otimes \zeta_{\tilde
      F}(\pi^* \dlog(m_i))$; and $\nabla \xi_{\tilde F}$ is zero on $S^{\geq 2}
    \sigma_{\tilde F}(F_{X/S}^* \Omega^1_{X'/S})$.  Similarly, $(\psi \xi_{\tilde
      F})(1) = 0$; $(\psi \xi_{\tilde F})(\sigma_{\tilde F}(\pi^* \dlog(m_j))) = 1
    \otimes \dlog(m_j)$; and $\psi \xi_{\tilde F}$ is zero on $S^{\geq 2} \sigma_{\tilde
      F}(F_{X/S}^* \Omega^1_{X'/S})$.  Therefore, $\nabla \xi_{\tilde F} = -(\id \otimes
    \zeta_{\tilde F}) (\psi \xi_{\tilde F})$.  Now consider an element $\phi \in
    \scrHom_{F^* \scrO_{\calG}} (\check \calK_{\calXS}, E)$, and let
    $e := \phi(\xi_{\tilde F})$.  Then
    \begin{align*}
      \nabla' e & = (\nabla \phi)(\xi_{\tilde F}) = \nabla(\phi(\xi_{\tilde F})) -
      (\phi \otimes \id)(\nabla \xi_{\tilde F}) \\
      & = \nabla e - (\phi \otimes \id) \left( -\sum_{i=1}^r \psi_{D_i '} \xi_{\tilde F} \otimes
        \zeta_{\tilde F}(\pi^* \dlog(m_i))
      \right) \\
      & = \nabla e + \sum_{i=1}^r \psi_{D_i '} e \otimes \zeta_{\tilde F}(\pi^* \dlog(m_i))
      \\
      & = \nabla e + (\id \otimes \zeta_{\tilde F}) (\psi e),
    \end{align*}
    as required.
    
    To prove (ii), note that $\zeta_{\tilde F}(\pi^* \dlog(m_j)) =
    \dlog(m_j) + d b_j$ for some $b_j \in \scrO_X$.  Thus, $\langle
    \zeta_{\tilde F}(\pi^* \dlog(m_j)), D_i \rangle = \delta_{ij} + D_i b_j$, so
    from (i) it follows that $\rho_{D_i}' = \rho_{D_i} + \bar \psi_{D_i'}$,
    where $\bar \psi$ is the residue of the $p$-curvature.  However,
    $\psi_{D_i'} = \nabla_{D_i}^p - \nabla_{D_i}$, so $\bar \psi_{D_i'} =
    \rho_{D_i}^p - \rho_{D_i}$, completing the proof.
  \end{proof}
\end{proposition}

\begin{remark}
  If $\nabla$ is locally nilpotent, the following local characterization of
  $E_0 \simeq (\calK_{\calXS} \otimes_{\scrO_X} E)^\theta$ is useful for calculations:
  suppose we have a lifting $\tilde F$ of $F_{X/S}$, which we use to
  identify $\calK_{\calXS}$ with $S^\cdot F_{X/S}^* \Omega^1_{X'/S}$, and a
  logarithmic system of coordinates $m_1, \ldots, m_r \in \MXgp$.  Then
  $(\calK_{\calXS} \otimes_{\scrO_X} E)^\theta$ is the set of elements of
  $\calK_{\calXS} \otimes_{\scrO_X} E$ of the form
  \[ \sum_{N \in \nats^r} (-1)^{|N|} (\pi^* \dlog m)^N \otimes \theta_{(D')^{[N]}} (e)
  \]
  for some $e \in E$, where $D_1', \ldots, D_r'$ is the basis of
  $\scrT_{X'/S}$ dual to the basis $\pi^* \dlog m_1, \ldots,
  \discretionary{}{}{} \pi^* \dlog m_r$
  of $\Omega^1_{X'/S}$.
\end{remark}

\begin{corollary}
  \label{cor:micpd0}
  Let $MIC_{PD}^0(X/S)$ be the full subcategory of $MIC_{PD}(X/S)$
  consisting of objects $(E, \nabla, \theta)$ such that the residue $\rho$ of
  $\nabla$ satisfies $\rho_D^p = 0$ for every $D \in \scrT_{X/S}$.  Define
  functors $C_{\calXS} : MIC_{PD}^0(X/S) \to HIG_{PD}(X'/S)$ and
  $C_{\calXS}^{-1} : HIG_{PD}(X'/S) \to MIC_{PD}^0(X/S)$ by
  \begin{align*}
    C_{\calXS}(E) & := \iota_* \scrHom_{\DXS^\gamma} (\check \calK_{\calXS}, E), \\
    C_{\calXS}^{-1}(E') & := \check \calK_{\calXS} \otimes_{\scrO_{\calG}}
    \iota_* E'
  \end{align*}
  \begin{enumerate}
  \item The functor $C_{\calXS}$ has left adjoint $C_{\calXS}^{-1}$.
  \item The unit $\eta : \id \to C_{\calXS} \circ C_{\calXS}^{-1}$ of this
    adjunction is an isomorphism, and the counit $\epsilon : C_{\calXS}^{-1}
    \circ C_{\calXS} \to \id$ is an epimorphism on every object of
    $MIC_{PD}^0(X/S)$.  (Note that this implies that $C_{\calXS}$ is
    faithful, and $C_{\calXS}^{-1}$ is fully faithful.)
  \item If $\Tor_1(E, \calR_{X/S}) = 0$, then $\epsilon_E$ is an
    isomorphism.  (Hence $E$ is in the essential image of
    $C_{\calXS}^{-1}$.)
  \end{enumerate}

  \begin{proof}
    Note that $\check \calK_{\calXS}$ is itself an object of
    $MIC_{PD}^0(X/S)$, so $C_{\calXS}^{-1}$ does indeed have image
    contained in $MIC_{PD}^0(X/S)$.  Now the adjunction in (i) is
    simply the standard adjunction between $\check\calK_{\calXS} \otimes \cdot$
    and $\scrHom_{\DXS^\gamma}(\check \calK_{\calXS}, \cdot)$ discussed in
    (\ref{rmk:azmadj}).

    For $E \in MIC_{PD}^0(X/S)$, letting $E_0 := \scrHom_{F^*
      \scrO_{\calG}} (\check \calK_{\calXS}, E)$, we see that the
    residue of $\nabla_{tot}$ on $E_0$ is zero.  Thus, we have a natural
    surjection $C_{\calXS}(E) \otimes_{\scrO_{X'}} \calB_{X/S} \to
    C_{\calXS}(E \otimes_{\scrO_X} \AXgp)$, which is an isomorphism if
    $\Tor_1(E_0, \calR_{X/S}) = 0$.  However, since $E_0$ is locally
    isomorphic to $E$, this is equivalent to the condition that
    $\Tor_1(E, \calR_{X/S}) = 0$.  Now applying $C_{\calXS}^{-1}$ to
    this map gives a surjection $C_{\calXS}^{-1}(C_{\calXS}(E))
    \otimes_{\scrO_X} \AXgp \to E \otimes_{\scrO_X} \AXgp$; however, it is
    straightforward to check that taking the degree zero
    part gives the counit $\epsilon$.
    
    Now for $E' \in HIG_{PD}(X' / S)$, let $E := C_{\calXS}^{-1}(E') =
    \check\calK_{\calXS} \otimes_{\scrO_{\calG}} E'$.  Then
    \[ E_0 := \scrHom_{F^* \scrO_{\calG}} (\check \calK_{\calXS}, E) \simeq
    \scrHom_{F^* \scrO_{\calG}} (\check \calK_{\calXS}, \check
    \calK_{\calXS} \otimes_{F^* \scrO_{\calG}} (\scrO_X \otimes_{\scrO_{X'}} E'))
    \simeq \scrO_X \otimes_{\scrO_{X'}} E' \] since $\check \calK_{\calXS}$ is
    invertible as an $F_{X/S}^* \scrO_{\calG}$-module.  Thus, $E_0^\nabla \simeq
    E'$, so $E_0 \simeq E_0^\nabla \otimes_{\scrO_{X'}} \scrO_X$, and we get an
    isomorphism $C_{\calXS}(E) \otimes_{\scrO_{X'}} \calB_{X/S} \simeq
    C_{\calXS}(E \otimes_{\scrO_X} \AXgp)$.  On the other hand,
    $C_{\calXS}^{-1}(E' \otimes_{\scrO_{X'}} \calB_{X/S}) \simeq E \otimes_{\scrO_X}
    \AXgp$, so $E' \otimes_{\scrO_{X'}} \calB_{X/S} \simeq C_{\calXS}(E
    \otimes_{\scrO_X} \AXgp)$.  Again, we see that taking the degree zero
    part of the isomorphism $E' \otimes_{\scrO_{X'}} \calB_{X/S} \simeq
    C_{\calXS}( C_{\calXS}^{-1} (E')) \otimes_{\scrO_{X'}} \calB_{X/S}$
    gives the unit $\eta$.
  \end{proof}
\end{corollary}

Similarly to before, we get isomorphisms $C_{\calXS}(E) \simeq (\check
\calK_{\calXS} \otimes_{\scrO_X} E)^{\nabla, \gamma}$ for $E \in MIC_{PD \cdot}^0 (X/S)$,
and $C_{\calXS}^{-1}(E') \simeq (\check \calK_{\calXS} \otimes_{\scrO_{X'}}
E')^\theta$ for $E' \in HIG_{PD \cdot}(X' / S)$.

\subsection{The Local Cartier Transform}

From the above calculations, we see that for $E \in MIC_{PD}^0(X' / S)$,
given a lifting $\tilde F$ of $F_{X/S}$ and a logarithmic system of
coordinates $m_1, \ldots, m_r \in \MXgp$, we have isomorphisms
\[ C_{\calXS}(E) \simeq (\iota_* E)^{\nabla '} = \{ x \in E : \nabla e = -(\id \otimes \zeta_{\tilde
  F}) (\psi e) \}, \] where $\nabla '$ is as in (\ref{lemma:nablap}), with the
PD Higgs field inherited from $-\psi$ on $E$.  Similarly, for $E' \in
HIG_{PD}(X'/S)$, the isomorphism $\check \calK_{\calXS} \simeq F_{X/S}^*
\scrO_{\calG}$ induced by $\tilde F$ gives an isomorphism
$C_{\calXS}^{-1}(E') \simeq F_{X/S}^* E'$, with the connection given by
\[ \nabla(e' \otimes f) = e' \otimes df + (\id \otimes \zeta_{\tilde F})(\psi' e') \otimes f. \]
Similar formulas hold for the equivalence between $MIC_{PD}^{\calA}(X
/ S)$ and $HIG_{PD}^{\calB}(X' / S)$.
This last formula is our motivation for the following definition.

\begin{definition}
  Let $\zeta : \Omega^1_{X'/S} \to F_{X/S *} \Omega^1_{X/S}$ be a splitting of the
  Cartier operator.  We then define a functor $\Psi : MIC^{\calA}(X / S) \to
  \FHIG^{\calA}(X / S)$ which maps $(E, \nabla)$ to the $p$-curvature of
  $\nabla$ on $E$.  We also define a functor $\Psi_\zeta^{-1} :
  HIG^{\calB}(X' / S) \to MIC^{\calA}(X / S)$ by
  $\Psi_\zeta^{-1} (E', \psi') := E' \otimes_{\calB_{X/S}} \AXgp$, with connection
  \[ \nabla(e' \otimes f) = e' \otimes df + (\id \otimes \zeta)(\psi' e') \otimes f. \]

  Similarly, we define functors $\Psi : MIC(X / S) \to \FHIG(X / S)$
  by $\Psi(E, \nabla) := (E, \psi(\nabla))$ and $\Psi_\zeta^{-1} : HIG(X' / S) \to MIC(X / S)$
  by $\Psi_\zeta^{-1}(E', \theta') := (F_{X/S}^* E', \nabla)$, where $\nabla$ is given by
  the same formula as above.
\end{definition}

We now calculate the $p$-curvature of this connection.  In order to
express the answer, we use the following notation: $\zeta$ gives a map
$F_{X/S}^* \Omega^1_{X'/S} \to \Omega^1_{X/S}$, so the adjoint $\zeta^*$ gives a map
$\scrT_{X/S} \to F_{X/S}^* \scrT_{X'/S}$.

\begin{lemma}
  Let $(E', \theta') \in HIG^{\calB}(X' / S)$ or $HIG(X' / S)$.
  \begin{enumerate}
  \item The connection $\nabla$ on $\Psi_\zeta^{-1}(E', \theta')$ is integrable.
  \item If $\psi$ is the $p$-curvature of $\nabla$, then $\psi_{\pi^* D} =
    (F_{X/S}^* \theta')_{(\zeta^* D)^p - \pi^* D}$ for $D \in \scrT_{X/S}$.
  \end{enumerate}

  \begin{proof}
    The integrability of $\nabla$ follows from the fact that the image of
    $\zeta$ consists of closed forms.  To prove the formula for $\psi$, since
    $\Psi_\zeta^{-1}(E') \simeq E' \otimes_{S^\cdot \scrT_{X'/S}} \Psi_\zeta^{-1}(S^\cdot
    \scrT_{X'/S})$, it suffices to verify the formula for $E' = S^\cdot
    \scrT_{X'/S}$.  The calculation uses the following two identities.

    \begin{lemma}
      \cite[1.4.1]{lorenzon}
      Let $\omega \in Z \Omega^1_{X/S}$ and $D \in \scrT_{X/S}$.  Then
      \[ F_{X/S}^* \langle C_{X/S} \omega, \pi^* D \rangle = \langle \omega, D^{(p)} \rangle - D^{p-1} \langle
      \omega, D \rangle. \]
    \end{lemma}

    \begin{lemma}
      Let $E'$ be an $\scrO_{X'}$-module, and let $\alpha, \beta \in
      \scrEnd_{\scrO_{X'}}(E')$; now define $\beta_n$ recursively by $\beta_0
      = \beta$ and $\beta_{n+1} = [\alpha, \beta_n]$.  Suppose that $\beta_0, \beta_1, \ldots,
      \beta_{p-1}$ commute pairwise.  Then
      \[ (\alpha + \beta)^p = \alpha^p + \beta_{p-1} + \beta^p. \]

      \begin{proof}
        Let $P_n'$ be the set of pairs $(S_0, \pi)$, where $S_0$ is
        a subset of $\{ 1, 2, \ldots, n \}$ and $\pi$ is a partition of $\{
        1, 2, \ldots, n \} - S_0$.  We then prove by induction on $n$ that
        \[ (\alpha + \beta)^n = \sum_{(S_0, \pi) \in P_n'} \left( \prod_{S \in \pi} \beta_{|S| -
            1} \right) \alpha^{|S_0|}. \] For $n = 0$, the statement is
        trivial.  Now for the inductive step, we use the identity $\alpha
        (\phi_1 \phi_2 \cdots \phi_n) = [\alpha, \phi_1] \phi_2 \cdots \phi_n + \phi_1 [\alpha, \phi_2] \cdots \phi_n + \cdots
        + \phi_1 \phi_2 \cdots [\alpha, \phi_n] + \phi_1 \phi_2 \cdots \phi_n \alpha$.  Applying this to a
        term in the sum above, we see that changing $\beta_{|S| - 1}$ to
        $[\alpha, \beta_{|S| - 1}] = \beta_{|S|}$ corresponds to adding $n + 1$ to
        $S$; multiplying on the right by $\alpha$ corresponds to adding $n
        + 1$ to $S_0$; and multiplying on the left by $\beta$ corresponds
        to adding $\{ n + 1 \}$ to $\pi$.  Thus
        \[ (\alpha + \beta) \sum_{(S_0, \pi) \in P_n'} \left( \prod_{S \in \pi} \beta_{|S| - 1}
        \right) \alpha^{|S_0|} = \sum_{(S_0, \pi) \in P_{n+1}'} \left( \prod_{S \in \pi}
          \beta_{|S| - 1} \right) \alpha^{|S_0|}, \]
        completing the induction.
        
        Now considering the action of the cycle $(1, 2, \ldots, p)$ on
        $P_p'$, we see that in the expression for $(\alpha + \beta)^p$, we can
        group the terms in the sum into groups of $p$ identical terms,
        except for those corresponding to the fixed points $(\{ 1, \ldots,
        p \}, \emptyset)$, $(\emptyset, \{ \{ 1, \ldots, p \} \})$, and $(\emptyset, \{ \{ 1 \} ,
        \{ 2 \}, \ldots, \{ p \} \})$.
      \end{proof}
    \end{lemma}

    To finish the calculation of the $p$-curvature of $\Psi_\zeta^{-1}(S^\cdot
    \scrT_{X'/S})$, we note that $\nabla_D = D \otimes \id + \zeta^* D$.  Thus,
    in the previous lemma we set $\alpha := D \otimes \id$, and $\beta := \zeta^* D$.
    Locally, we see that $\beta = \sum_{i=1}^r \langle \zeta(\pi^* \dlog m_i), D \rangle \otimes
    D_i'$, so $\beta_n = \sum_{i=1}^r D^n . \langle \zeta(\pi^* \dlog m_i), D \rangle \otimes D_i'$,
    and in particular $\beta_0, \ldots, \beta_{p-1}$ commute pairwise.  Also,
    \begin{align*}
      \beta_{p-1} & = \sum_{i=1}^r D^{p-1} . \langle \zeta(\pi^* \dlog m_i), D \rangle \otimes D_i '
      \\
      & = \sum_{i=1}^r [ \langle \zeta(\pi^* \dlog m_i), D^{(p)} \rangle - F_{X/S}^* \langle \pi^*
      \dlog m_i, \pi^* D \rangle ] \otimes D_i' = \zeta^* (D^{(p)}) - \pi^* D.
    \end{align*}
    Therefore, $\nabla_D^p = D^{(p)} \otimes \id + \zeta^* (D^{(p)}) - \pi^* D + (\zeta^*
    D)^p$, and $\nabla_{D^{(p)}} = D^{(p)} \otimes \id + \zeta^* (D^{(p)})$.  Hence
    $\psi_{\pi^* D} = (\zeta^* D)^p - \pi^* D$, as required.
  \end{proof}
\end{lemma}

We may express this formula geometrically as follows: let $\bfT_{X/S}
= \bfSpec S^\cdot \scrT_{X/S}$ be the cotangent bundle on $X'$.  Then
$\pi_{X/S}^* (\zeta^*)$ gives a map $\scrT_{X'/S} \to F_{X'}^* \scrT_{X'/S}$,
which corresponds to a morphism $\phi' : \bfT_{X'/S}^{(X')} \to
\bfT_{X'/S}$.  Composing with $F_{\bfT / X'} : \bfT_{X'/S} \to
\bfT_{X'/S}^{(X')}$ gives a map $h_\zeta : \bfT_{X'/S} \to \bfT_{X'/S}$, and
we see $h_\zeta^* (\pi^* D) = (\zeta^* D)^p$ for $D \in \scrT_{X/S}$.  Thus,
letting $\alpha_\zeta := \id - h_\zeta$, we get
a commutative diagram
\[
\begin{CD}
  HIG(X'/S) @> \alpha_{\zeta *} \iota_* >> HIG(X'/S) \\
  @VV \Psi_\zeta^{-1} V @VV F_{X/S}^* V \\
  MIC(X/S) @> \Psi >> \FHIG(X/S),
\end{CD}
\]
and similarly replacing $HIG(X'/S)$, $MIC(X/S)$, and $\FHIG(X/S)$ by
$HIG^{\calB}(X'/S)$, $MIC^{\calA}(X/S)$, and $\FHIG^{\calA}(X/S)$,
respectively.  (We choose this sign for $\alpha_\zeta$ so that we may think of
it as a perturbation of the identity map.)

Now the restriction $\hat \alpha_\zeta$ of $\alpha_\zeta$ to the completion $\hat\bfT_{X'/S}$ of
$\bfT_{X'/S}$ over the zero section is an isomorphism, with inverse
\[ \hat \alpha_\zeta^{-1} = \id + h_\zeta + h_\zeta^2 + \cdots. \]
This allows us to construct a splitting module for $(\tildeDXS)
\hat{\,} := \tildeDXS \otimes_{S^\cdot \scrT_{X'/S}} \hat S^\cdot \scrT_{X'/S}$ over
$\tilde\frakZ \hat{\,} := \tilde\frakZ \otimes_{S^\cdot \scrT_{X'/S}} \hat S^\cdot
\scrT_{X'/S} \simeq \calB_{X/S} \otimes_{\scrO_{X'}} \hat S^\cdot \scrT_{X'/S}$,
namely $\check\calK_\zeta^{\calA} := \Psi_\zeta^{-1} (\hat \alpha_\zeta^{-1})_*
\iota_* \tilde\frakZ \hat{\,}$.  By the above diagram we see that the
$p$-curvature on $\check\calK_\zeta^{\calA} \simeq \tilde\frakZ \hat{\,}
\otimes_{\calB_{X/S}} \AXgp$ is just the natural action of $S^\cdot \scrT_{X/S}$
on $\tilde\frakZ \hat{\,}$, so the action of $\hat S^\cdot \scrT_{X'/S}$
on $\tilde\frakZ \hat{\,}$ and the connection on
$\check\calK_\zeta^{\calA}$ extend to a $(\tildeDXS) \hat{\,}$-module
structure.  Also, since $\AXgp$ is locally free of rank $p^d$ over
$\calB_{X/S}$, $\check\calK_\zeta^{\calA}$ is locally free of rank $p^d$
over $\tilde\frakZ \hat{\,}$, showing that $\check\calK_\zeta^{\calA}$ is
indeed a splitting module.  However, since $\hat \alpha_\zeta^{-1}$ is an
isomorphism, we have a canonical isomorphism $\check\calK_\zeta^{\calA} \simeq
\tilde\frakZ \hat{\,} \otimes_{\calB_{X/S}} \AXgp \simeq \AXgp \otimes_{\scrO_X}
F_{X/S}^* \hat S^\cdot \scrT_{X'/S}$.

\begin{theorem}
  \label{thm:local-rh}
  Let $\widehat{MIC}^{\calA} (X / S)$ be the category of
  $(\tildeDXS) \hat{\,}$-modules, and let $\widehat{HIG}^{\calB} (X' /
  S)$ be the category of $\tilde\frakZ \hat{\,}$-modules.  Thus
  $\widehat{MIC}^{\calA} (X / S)$ is equivalent to the category of
  $\AXgp$-modules with integrable and admissible connections, with an
  extension of the $p$-curvature to a $\hat \bfT_{X'/S}$-Higgs field
  which commutes with $\AXgp$, and $\widehat{HIG}^{\calB} (X' / S)$ is
  equivalent to the $\calB_{X/S}$-modules with $\hat
  \bfT_{X'/S}$-Higgs fields which commute with $\calB_{X/S}$.

  \begin{enumerate}
  \item The functors
    \begin{align*}
      C_\zeta : \widehat{MIC}^{\calA} (X / S) \to \widehat{HIG}^{\calB} (X'
      / S), & ~ E \mapsto \iota_* \scrHom_{(\tildeDXS) \hat{\,}}
      (\check\calK_\zeta^{\calA}, E), \\
      C_\zeta^{-1} : \widehat{HIG}^{\calB} (X' / S) \to \widehat{MIC}^{\calA}
      (X / S), & ~ E' \mapsto \check\calK_\zeta^{\calA} \otimes_{\tilde\frakZ
        \hat{\,}} \iota_* E',
    \end{align*}
    are quasi-inverse equivalences of categories.
  \item Let $E \in \widehat{MIC}^{\calA} (X / S)$, and let $E'$ be the
    corresponding element of $\widehat{HIG}^{\calB} (X' / S)$.  Then
    we have isomorphisms
    \[ E \simeq \Psi_\zeta^{-1} (\hat \alpha_\zeta^{-1})_* E' \simeq \Psi_\zeta^{-1} \hat \alpha_\zeta^* E'
    \]
    and
    \[ \Psi(E) \simeq \iota_* E' \otimes_{\calB_{X/S}} \AXgp. \]
    Here in the second formula, $S^\cdot \scrT_{X'/S}$ acts on the right
    hand side via its action on $\iota_* E'$.
  \end{enumerate}

  \begin{proof}
    The first part follows directly from the fact that $\check\calK_\zeta
    ^{\calA}$ is a splitting module.  For the second part, since
    $\Psi_\zeta^{-1} E' \simeq \Psi_\zeta^{-1} S^\cdot \scrT_{X'/S} \otimes_{S^\cdot \scrT_{X'/S}} E'$,
    a base extension gives $\Psi_\zeta^{-1} E' \simeq \Psi_\zeta^{-1} \iota_* \tilde\frakZ
    \hat{\,} \otimes_{\tilde\frakZ \hat{\,}} E'$.  Thus,
    we have
    \[ E \simeq \Psi_\zeta^{-1} (\hat \alpha_\zeta^{-1})_* \iota_* \tilde \frakZ \hat{\,}
    \otimes_{\tilde \frakZ \hat{\,}} \iota_* E' \simeq \Psi_\zeta^{-1} (\hat \alpha_\zeta)^* \iota^*
    \tilde \frakZ \hat{\,} \otimes_{\tilde \frakZ \hat{\,}} \iota^* E'
    \simeq \Psi_\zeta^{-1} (\hat \alpha_\zeta)^* E' \simeq \Psi_\zeta^{-1} (\hat \alpha_\zeta^{-1})_* E', \]
    which proves the first isomorphism.  The second follows from the
    commutative diagram above.
  \end{proof}
\end{theorem}

\begin{remark}
$\widehat{MIC}^{\calA}(X / S)$ is {\em not} equivalent to
$MIC_\cdot^{\calA} (X / S)$.  (Indeed, $\check \calK_\zeta^{\calA}$ is itself
a counterexample.)  However, $MIC_\cdot^{\calA} (X / S)$ and
$HIG_\cdot^{\calB} (X' / S)$ can be regarded as full subcategories of
$\widehat{MIC}^{\calA} (X / S)$ and $\widehat{HIG}^{\calB} (X' / S)$.
Similarly to before, when $C_\zeta$ and $C_\zeta^{-1}$ are restricted to these
subcategories, we have isomorphisms $C_\zeta(E) \simeq (\calK_\zeta^{\calA}
\otimes_{\AXgp} E)^{\nabla, \gamma}$ and $C_\zeta^{-1}(E') \simeq (\calK_\zeta^{\calA}
\otimes_{\calB_{X/S}} E')^\theta$, where $\calK_\zeta^{\calA} := \Psi_\zeta^{-1}
(\alpha_\zeta^{-1})_* (\AXgp \otimes_{\scrO_{X'}} \Gamma_\cdot F_{X/S}^* \Omega^1_{X'/S}) \simeq
\scrHom_{\AXgp} (\check \calK_\zeta^{\calA}, \AXgp)$.
\end{remark}

Similarly, let $\widehat{MIC}^0(X / S)$ be the category of $\hat
\calD_{X/S} := \DXS \otimes_{S^\cdot \scrT_{X'/S}} \hat S^\cdot
\scrT_{X'/S}$-modules with residue $\rho$ satisfying $\rho_D^p = 0$ for all
$D \in \scrT_{X/S}$, and $\widehat{HIG}(X' / S)$ be the category of
$\hat S^\cdot \scrT_{X'/S}$-modules.  Also, let $\check\calK_\zeta := \Psi_\zeta^{-1}
(\hat \alpha_\zeta^{-1})_* \iota_* \hat S^\cdot \scrT_{X'/S}$.  Then $C_\zeta^{-1} : E' \mapsto
\iota_* E' \otimes_{\hat S^\cdot \scrT_{X'/S}} \check\calK_\zeta$ gives an ``almost
equivalence'' $\widehat{HIG}(X' / S) \to \widehat{MIC}^0(X / S)$ in the
sense of (\ref{cor:micpd0}), with left inverse $C_\zeta : E \mapsto
\scrHom_{\hat \calD_{X/S}} (\check\calK_\zeta, \iota_* E)$.  Yet again, letting
$\calK_\zeta := \Psi_\zeta^{-1} (\hat \alpha_\zeta^{-1})_* \Gamma_\cdot \Omega^1_{X'/S} \simeq
\scrHom_{\scrO_X} (\check \calK_\zeta, \scrO_X)$, we get that $C_\zeta(E) \simeq
(\calK_\zeta \otimes_{\scrO_X} E)^{\nabla, \gamma}$ for $E \in MIC_\cdot^0(X / S)$, and
$C_\zeta^{-1}(E') \simeq (\calK_\zeta \otimes_{\scrO_{X'}} E')^\theta$ for $E' \in HIG_\cdot(X' /
S)$.

\subsection{The Image of the $p$-curvature Functor}

In this section we give an application of the local theory developed
in the last section to characterize the essential image of the functor
$\Psi:MIC(X/S) \to \FHIG(X/S)$.

\begin{theorem}
  Let $X \to S$ be a smooth morphism of fine log schemes of
  characteristic $p$.
  \begin{enumerate}
  \item Let $(E, \nabla)$ be a finitely generated $\AXgp$-module with
    integrable connection on $X$.  Then \'etale locally on the
    underlying scheme of $X$, there exists a $\calB$-Higgs field $(E',
    \theta')$ on $X'$ and an isomorphism $(E, \psi(\nabla)) \simeq
    (E', \theta') \otimes_{\calB_{X/S}} \AXgp$, where $\psi(\nabla)$
    is the $p$-curvature of $\nabla$.
  \item Let $(E', \theta')$ be a finitely generated $\calB$-Higgs
    field on $X'$.  Then \'etale locally on the underlying scheme on
    $X$, there exists an $\AXgp$-module with integrable connection
    $(E, \nabla)$ on $X$ and an isomorphism $(E, \psi(\nabla)) \simeq
    (E', \theta') \otimes_{\calB_{X/S}} \AXgp$.
  \end{enumerate}

  \begin{proof}
    We begin by stating the lemma we will use to construct the
    appropriate \'etale morphisms.

    \begin{lemma}
      Let $Y \to X$ be a morphism of noetherian schemes, $\tilde Y \to Y$
      a surjective \'etale morphism, and $Z \to Y$ a closed immersion
      such that the composition $Z \to X$ is finite.  Then \'etale
      locally on $X$, there exists a section of $\tilde Y \to Y$ over
      $Z$.

      \begin{proof}
        First, we may replace $Y$ by $Z$ and $\tilde Y$ by $\tilde Z
        := \tilde Y \times_Y Z$.  Thus we may assume that $Y \to X$ is
        finite.  Suppose that $X$ is the spectrum of a strictly
        henselian local ring $A$.  Then by Hensel's lemma, $Z$ splits
        up as a finite direct sum of schemes $Z_i$ where each $Z_i$ is
        the spectrum of a strictly henselian local ring.  Then each
        $\tilde Z_i := \tilde Z \times_Z Z_i$ is \'etale and surjective
        over $Z_i$, and hence admits a section.  Thus $\tilde Z \to Z$
        also admits a section.

        For the general case, let $A$ be a ring endowed with a map
        $\Spec A \to X$, and let $F(A)$ denote the set of sections of
        $\tilde Z \to Z$ over $Z \times_X \Spec A$.  Since $\tilde Z \to X$ is
        of finite presentation, the functor $F$ commutes with direct
        limits.  (We may assume that $X$ is affine, say $\Spec R$.
        Then $Z$ is also affine, say $\Spec B$, and $Z \times_R A = \Spec(B
        \otimes_R A)$, and the functor $A \mapsto B \otimes_R A$ commutes with direct
        limits.)  We have seen that $F(A)$ is nonempty if $A$ is the
        strict henselization of the local ring of $X$ at any point
        $x$.  But this $A$ is a direct limit of the rings
        corresponding to \'etale neighborhoods of $x$.
      \end{proof}
    \end{lemma}

    Now to prove (ii), let $Z \subseteq \bfT_{X'/S}$ be the support of
    $E'$, which is finite over $X'$ since $E'$ is finitely generated
    (as $\calB_{X/S}$ is finitely generated over $\calA_{X'}^{\gp}$).
    Locally, we may choose a splitting $\zeta$ of $C_{X/S}$.  Then
    $\alpha_\zeta : \bfT_{X'/S} \to \bfT_{X'/S}$ is surjective
    \'etale, so \'etale locally on $X$ there exists a lifting $g : Z
    \to \bfT_{X'/S}$ such that $\alpha_\zeta \circ g : Z \to
    \bfT_{X'/S}$ is the closed immersion.  Let $(E, \nabla) :=
    \Psi_{\zeta}^{-1} (\iota_* g_* E')$.  Then $(E, \nabla)$ has
    $p$-curvature $F_{X/S}^* (\alpha_{\zeta *} g_* E') = F_{X/S}^*
    E'$.

    Similarly, for (i), let $Z \subseteq \bfT_{X'/S}$ be the support
    of the $p$-curvature of $\nabla$, which is again finite over
    $X'$.  Then as before, we may \'etale locally choose a splitting
    $\zeta$ such that $\alpha_\zeta : \bfT_{X'/S} \to \bfT_{X'/S}$ has
    a section $g$ over $Z$.  Then $g_* (E, \psi)$ is a
    $\calB$-Higgs field compatible with $(E, \psi)$ via
    $\alpha_\zeta$.  Thus, since $\alpha_\zeta$ splits the Azumaya
    algebra $\tildeDXS$, $(E, \psi) \simeq \Psi_\zeta^{-1}(E',
    \theta')$ for some $\calB$-Higgs field $(E', \theta')$.
    Therefore, $(E, \psi) \simeq (\alpha_{\zeta*} \iota_* E')
    \otimes_{\calB_{X/S}} \AXgp$.
  \end{proof}
\end{theorem}

\begin{example}
  For the case of $E = \scrO_X$, we can in fact say more: \'etale
  locally, $(\scrO_X, \omega)$ is the $p$-curvature of some integrable
  connection on $\scrO_X$ if and only if $\omega = F_{X/S}^* \omega'$
  for some $\omega' \in \Omega^1_{X'/S}$.  Indeed, for any closed
  1-form $\omega$ on $X$, let $\nabla_\omega$ be the connection on
  $\scrO_X$ such that $\nabla_\omega(1) = \omega$.  Then the
  $p$-curvature of $\nabla_\omega$ is multiplication by $F_{X/S}^*
  (\pi^* - C_{X/S})(\omega)$, where $C_{X/S} : Z^1_{X/S} \to
  \Omega^1_{X'/S}$ is the map induced by the Cartier isomorphism.

  Conversely, for any $f \in \scrO_{X'}$, $m \in \MXgp$, let $\omega'
  = f \dlog(\pi^* m)$, and let $\theta' = \cdot \omega'$ be the Higgs
  field on $\scrO_{X'}$ induced by $\omega'$.  Now let $g \in \scrO_X$
  be a root of $g^p - g - F_{X/S}^* f = 0$, which clearly exists
  \'etale locally on $X$.  Then since $g = F_{X/S}^* (\pi^* g - f)$,
  \[ C_{X/S} (g \dlog m) = (\pi^* g - f) \dlog(\pi^* m) = \pi^*(g
  \dlog m) - f \dlog(\pi^* m). \] Therefore, letting $\omega = g \dlog
  m$, $\nabla_\omega$ has $p$-curvature $F_{X/S}^* \omega'$.  Since
  $\pi^* - C_{X/S}$ is additive, this shows that $\pi^* - C_{X/S} :
  Z^1_{X/S} \to \Omega^1_{X'/S}$ is surjective with respect to the
  \'etale topology on the underlying scheme of $X$.
\end{example}

\section{De Rham and Higgs Cohomology}

Let $\calXS := (X / S, \tilde X' / \tilde S)$ be as in the previous
chapters.  The main result of this section is as follows:

\begin{theorem}
  \label{thm:cohom-main}
  Let $\ell < p$ be a natural number, let $(E, \nabla)$ be an object of
  $MIC_\ell^{\calA}(X / S)$, i.e.{} an object of $MIC^{\calA}(X / S)$
  whose $p$-curvature is nilpotent of level $\leq \ell$, and let $(E', \theta')$
  be its Cartier transform.  Let $m := p - \ell - 1$.  Then there is a
  canonical isomorphism in the derived category of complexes of
  $\calB_{X/S}$-modules
  \[ \tau_{\leq m} (E \otimes \Omega^\cdot_{X/S}, \nabla) \simeq \tau_{\leq m} (E' \otimes \Omega^\cdot_{X'/S}, \theta'). \]
\end{theorem}

Note that if $\ell + \dim(X / S) < p$, this implies that $F_{X/S *} (E \otimes
\Omega^\cdot_{X/S}, \nabla) \simeq (E' \otimes \Omega^\cdot_{X'/S}, \theta')$ in the derived category.

In fact, we shall make this isomorphism explicit as follows: let
\[ \calK_{\calXS}^{\calA, ij} (E) := E \otimes \calK_{\calXS}^{\calA} \otimes \Omega^i_{X'/S}
\otimes \Omega^j_{X/S}. \]
Then for fixed $i$, the de Rham complex of $E \otimes \calK_{\calXS}^{\calA}
\otimes \Omega^i_{X'/S}$ with the total connection gives a complex
\[ \calK_{\calXS}^{\calA, i0}(E) \overset{d}\longrightarrow \calK_{\calXS}^{\calA, i1}(E)
\overset{d}\longrightarrow \cdots \] Similarly, for fixed $j$, the Higgs complex of
$(\calK_{\calXS}^{\calA}, \psi_{\calK^{\calA}})$ tensored with the
identity on $E \otimes \Omega^j_{X/S}$ gives a complex
\[ \calK_{\calXS}^{\calA, 0j}(E) \overset{d'}\longrightarrow \calK_{\calXS}^{\calA, 1j}(E)
\overset{d'}\longrightarrow \cdots \] The differentials $d$ and $d'$ commute and
therefore form a double complex $\calK_{\calXS}^{\calA \cdot \cdot}(E)$.  In fact,
we see that $d$ preserves the subcomplex $E \otimes N_n
\calK_{\calXS}^{\calA} \otimes \Omega^i_{X'/S} \otimes \Omega^j_{X/S}$, and that $d'$ maps
$E \otimes N_n \calK_{\calXS}^{\calA} \otimes \Omega^i_{X'/S} \otimes \Omega^j_{X/S}$ to $E \otimes
N_{n-1} \calK_{\calXS}^{\calA} \otimes \Omega^{i+1}_{X'/S} \otimes \Omega^j_{X/S}$.  Thus,
we get a subcomplex $N_n \calK_{\calXS}^{\calA \cdot \cdot}(E)$ defined by
\[ N_n \calK_{\calXS}^{\calA, ij}(E) := E \otimes N_{n-i} \calK_{\calXS}^{\calA} \otimes
\Omega^i_{X'/S} \otimes \Omega^j_{X/S}. \]
This double complex looks like:
\[
\begin{CD}
  @AAA @AAA @AAA \\
  E \otimes \calK_n^{\calA} \otimes \Omega^2_{X/S} @> d' >> E \otimes \calK_{n-1}^{\calA} \otimes
  \Omega^1_{X'/S} \otimes \Omega^2_{X/S} @> d' >> E \otimes \calK_{n-2}^{\calA} \otimes
  \Omega^2_{X'/S} \otimes \Omega^2_{X/S} @>>> \\
  @AA d A @AA d A @AA d A \\
  E \otimes \calK_n^{\calA} \otimes \Omega^1_{X/S} @> d' >> E \otimes \calK_{n-1}^{\calA} \otimes
  \Omega^1_{X'/S} \otimes \Omega^1_{X/S} @> d' >> E \otimes \calK_{n-2}^{\calA} \otimes
  \Omega^2_{X'/S} \otimes \Omega^1_{X/S} @>>> \\
  @AA d A @AA d A @AA d A \\
  E \otimes \calK_n^{\calA} @> d' >> E \otimes \calK_{n-1}^{\calA} \otimes \Omega^1_{X'/S} @> d' >>
  E \otimes \calK_{n-2}^{\calA} \otimes \Omega^2_{X'/S} @>>>
\end{CD}
\]
We denote by $N_n \calK_{\calXS}^{\calA \cdot}(E)$ the total complex
associated to this double complex.

For any natural number $n$, we have a morphism of complexes
\[ E \otimes \Omega^\cdot_{X/S} \to E \otimes N_n \calK_{\calXS}^{\calA} \otimes \Omega^\cdot_{X/S} \]
via the inclusion $\AXgp \to N_n \calK_{\calXS}^{\calA}$, which when
composed with $d'$ gives zero.  This induces a map $b : E \otimes \Omega^\cdot_{X/S}
\to N_n \calK_{\calXS}^{\calA \cdot} (E)$.  Similarly, if $E$ has level $\ell$,
then $E' \simeq (\calK_{\calXS}^{\calA} \otimes E)^{\nabla, \gamma} \simeq (N_n
\calK_{\calXS}^{\calA} \otimes E)^{\nabla, \gamma}$ for $n \geq \ell$.  Thus, the inclusion
$E' \to N_{n - i} \calK_{\calXS}^{\calA} \otimes E$ gives a morphism of
complexes
\[ \tau_{\leq (n - \ell)} (E' \otimes \Omega^\cdot_{X'/S}) \to \tau_{\leq (n - \ell)} (E \otimes N_{n - \cdot}
\calK_{\calXS}^{\calA} \otimes F_{X/S}^* \Omega^\cdot_{X'/S}), \] which when composed
with $d$ gives zero.  This induces a map $a : \tau_{\leq (n - \ell)} (E' \otimes
\Omega^\cdot_{X'/S}) \to \tau_{\leq (n - \ell)} N_n \calK_{\calXS}^{\calA \cdot}(E)$.

If there is a splitting $\zeta : \Omega^1_{X'/S} \to F_{X/S *} \Omega^1_{X/S}$ of the
Cartier operator $C_{X/S}$, then we may make analogous constructions
with $\calK_\zeta^{\calA}$ in place of $\calK_{\calXS}^{\calA}$.

\begin{theorem}
  \label{thm:cohom-detail}
  Let $(E, \nabla)$ be an object of $MIC_\ell^{\calA}(X / S)$.
  \begin{enumerate}
  \item Let $n < p$ be a natural number with $n \geq \ell$, and let $(E',
    \theta') := C_{\calXS}(E, \nabla)$.  Then the natural maps
    \[ \tau_{\leq (n - \ell)} (E \otimes \Omega^\cdot_{X/S}, \nabla) \overset{b}\longrightarrow \tau_{\leq (n - \ell)}
    (N_n \calK_{\calXS}^{\calA \cdot} (E)) \overset{a}\longleftarrow \tau_{\leq (n - \ell)} (E' \otimes
    \Omega_{X'/S}^\cdot, \theta') \]
    are quasi-isomorphisms.
  \item Assume that there exists a splitting $\zeta : \Omega^1_{X'/S} \to F_{X/S
      *} \Omega^1_{X/S}$ of $C_{X/S}$.  Then the analogs
    \[ (E \otimes \Omega^\cdot_{X/S}, \nabla) \overset{g}\longrightarrow \calK_\zeta^{\calA \cdot}(E)
    \overset{f}\longleftarrow (E' \otimes \Omega^\cdot_{X'/S}, \theta') \] of $a$ and $b$ are
    quasi-isomorphisms.  In fact, for any integer $n \geq \ell$, the analogs
    \[ \tau_{\leq (n - \ell)} (E \otimes \Omega^\cdot_{X/S}, \nabla) \overset{g_n}\longrightarrow \tau_{\leq (n - \ell)}
    (N_n \calK_\zeta^{\calA \cdot}(E)) \overset{f_n}\longleftarrow \tau_{\leq (n - \ell)} (E' \otimes
    \Omega^\cdot_{X'/S}, \theta') \]
    of $a$ and $b$ are quasi-isomorphisms.
  \end{enumerate}
\end{theorem}
Note that (\ref{thm:cohom-main}) follows easily from the case $n = p -
1$ in the first part.

\begin{example}
  Let $(E, \nabla) := (\AXgp, d)$, so that $(E', \theta') = (\calB_{X/S}, 0)$.
  Then for $m \in \MXgp$, the image of $\dlog m \in \Omega^1_{X/S}$ under $b$ is $1 \otimes
  \dlog m \in \calK^{\calA} \otimes \Omega^1_{X/S}$; similarly, the image of $\pi^*
  \dlog m \in \Omega^1_{X'/S}$ under $a$ is $1 \otimes \pi^* \dlog m \in \calK^{\calA}
  \otimes \Omega^1_{X'/S}$.  However, if
  we have a lifting $\tilde \pi$ of $\pi_{X/S}$ and a lifting $\tilde m \in
  \scrM_{\tilde X}^{\gp}$ of $m$, then
  \[ (d \oplus d') \beta_{\tilde \pi}(\tilde m) = 1 \otimes \dlog m - 1 \otimes \pi^* \dlog m \in
  (\calK^{\calA} \otimes \Omega^1_{X/S}) \oplus (\calK^{\calA} \otimes \Omega^1_{X'/S}). \] Hence
  $a(\pi^* \dlog m)$ and $b(\dlog m)$ have the same image in $\calH^1
  \calK^{\calA \cdot}_{\calXS}$, so the isomorphism above is compatible
  with the standard Cartier isomorphism.  In general, if we only have
  a local lifting $\tilde F$ of $F_{X/S}$, then the difference is
  locally the image of $g + \sigma_{\tilde F}(1 \otimes \pi^* \dlog m)$, where we
  construct $g \in \scrO_X$ by taking a local lifting $\tilde m' \in
  \scrM_{\tilde X'}^{\gp}$ of $\pi^* m$ and letting $g$ satisfy $1 + [p]
  g = \alpha_{\tilde X}(\tilde F^* \tilde m' - p \tilde m)$.
\end{example}

To prove (\ref{thm:cohom-detail}), we begin with the case $(E, \nabla) =
(\AXgp, d)$.

\begin{lemma}
  Assume that $\zeta$ is a splitting of $C_{X/S}$.
  \label{lemma:h-dr-res}
  \begin{enumerate}
  \item \label{lemma:higgs-res}
    The Higgs complex $(\calK_\zeta^{\calA} \otimes F_{X/S}^* \Omega^\cdot_{X'/S}, \psi)$ is a
    resolution of $\AXgp$.  In fact, for any $n \geq 0$,
    \[ N_n \calK_\zeta^{\calA} \overset{\psi}{\longrightarrow} N_{n-1} \calK_\zeta^{\calA} \otimes
    \Omega^1_{X'/S} \overset{\psi}{\longrightarrow} N_{n-2} \calK_\zeta^{\calA} \otimes \Omega^2_{X'/S}
    \overset{\psi}{\longrightarrow} \cdots \]
    is also a resolution of $\AXgp$.
  \item \label{lemma:derham-res}
    The de Rham complex $(\calK_\zeta^{\calA} \otimes \Omega^\cdot_{X/S}, \nabla)$ is a
    resolution of $\calB_{X/S}$.
  \end{enumerate}

  \begin{proof}
    In (i), the statement for $(\calK_\zeta^{\calA} \otimes \Omega^\cdot_{X'/S}, \psi)$
    follows from the statement for $(N_{n - \cdot} \calK_\zeta^{\calA} \otimes
    \Omega^\cdot_{X'/S}, \psi)$ by taking direct limits.  Since the latter complex
    is just $\AXgp$ for $n = 0$, it suffices to show that the
    successive quotients are exact.  However, these quotients are
    equivalent to
    \[ \Gamma_n (\AXgp \otimes \Omega^1_{X'/S}) \to \Gamma_{n-1} (\AXgp \otimes \Omega^1_{X'/S}) \otimes
    \Omega^1_{X'/S} \to \Gamma_{n-2} (\AXgp \otimes \Omega^1_{X'/S}) \otimes
    \Omega^2_{X'/S} \to \cdots, \]
    with the maps given by tensoring the maps in the divided
    power de Rham complex of $\Omega^1_{X'/S}$ by $\AXgp$.  Thus, the
    result follows by the divided power Poincar\'e lemma.

    To prove (ii), we note that the connection on $\calK_\zeta^{\calA}$
    preserves the filtration $N_\cdot$, and we consider the spectral
    sequence of the filtered complex $(N_\cdot \calK_\zeta^{\calA} \otimes
    \Omega^\cdot_{X/S}, \nabla)$.  In order to be compatible with standard notation,
    let us denote $N^i := N_{-i}$.  We begin with an indexed version
    of \cite[5.1.1]{ogus-higgs}, which in fact implies this version by
    applying it in each degree.

    \begin{lemma}
      Let $(E, \nabla) \in MIC^{\calA}(X / S)$, and let $\psi : E \to E
      \otimes_{\scrO_{X'}} \Omega^1_{X'/S}$ be the $p$-curvature of $\nabla$.  Suppose
      we have a horizontal filtration $N^\cdot$ on $E$ such that $\psi$
      vanishes on $\Gr^i_N E$ for each $E$.  Then in the spectral
      sequence of the filtered complex $(N^\cdot E \otimes \Omega^\cdot_{X/S}, \nabla)$, there
      is a commutative diagram
      \[
      \begin{CD}
        E_1^{i,j} @> d_1^{i,j} >> E_1^{i+1, j} \\
        @VV \sim V @VV \sim V \\
        \Gr^i_N(E)^\nabla \otimes \Omega^{i+j}_{X'/S} @> -\bar \psi >> \Gr^{i+1}_N(E)^\nabla \otimes
        \Omega^{i+j+1}_{X'/S},
      \end{CD}
      \]
      where $\bar \psi$ is the map induced by $\psi$.
      \begin{proof}
        For any fixed $i,j$, by replacing $E$ with $E / N^k E$ for $k$
        sufficiently large, we may assume $E$ is locally nilpotent.
        Then given a local lifting $\zeta$ of $C_{X/S}$, by
        (\ref{thm:local-rh}), $(E, \nabla)$ is in the essential image of
        $\Psi_\zeta^{-1}$; say $(E, \nabla) \simeq \Psi_\zeta^{-1} (E', \theta')$.  Also, the
        filtration on $E$ descends to a filtration on $E'$.  Now we
        have $\psi = \alpha_{\zeta *} \theta' \otimes \id : E' \otimes_{\calB_{X/S}} \AXgp \to E'
        \otimes_{\calB_{X/S}} \AXgp \otimes_{\scrO_{X'}} \Omega^1_{X'/S}$.  However,
        since $E'$ is locally nilpotent and $\hat \alpha_\zeta$ is an
        isomorphism, we get that $\theta'$ vanishes on $N^i E'$.

        Now by definition, $\nabla = d_{\calA} + (\id \otimes \zeta) \circ \theta'$; however,
        $(\id \otimes \zeta) \circ \theta'$ maps $N^i (E \otimes \Omega^{i+j}_{X/S})$ to $N^{i+1} (E
        \otimes \Omega^{i+j+1}_{X/S})$.  Therefore, $E_1^{i,j} \simeq \calH^{i+j} (\Gr^i_N
        E' \otimes \AXgp \otimes \Omega^\cdot_{X/S}, d_{\calA})$, and $d_1^{i,j}$ is the
        map induced by $(\id \otimes \zeta) \circ \theta'$.  However, since each term of
        the complex $(\AXgp \otimes \Omega^\cdot_{X/S}, d)$ is locally free over
        $\AXgp$ and therefore over $\calB_{X/S}$, and the differential
        $d$ is $\calB_{X/S}$-linear, from (\ref{prop:cart-isom}) we
        conclude
        \[ E_1^{i,j} \simeq \Gr^i_N(E') \otimes \Omega^{i+j}_{X'/S} \simeq E_1^{0,j} \otimes
        \Omega^{i+j}_{X'/S} \simeq \Gr^i_N(E)^\nabla \otimes \Omega^{i+j}_{X'/S}. \] This gives
        the vertical isomorphisms.  Also, since $C_{X/S} \circ \zeta = \id$,
        it is easy to see that the map $\Gr^i_N(E') \otimes \Omega^{i+j}_{X'/S} \to
        \Gr^{i+1}_N(E') \otimes \Omega^{i+j+1}_{X'/S}$ induced by $(\id \otimes \zeta) \circ
        \theta'$ is just $\bar \theta'$.  However, $\psi = \alpha_{\zeta *} \iota_* \theta' \otimes \id$,
        while $h_\zeta$ maps $N^i E' \otimes \AXgp \otimes \Omega^{i+j}_{X/S}$ to $N^{i +
          p} E' \otimes \AXgp \otimes \Omega^{i+j+1}_{X/S}$, so $\bar \theta' = -\bar \psi$.
      \end{proof}
    \end{lemma}

    In our situation with $E = \calK_\zeta^{\calA} \simeq \AXgp \otimes \Gamma_\cdot
    \Omega^1_{X'/S}$, we have $\Gr^{-i}_N(E) \simeq \AXgp \otimes \Gamma_i \Omega^1_{X'/S}$ with
    connection $d_{\calA}$, so $E_1^{-i,j} \simeq \calB_{X/S} \otimes \Gamma_i
    \Omega^1_{X'/S} \otimes \Omega^{j-i}_{X'/S}$.  Also, the $p$-curvature of
    $\calK_\zeta^{\calA}$ is exactly the map we get by tensoring the
    standard differential $d$ in the divided power de Rham complex of
    $\Gamma_\cdot \Omega^1_{X'/S}$ with $\AXgp$, so that we get a commutative
    diagram
    \[
    \begin{CD}
      E_1^{-i,j} @> d_1^{-i,j} >> E_1^{-i+1, j} \\
      @VV \sim V @VV \sim V \\
      \calB_{X/S} \otimes \Gamma_i \Omega^1_{X'/S} \otimes \Omega^{j-i}_{X'/S} @> -\id \otimes d >>
      \calB_{X/S} \otimes \Gamma_{i-1} \Omega^1_{X'/S} \otimes \Omega^{j-i+1}_{X'/S}.
    \end{CD}
    \]
    Thus, by the divided power Poincar\'e lemma, we get $E_2^{i,j}
    = 0$ unless $i = j = 0$.  Hence the cohomology of
    $(\calK_\zeta^{\calA} \otimes \Omega^\cdot_{X/S}, \nabla)$ vanishes in positive degree,
    and we have an isomorphism
    \[ \calH^0(\calK_\zeta^{\calA} \otimes \Omega^\cdot_{X/S}, \nabla) \simeq E_2^{0,0} \simeq
    \calB_{X/S}. \]
  \end{proof}
\end{lemma}

\begin{lemma}
  Let $(E, \nabla) \in MIC^{\calA}_\cdot(X / S)$, and let $(E', \theta') := C_\zeta(E, \nabla) \in
  HIG^{\calB}_\cdot(X' / S)$ be its Cartier transform.

  \begin{enumerate}
  \item The Higgs complex $(E \otimes \calK_\zeta^{\calA} \otimes \Omega^\cdot_{X'/S},
    \psi_{\calK})$ is a resolution of $E$.  In fact, for any natural
    number $n$, the complex
    \[ E \otimes N_n \calK_\zeta^{\calA} \overset{\psi_{\calK}} \longrightarrow E \otimes
    N_{n-1} \calK_\zeta^{\calA} \otimes \Omega^1_{X'/S} \overset{\psi_{\calK}}\longrightarrow E \otimes
    N_{n-2} \calK_\zeta^{\calA} \otimes \Omega^2_{X'/S} \overset{\psi_{\calK}}\longrightarrow \cdots \]
    is a resolution of $E$.
  \item The de Rham complex $(E \otimes \calK_\zeta^{\calA} \otimes \Omega^\cdot_{X/S},
    \nabla_{tot})$ is a resolution of $E'$.
  \end{enumerate}

  \begin{proof}
    Note that the statement for (i) does not involve the connection on
    $E$.  Thus, since each term in the complex $\calK_\zeta^{\calA} \otimes
    \Omega^\cdot_{X'/S}$ is a locally free $\AXgp$-module, and the differential
    is $\AXgp$-linear, the result follows from
    (\ref{lemma:h-dr-res}.\ref{lemma:higgs-res}).  The proof for $E \otimes
    N_{n - \cdot} \calK_\zeta^{\calA} \otimes \Omega^\cdot_{X'/S}$ is similar.
    
    To prove (ii), recall that since $E$ is locally nilpotent, $E \simeq (E'
    \otimes_{\calB_{X/S}} \calK_\zeta^{\calA})^\theta$.  We thus have a map $E
    \otimes_{\AXgp} \calK_\zeta^{\calA} \to E' \otimes_{\calB_{X/S}} \calK_\zeta^{\calA}$
    induced by the multiplication morphism $\calK_\zeta^{\calA} \otimes
    \calK_\zeta^{\calA} \to \calK_\zeta^{\calA}$, which is an isomorphism by a
    result in the appendix to \cite{ogus-vol}.  This isomorphism is
    horizontal since the connection on $\calK_\zeta^{\calA}$ satisfies the
    Leibniz rule and since the connection on $E \simeq (E' \otimes
    \calK_\zeta^{\calA})^\theta$ comes from the action of $\nabla$ on
    $\calK_\zeta^{\calA}$.  Therefore, the de Rham complex $(E \otimes
    \calK_\zeta^{\calA} \otimes \Omega^\cdot_{X/S}, \nabla_{tot})$ comes from tensoring the de
    Rham complex $(\calK_\zeta^{\calA} \otimes \Omega^\cdot_{X/S}, \nabla_{tot})$ with $E'$
    over $\calB_{X/S}$.  Since each term in the latter complex is a
    locally free $\calB_{X/S}$-module (using the fact that $\AXgp$
    itself is), and the differential is $\calB_{X/S}$-linear, the
    desired result follows from
    (\ref{lemma:h-dr-res}.\ref{lemma:derham-res}).
  \end{proof}
\end{lemma}

\begin{proof}[Proof of \ref{thm:cohom-detail}]
  Clearly, (i) follows from (ii) since locally we may choose a
  splitting $\zeta$ of $C_{X/S}$, and then $N_n \calK_{\calXS}^{\calA} \simeq
  N_n \calK_\zeta^{\calA}$ for $n < p$.

  To prove (ii), observe that we have a commutative diagram
  \[
  \begin{diagram}
    \node{\tau_{\leq (n - \ell)} (E \otimes \Omega^\cdot_{X/S}, \nabla)} \arrow{s,r}{g_n}
    \arrow{se,t}{g} \\
    \node{\tau_{\leq (n - \ell)} N_n \calK_\zeta^{\calA \cdot}(E)} \arrow{e,t}{c}
    \node{\tau_{\leq (n - \ell)} \calK_\zeta^{\calA \cdot}(E)} \\
    \node{\tau_{\leq (n - \ell)} (E' \otimes \Omega^\cdot_{X'/S}, \theta').} \arrow{n,r}{f_n}
    \arrow{ne,b}{f}
  \end{diagram}
  \]
  By the previous lemma, the rows of $N_n \calK_\zeta^{\calA \cdot \cdot}(E)$ and of
  $\calK_\zeta^{\calA \cdot \cdot}(E)$ are resolutions of $E \otimes \Omega^i_{X/S}$; therefore,
  $g_n$ and $g$ are quasi-isomorphisms.  Similarly, the columns of
  $\calK_\zeta^{\calA \cdot \cdot}(E)$ are resolutions of $E' \otimes \Omega^i_{X'/S}$, which
  implies that $f$ is a quasi-isomorphism.  From the above diagram, it
  follows that $c$ and $f_n$ are also quasi-isomorphisms.
\end{proof}

If $E \in MIC_\ell(X / S)$ for $\ell < p$, and we apply this result to $E
\otimes_{\scrO_X} \AXgp \in MIC_\ell^{\calA}(X / S)$, then take the degree zero
parts, we get the following.

\begin{corollary}
  Let $(E, \nabla) \in MIC_\ell(X / S)$ for $\ell < p$, and let $E' := (\calK_{\calXS}
  \otimes_{\scrO_X} E)^{\nabla, \gamma}$, with Higgs field $\theta' := \psi_{\calK}$.  (Thus
  $(E', \theta') \simeq C_{\calXS}(E, \nabla)$ in the case that $(E, \nabla) \in MIC_\ell^0(X /
  S)$.)  Let $m := p - \ell - 1$.  Then there is a quasi-isomorphism in
  the derived category of complexes of $\scrO_{X'}$-modules
  \[ \tau_{\leq m} (E \otimes_{\scrO_X} \Omega^\cdot_{X/S}, \nabla) \simeq \tau_{\leq m} (E' \otimes_{\scrO_{X'}}
  \Omega^\cdot_{X'/S}, \theta'). \]
\end{corollary}

\appendix
\section{A Brief Proof of the $F$-linearity of $p$-curvature}

Fix the following notation: let $f : X\to S$ be a morphism of fine log
schemes of characteristic $p$.  Let $Y$ be the logarithmic formal
neighborhood of the diagonal in $X \times_S X$, with the exact closed
immersion $\Delta : X \to Y$ and projections $p_1, p_2 : Y\to X$.  Let
$\regpow_{X/S} = \scrO_Y$, let $J$ denote the ideal of $\Delta$, and let
$\regpow_{X/S}^n = \regpow_{X/S} / J^{n+1}$ denote the structure sheaf
of the $n$th log infinitesimal neighborhood of $\Delta$.  Similarly, let
$\divpow_{X/S}(1)$ denote the structure sheaf of the log divided power
envelope of $\Delta$, with PD ideal $\Jbar$, and $\divpow_{X/S}^n(1) =
\divpow_{X/S}(1) / \Jbar^{[n+1]}$ the structure sheaf of the $n$th log
infinitesimal divided power envelope.  We consider
$\regpow_{X/S}$, $\regpow_{X/S}^n$, $\divpow_{X/S}(1)$, and
$\divpow_{X/S}^n(1)$ as $\scrO_X$-algebras via $p_1^*$, i.e.\ by
multiplication on the left.

For $m$ a section of $\scrM_X^{\gp}$, we define $\eta_m :=
\alpha_Y(p_2^* m - p_1^* m) - 1 \in J$, and $\zeta_m := \log(1 + \eta_m) =
\eta_m - \eta_m^{[2]} + 2! \eta_m^{[3]} - \cdots + (p-1)! \eta_m^{[p]} \in \Jbar$.
(The reader who is unfamiliar with log geometry can skip to Lemma
\ref{lemma:compval}, and replace $\zeta_m$ by $\xi_x := 1 \otimes x - x \otimes 1 \in
\Jbar$ for sections $x$ of $\scrO_X$.)

We begin with a logarithmic generalization of a result which is
well-known in the nonlogarithmic case.

\begin{theorem}
  \label{thm:pthpower}
  Let $D : \regpow^1 \to \scrO_X$ be a \diffop{} of order $\leq 1$.  Then
  $D^p : \regpow^p \to \scrO_X$ is also a \diffop{} of order $\leq 1$.
\end{theorem}

In the case of trivial log structure, we prove this by showing that
$[D^p, a]$ is a \diffop{} of order $\leq 0$ for each $a \in \scrO_X$.  We
will essentially adapt this proof to the logarithmic case by
substituting $a = \alpha(m)$ for $m \in \scrM_X$ and dividing each step by
$\alpha(m)$.  This leads to the following construction:

\begin{definition}
  For $m \in \scrM_X^{\gp}$, and $\phi:\regpow \to \scrO_X$, we define
  $\phi_m : \regpow \to \scrO_X$ by $\phi_m(\tau) = \phi(\tau (1 + \eta_m))$ for
  each $\tau \in \regpow$.
\end{definition}

Note that for $m \in \scrM_X$, since $\alpha(m) (1 + \eta_m) = p_2^*(\alpha(m))$,
we have $\alpha(m) \phi_m = \phi \circ \alpha(m)$.  Thus, we may view $\phi_m$ as
conjugation of $\phi$ by $\alpha(m)$.

\begin{lemma}
  The collection of maps $\Diff(\scrO_X, \scrO_X) \to \Diff(\scrO_X,
  \scrO_X)$ defined by $\phi \mapsto \phi_m$ induces a group
  homomorphism $\scrM_X^{\gp} \to \Aut_{\scrO_X}(\Diff(\scrO_X,
  \scrO_X))$.

  \begin{proof}
    Since $(1 + \eta_m) (1 + \eta_{m'}) = 1 + \eta_{m + m'}$, we see that
    $(\phi_m)_{m'} = \phi_{m + m'}$.  We also have $\eta_0 = 0$, so $\phi_0 =
    \phi$.  Since $\phi \mapsto \phi_m$ is clearly $\scrO_X$-linear, all we have
    left to prove is that $(\phi \circ \psi)_m = \phi_m \circ \psi_m$.
    
    Thus, suppose $\tau \in \regpow$.  By definition, we calculate $(\phi
    \circ \psi)(\tau (1 + \eta_m))$ by first taking $\delta(\tau (1 + \eta_m)) = \delta(\tau)
    \delta(1 + \eta_m)$.  However, $\delta(1 + \eta_m) = (1 + \eta_m) \otimes (1 +
    \eta_m)$, so $\delta(\tau (1 + \eta_m))$ gets mapped by $\id \otimes \psi$ to $(1 +
    \eta_m) \cdot (\id \otimes \psi_m)(\delta(\tau))$.  Applying $\phi$, we get $(\phi_m \circ
    \psi_m)(\tau)$, as desired.
  \end{proof}
\end{lemma}

In the case of a scheme with trivial log structure, we have the result
that a differential operator $\phi : \regpow\to \scrO_X$ is of
order $\leq k$ if and only if $[\phi, a]$ is of order $\leq (k - 1)$ for
each $a \in \scrO_X$.  The following is a logarithmic analogue of this
fact.

\begin{lemma}
  A differential operator $\phi : \regpow\to \scrO_X$ (of finite order)
  is of order $\leq k$ if and only if $\phi_m - \phi$ is of order $\leq (k -
  1)$ for each $m \in \scrM_X^{\gp}$.

  \begin{proof}
    Suppose $\phi$ is of order $\leq \ell$, so that $\phi$ induces a map
    $\regpow^\ell \to \scrO_X$.  Then for $\tau\in \regpow^\ell$, $(\phi_m -
    \phi)(\tau) = \phi(\tau \eta_m)$; thus, $\phi_m - \phi$ is of order $\leq (k - 1)$
    if and only if $\phi(\tau \eta_m) = 0$ whenever $\tau \in J^k$.  However,
    since $\{ \eta_m : m \in \scrM_X \}$ generates $J$ as an
    $\scrO_X$-module in $\regpow^\ell$, this condition is sufficient
    (and obviously necessary) to get $\phi(J^{k+1}) = 0$.
  \end{proof}
\end{lemma}

\begin{proof}[Proof of Theorem \ref{thm:pthpower}]
  First consider the case in which $D(1) = 0$.  By hypothesis together
  with the last lemma, we see that $D_m - D = a$ for some $a \in
  \scrO_X$.  Therefore, $(D^p)_m = (D_m)^p = (D + a)^p$.  However, a
  straightforward induction, using the formula $(D + a) \circ x = (D +
  a)(x) + x D$ for $x \in \scrO_X$, shows that
  \[ (D+a)^k = \sum_{i=0}^k \binom{k}{i} b_i D^{k-i} \]
  for every $k$, where $b_i = (D + a)^i(1)$.  In particular, for $k =
  p$, we get $(D+a)^p = D^p + b_p$, so $(D^p)_m - D^p = b_p$ is a
  \diffop{} of order $\leq 0$.
  
  Now in the general case, the previous paragraph shows that $(D -
  D(1))^p$ is a \diffop{} of order $\leq 1$.  But using the above
  formula again, $D^p = (D - D(1))^p + b_p'$ for some $b_p' \in
  \scrO_X$, so $D^p$ is also a \diffop{} of order $\leq 1$.
\end{proof}

\begin{corollary}
  If $D = (\delta, \partial) : (\scrM_X, \scrO_X) \to \scrO_X$ is a log derivation
  over $S$, then so is $D^{(p)} := (\partial^{p-1} \circ \delta + F_X^* \circ \delta,
  \partial^{(p)})$.

  \begin{proof}
    Corresponding to the given log derivation, we have a \PDop{} $D$
    of order $\leq 1$ with $D(1) = 0$.  Then by Theorem
    \ref{thm:pthpower}, $(D^p)^\flat : \regpow \to \scrO_X$ is a \diffop{}
    of order $\leq 1$, so $(D^p)^\flat : \regpow / J^2 \to \scrO_X$
    restricts to a log derivation $J / J^2 \simeq \Omega^1_{X/S}\to \scrO_X$.
    We claim that this log derivation is exactly $\partial^{p-1} \circ \delta +
    F_X^* \circ \delta$ on $\scrM_X$ and $\partial^{(p)}$ on $\scrO_X$.
    
    The latter statement follows exactly as in the case of trivial log
    structure.  Now by the identification of $J / J^2$ with
    $\Omega^1_{X/S}$, we see the log derivation sends $m \in \scrM_X$ to
    $D^p(\eta_m)$.  However, $\eta_m \equiv \zeta_m + \zeta_m^{[2]} + \cdots +
    \zeta_m^{[p]} \pmod{\Jbar^{[p+1]}}$.  On the other hand, for $1 < k <
    p$, $\zeta_m^{[k]}$ is in the image of the natural map $J^2 \subseteq
    \regpow \to \divpow(1)$, so $D^p(\zeta_m^{[k]}) = 0$.  Therefore,
    $D^p(\eta_m) = D^p(\zeta_m) + D^p(\zeta_m^{[p]})$.
    
    We calculate $D^p(\zeta_m)$ as $(D^{p-1} \circ D)(\zeta_m)$.  We have
    $\delta(\zeta_m) = \zeta_m \otimes 1 + 1 \otimes \zeta_m$, which is mapped by $\id \otimes D$
    to $p_2^*(\delta m)$ since $D(1) = 0$ and $D(\zeta_m) = D(\eta_m) = \delta m$.
    Now applying $D^{p-1}$ gives $\partial^{p-1}(\delta m)$.  The following
    lemma will give that $D^p(\zeta_m^{[p]}) = D(\zeta_m)^p = (\delta m)^p$.
  \end{proof}
\end{corollary}

\begin{lemma}
  \label{lemma:compval}
  Let $\phi, \psi : \divpow(1) \to \scrO_X$ be \PDop{}s of order $\leq k,
  \ell$, respectively.  Then for any section $m$ of $\scrM_X$, we have
  $(\phi \circ \psi)(\zeta_m^{[k + \ell]}) = \phi(\zeta_m^{[k]}) \psi(\zeta_m^{[\ell]})$.

  \begin{proof}
    By definition, we calculate $\phi \circ \psi$ by first taking $\delta^{k,
      \ell}(\zeta_m^{[k+\ell]}) = \sum_{j = 0}^{k + \ell} \zeta_m^{[j]} \otimes \zeta_m^{[k
      + \ell - j]} = \zeta_m^{[k]} \otimes \zeta_m^{[\ell]}$ in $\divpow^k(1) \otimes
    \divpow^\ell(1)$.  We now apply $\id \otimes \psi$, which gives $\zeta_m^{[k]}
    p_2^* \psi(\eta_m^{[\ell]})$.  However, since $\zeta_m^{[k]} \in
    \Jbar^{[k]}$ and $p_2^* \psi(\zeta_m^{[\ell]}) - p_1^* \psi(\zeta_m^{[\ell]})\in
    \Jbar$, we see that $\eta_m^{[k]} p_2^* \psi(\zeta_m^{[\ell]}) =
    \psi(\zeta_m^{[\ell]}) \zeta_m^{[k]}$ in $\divpow^{k+\ell}(1)$.  Thus,
    applying $\phi$ gives $\psi(\zeta_m^{[\ell]}) \phi(\zeta_m^{[k]})$, as desired.
  \end{proof}
\end{lemma}

We will now show the linearity properties of the $p$-curvature by
showing directly that a construction of Mochizuki agrees with the
standard definition of the $p$-curvature.  Mochizuki's construction
begins with the following observation:

\begin{proposition}
  Let $(D(1), \bar I, \gamma)$ denote the logarithmic PD envelope of the
  diagonal in $Y$ (the logarithmic formal neighborhood of the diagonal
  in $X \times_S X$), and $I$ the ideal of the diagonal in $Y$.  Then there
  is a unique isomorphism
  \[ \alpha : F_{X/S}^* \Omega^1_{X'/S} \to \bar I / (\bar I^{[p+1]} + I \scrO_{D(1)})
  \]
  such that for any $\xi \in I$ with image $\omega \in I / I^2 \simeq \Omega^1_{X/S}$, 
  \[ \alpha(1 \otimes \pi^* \omega) = \xi^{[p]}. \]

  \begin{proof}
    First, the corresponding map $\Omega^1_{X/S} \to F_{X *} [\bar I /
    (\bar I^{[p+1]} + I \scrO_{D(1)}) ]$ is additive since
    \[ (\xi + \tau)^{[p]} = \xi^{[p]} + \tau^{[p]} + \sum_{i=1}^{p-1} \frac{1}{i!
      (p-i)!} \xi^{[p-i]} \tau^{[i]}, \]
    where the last term is in $I \scrO_{D(1)}$.  It is also $\scrO_X$-linear
    since
    \[ (a \xi)^{[p]} = a^p \xi^{[p]}. \]
    Finally, for $\xi, \tau \in I$,
    \[ (\xi \tau)^{[p]} = \xi^p \tau^{[p]} \in I \scrO_{D(1)}, \]
    so the map annihilates $I^2$.  We thus get a well-defined map
    \[ F_{X/S}^* \Omega^1_{X'/S} \simeq F_X^* \Omega^1_{X/S} \to \bar I / (\bar
    I^{[p+1]} + I \scrO_{D(1)}). \]

    To see this is an isomorphism, we work locally; thus, assume we
    have a logarithmic system of coordinates $m_1, \ldots, m_r \in \scrM_X$.
    Then $\bar I / (\bar I^{[p+1]} + I \scrO_{D(1)})$ has basis $\eta^{[p
      \epsilon_1]}, \ldots, \eta^{[p \epsilon_r]}$, which is the image under $\alpha$ of the
    basis $1 \otimes \pi^*(\dlog m_1), \ldots, 1 \otimes \pi^*(\dlog m_r)$ of $F_{X/S}^*
    \Omega^1_{X'/S}$.
  \end{proof}
\end{proposition}

Now let $E$ be a crystal of $\scrO_{X/S}$-modules, and let $p_1, p_2 :
D(1) \to X$ be the canonical projections.  Then we get an isomorphism
\[ \epsilon : p_2^* E \overset{\theta_{p_2}}{\longrightarrow} E_{D(1)} \overset{\theta_{p_1}^{-1}}{\longrightarrow}
p_1^* E. \]
Now for $e \in E$, $\epsilon(p_2^* e) - p_1^* e \in E \otimes_{\scrO_X} \bar I$, so this induces a
map
\[ \psi : E \to E \otimes_{\scrO_X} [\bar I / (\bar I^{[p+1]} + I \scrO_{D(1)})]
\simeq E \otimes_{\scrO_X} F_{X/S}^* \Omega^1_{X'/S}. \]

\begin{theorem}
  Let $D \in \scrT_{X/S}$.  Then
  \[ \psi_{\pi^* D} : E \overset{\psi}{\to} E \otimes_{\scrO_X} F_{X/S}^* \Omega^1_{X'/S}
  \overset{\id \otimes F_{X/S}^* \pi^* D}{\longrightarrow} E \]
  is equal to $\nabla_D^p - \nabla_{D^{(p)}}$.
  \begin{proof}
    First, $D^p - D^{(p)} : \scrO_{D(1)} \to \scrO_X$ annihilates $I
    \scrO_{D(1)}$ by the definition of $D^{(p)}$, and it annihilates
    $\bar I^{[p+1]}$ since it is a PD differential operator of order
    $\leq p$.  Therefore, $D^p - D^{(p)}$ induces a well-defined map
    \[ F_{X/S}^* \Omega^1_{X'/S} \simeq \bar I / (\bar I^{[p+1]} + I
    \scrO_{D(1)}) \to \scrO_X. \]
    We claim that in fact, this map agrees with $\id \otimes F_{X/S}^* \pi^*
    D$.  It suffices to check this for $1 \otimes \pi^* \omega$ for $\omega \in
    \Omega^1_{X/S}$.  However, taking a preimage $\xi \in I$ of $\omega$, we get
    $\alpha(1 \otimes \pi^* \omega) = \xi^{[p]}$, which is mapped by $D^p$ to $(D \omega)^p =
    F_X^* (D \omega)$ and by $D^{(p)}$ to zero.  On the other hand, $\id \otimes
    F_{X/S}^* \pi^* D$ maps $1 \otimes \pi^* \omega$ to $F_X^* (D \omega)$ also.

    Therefore, we may rewrite $\psi_{\pi^* D}$ as
    \[ E \overset{\bar \psi}{\to} E \otimes_{\scrO_X} \bar I \to E \otimes_{\scrO_X}
    (\bar I / (\bar I^{[p+1]} + I \scrO_{D(1)})) \overset{\id \otimes (D^p -
      D^{(p)})}{\longrightarrow} E, \]
    where $\bar \psi : E \to E \otimes_{\scrO_X} \bar I$ maps $e$ to $\epsilon(p_2^* e)
    - p_1^* e$.  However, this is exactly $\nabla_{D^p - D^{(p)}}$.
  \end{proof}
\end{theorem}

The immediate corollary of this is:

\begin{corollary}
  If $\nabla : E \to E \otimes_{\scrO_X} \Omega^1_{X/S}$ is an integrable connection,
  then the $p$-curvature map $\psi(\nabla) : \scrT_{X/S} \to F_{X*}
  \scrEnd_{\scrO_X}(E)$, $D \mapsto \nabla_D^p - \nabla_{D^{(p)}}$, is $F$-linear.
\end{corollary}

A more concrete proof of the linearity starts with the following
criterion:

\begin{lemma}
  Let $\phi : \divpow^p(1) \to \scrO_X$ be a \PDop{} of order $\leq p$.
  Then $\phi$ is of order $\leq 1$ if and only if $\phi^\flat : \regpow^p \to
  \scrO_X$ is a \diffop{} of order $\leq 1$, and $\phi(\zeta_m^{[p]}) = 0$
  for every $m \in \scrM_X$.

  \begin{proof}
    The forward implication is trivial.  For the reverse implication,
    since the statement is local, we may choose a logarithmic system
    of coordinates $m_1, \ldots, m_r \in \scrM_X$, and let $\zeta^{[k]}$ ($k
    \in \nats^r$) be the corresponding basis for $\divpow(1)$.  Let
    $\epsilon_1, \ldots, \epsilon_r$ denote the canonical basis of $\nats^r$.  Since
    $\phi$ is of order $\leq p$, $\phi(\zeta^{[k]}) = 0$ for $|k| > p$.  If $2
    \leq |k| \leq p$ but $k \neq p \epsilon_i$ for any $i$, then in fact
    $\zeta^{[k]}$ is in the image of the natural map $J^2 \subseteq \regpow \to
    \divpow(1)$.  Since $\phi^\flat$ is of order $\leq 1$, this implies
    $\phi(\zeta^{[k]}) = 0$.  Finally, if $k = p \epsilon_i$, then $\phi(\zeta^{[k]})
    = \phi(\zeta_{m_i}^{[p]}) = 0$.
  \end{proof}
\end{lemma}

Using this criterion along with (\ref{lemma:compval}), it is easy to
show that both $(D_1 + D_2)^p - D_1^p - D_2^p$ and $(aD)^p - a^p D^p$
are PD differential operators of order $\leq 1$.  Therefore, for
instance, $(aD)^p - a^p D^p = (aD)^{(p)} - a^p D^{(p)}$, which implies
$\psi_{\pi^*(aD)} = a^p \psi_{\pi^* D}$, and similarly for the additivity.

\bibliographystyle{amsalpha}
\bibliography{loghodge}

\providecommand{\bysame}{\leavevmode\hbox to3em{\hrulefill}\thinspace}
\providecommand{\MR}{\relax\ifhmode\unskip\space\fi MR }
\providecommand{\MRhref}[2]{%
  \href{http://www.ams.org/mathscinet-getitem?mr=#1}{#2}
}
\providecommand{\href}[2]{#2}
\begin{thebibliography}{Mon02}

\bibitem[Kat70]{katz}
Nicholas Katz, \emph{Nilpotent connections and the monodromy theorem:
  Applications of a result of {T}urrittin}, Inst. Hautes {\'E}tudes Sci. Publ.
  Math. \textbf{39} (1970), 175--232.

\bibitem[Kat88]{kato}
Kazuya Kato, \emph{Logarithmic structures of {F}ontaine-{I}llusie}, Algebraic
  Analysis, Geometry, and Number Theory (Jun-Ichi Igusa, ed.), The Johns
  Hopkins University Press, 1988, pp.~191--224.

\bibitem[Lor00]{lorenzon}
Pierre Lorenzon, \emph{Indexed algebras associated to a log structure and a
  theorem of $p$-descent on log schemes}, Manuscripta Mathematica \textbf{101}
  (2000), 271--299.

\bibitem[Mon02]{montagnon}
Claude Montagnon, \emph{G{\'e}n{\'e}ralisation de la th{\'e}orie
  arithm{\'e}tique des $d$-modules {\`a} la g{\'e}om{\'e}trie logarithmique},
  Ph.D. thesis, L'universit{\'e} de Rennes I, 2002.

\bibitem[Ogu04]{ogus-higgs}
Arthur Ogus, \emph{Higgs cohomology, $p$-curvature, and the {C}artier
  isomorphism}, Compositio Mathematica \textbf{140} (2004), 145--164.

\bibitem[OV]{ogus-vol}
A.~Ogus and V.~Vologodsky, \emph{Nonabelian {H}odge theory in characteristic
  $p$}, in preparation.

\bibitem[Sim92]{simpson}
Carlos Simpson, \emph{Higgs bundles and local systems}, Inst. Hautes {\'E}tudes
  Sci. Pub. Math. \textbf{75} (1992), 5--95.

\bibitem[Tsu96]{tsuji}
T.~Tsuji, \emph{Syntomic complexes and $p$-adic vanishing cycles}, J. Reine
  Angew. Math. \textbf{472} (1996), 69--138.

\end{thebibliography}
\end{document}